\documentclass[12pt,fullpage]{article}

\usepackage{amsfonts}
\usepackage{xr}
\usepackage{graphicx}
\usepackage{amsmath}
\usepackage{epsfig}

\def\Z{\mathbb{Z}}
\def\R{\mathbb{R}}
\def\N{\mathbb{N}}
\def\C{\mathbb{C}}
\def\S{\mathcal{S}}
\def\c{\mathcal{C}}

\def\D{\kappa}
\def\epsilon{\varepsilon}
\def\trait (#1) (#2) (#3){\vrule width #1pt height #2pt depth #3pt}
\def\fin{\hfill\trait (0.1) (5) (0) \trait (5) (0.1) (0) \kern-5pt 
\trait (5) (5) (-4.9) \trait (0.1) (5) (0)}

\newcommand{\SE}{\setcounter{equation}{0} \section}
\newcommand{\be}{\begin{equation}}
\newcommand{\ee}{\end{equation}}
\newcommand{\baa}{\begin{array}}
\newcommand{\eaa}{\end{array}}
\newcommand{\ba}{\begin{eqnarray}}
\newcommand{\ea}{\end{eqnarray}}

\newtheorem{theo}{\bf Theorem}[section]
\newtheorem{lem}[theo]{\bf Lemma}
\newtheorem{pro}[theo]{\bf Proposition}
\newtheorem{cor}[theo]{\bf Corollary}
\newtheorem{defi}[theo]{\bf Definition}
\newtheorem{rem}[theo]{\bf Remark}

\newtheorem{fact}[theo]{\bf Fact}

\begin{document}
\date{}
\title{\bf{Hardy spaces of differential forms on Riemannian 
manifolds}}
\author{Pascal Auscher$^{\hbox{\small{a}}}$, Alan 
McIntosh$^{\hbox{\small{b}}}$ and Emmanuel Russ$^{\hbox{\small{c}}}$\\
\\
\footnotesize{$^{\hbox{a }}$Universit\'e de Paris-Sud, Orsay et CNRS 
UMR 8628, 91405 Orsay Cedex, France}\\
\footnotesize{$^{\hbox{b }}$Centre for Mathematics and its 
Applications, Mathematical Sciences Institute}\\
\footnotesize{Australian National University, Canberra ACT 0200, 
Australia}\\
\footnotesize{$^{\hbox{c }}$Universit\'e Paul C\'ezanne, LATP, CNRS 
UMR 6632,
Facult\'e des Sciences et Techniques, Case cour A}\\
\footnotesize{Avenue Escadrille Normandie-Ni\'emen, F-13397 Marseille 
Cedex 20, France}\\
\footnotesize{Pascal.Auscher@math.u-psud.fr, 
Alan.McIntosh@maths.anu.edu.au, emmanuel.russ@univ-cezanne.fr}}

\maketitle

\noindent{\small{{\bf Abstract.} }} Let $M$ be a complete connected
Riemannian 
manifold. Assuming that the Riemannian measure is doubling, we define 
Hardy 
spaces $H^p$ of differential forms on $M$ and give various 
characterizations of them, 
including an atomic decomposition. As a consequence, we derive the 
$H^p$-boundedness for Riesz transforms on $M$, generalizing 
previously known results. 
Further applications, in particular to $H^{\infty}$ functional 
calculus 
and Hodge decomposition, are given.

\noindent{\small{{\bf AMS numbers 2000: }}} Primary: 42B30. Secondary:  58J05, 47F05, 47A60.


\tableofcontents
%%%%%%%%%%%%%%%%%%%%%%%%%%%%%%%%%%%%%%%%%%%%%%%%%%%%%%%%%%%%%%%%%%%%%%%%%%%%%%%

%%%%%%%%%%%%%%%%%%%%%%%%%%%%%%%%%%%%%%%%%%%%%%%%%%%%%%%%%%%%%%%%%%%%%%%%%%%%%%%

\SE{Introduction and main results}\label{intro}

The study of Hardy spaces started in the 1910's and was 
closely related to Fourier 
series and complex analysis in one variable (see \cite{zygmund}, 
Chapters 7 and 14). In the 1960's, an essential feature of the 
development 
of real analysis in several variables was the theory of real Hardy 
spaces 
$H^p(\R^n)$, and in particular $H^1(\R^n)$, which began with the 
paper of Stein and Weiss \cite{steinweiss}. In this work, Hardy 
spaces were defined and studied 
by means of Riesz transforms and harmonic functions. The celebrated 
paper of Fefferman and Stein \cite{feffermanstein} provided many 
characterizations of Hardy 
spaces on $\R^n$, in particular in terms of suitable maximal 
functions. The dual space of $H^1(\R^n)$ was also identified as 
$BMO(\R^n)$. An important step was 
the atomic decomposition of $H^1(\R^n)$, due to Coifman (for $n=1$, 
\cite{coifman}) and to Latter (for $n\geq 2$, \cite{latter}).  A detailed review and 
bibliography on these topics may 
be found in \cite{semmes}. Hardy spaces have been generalized to various geometric settings. See for example the work of Strichartz \cite{strihardy} for compact manifolds with a characterization via pseudo-differential operators, and more generally, starting from the point of view of the atomic decomposition, the work of Coifman and Weiss \cite{coifmanweiss} for spaces of homogeneous type,  which are known to be a relevant 
setting for most tools in harmonic analysis such as Hardy-Littlewood 
maximal function, 
covering lemmata, Whitney decomposition, Calder\'on-Zygmund 
decomposition and singular integrals (see \cite{coifmanweiss}, 
\cite{coifmanweissspringer}, 
\cite{davidsemmes}).

The connection between Hardy spaces and area functionals will be most important to us thanks to the theory of tent spaces developed by Coifman, Meyer and Stein in \cite{coifmanmeyerstein}.  Let us recall the main line of ideas. For suitable functions $f$ on $\R^n$ and all $x\in \R^n$, define
\[
\S f(x)=\left(\iint_{\left\vert y-x\right\vert<t} \left\vert 
t\sqrt{\Delta}e^{-t\sqrt{\Delta}}f(y)\right\vert^2 \, \frac{dydt}{t^{n+1}}\right)^{1/2}
\]
where $\Delta = - \sum_{j=1}^n \frac{\partial^2}{\partial x_j^2}$. Hence, $e^{-t\sqrt{\Delta}}$ is nothing but the Poisson semigroup.
It is proved in \cite{feffermanstein}
that, if $f\in H^1(\R^n)$, then $\S f\in L^1(\R^n)$ and
\[
\left\Vert \S f\right\Vert_{L^1(\R^n)}\leq C\left\Vert 
f\right\Vert_{H^1(\R^n)}.
\]
This fact exactly means that for $f\in H^1(\R^n)$ the function $F$ defined by
$F(t,x)=t\sqrt{\Delta}e^{-t\sqrt{\Delta}}f(x)$ belongs to the tent space 
$T^{1,2}(\R^n)$. Conversely, for  any $F\in T^{1,2}(\R^n)$ and $F_t(x)=F(t,x)$, the function $f$ given
by
\begin{equation} \label{projection}
f=\int_0^{+\infty} t\sqrt{\Delta}e^{-t\sqrt{\Delta}} F_t\, \frac{dt}t
\end{equation}
is in the Hardy space $H^1(\R^n)$ with appropriate estimate. The round trip is granted by
the Calder\'on reproducing formula
\[
f=4\int_0^{+\infty} t\sqrt{\Delta}e^{-t\sqrt{\Delta}} 
t\sqrt{\Delta}e^{-t\sqrt{\Delta}} f\, \frac{dt}t.
\]
It is proved in \cite{coifmanmeyerstein} that the tent space $T^{1,2}(\R^n)$ has an atomic decomposition. Applying this to the function $F(t,x)=t\sqrt{\Delta}e^{-t\sqrt{\Delta}}f(x)$ and plugging this into the reproducing formula, one obtains a decomposition of $f$, not  as a sum of atoms but of so-called molecules (see \cite{coifmanweiss}). Molecules do not have compact support but decay sufficiently fast that they can be used in place of atoms  for many purposes. 
Had we changed the operator $t\sqrt{\Delta}e^{-t\sqrt{\Delta}}$ to an appropriate convolution operator with compactly supported kernel, then  the same strategy would give an atomic decomposition of $f$.  
 The tent space method can be used in different contexts, for instance to obtain 
an atomic decomposition for Hardy spaces defined by maximal functions 
involving second order elliptic operators, see \cite{grenoble}. See 
also \cite{wilson} for a variant of this argument. Let us also 
mention that the duality for tent spaces provides an alternative 
proof of the $H^1-BMO$ duality (see \cite{feffermanstein}, \cite{stein}). 

\bigskip

If one wants to replace functions by forms given some differential structure, then the first thing 
that changes is the mean value condition. A function $f$ in the Hardy space $H^1(\R^n)$ has vanishing mean, that is $\int f=0$. For general forms, the integral has no meaning. 
However, an appropriate atomic decomposition of the Hardy space of divergence free vector fields was proved in \cite{ghl}. This space turned out to be a  specific case of the Hardy spaces of 
exact differential forms 
on $\R^n$ defined  by Lou and the second author in 
\cite{loumcintosh} via the tent spaces approach. There, atomic decompositions, duality results, among other things, were obtained.  What replaces the mean value property in \cite{loumcintosh}  is the fact that atoms are \textit{exact} forms (see section \ref{mol}).

One motivation for studying Hardy spaces of forms is the Riesz transforms. Indeed, the Riesz transforms $R_j= \partial_{x_j} (\sqrt \Delta) ^{\, -1}$  are well-known bounded operators on $H^1(\R^n)$ (and $L^p(\R^n)$, $1<p<\infty$). However, the vector  map $(R_1,R_2,..., R_n)= \nabla (\sqrt \Delta) ^{\, -1} $ is geometrically meaningful as its target space is a space of gradient vector fields (in a generalized sense). This observation is valid in any Riemannian manifold.  On such manifolds 
the understanding of the $L^p$ boundedness property of  Riesz transforms was proposed by Strichartz in \cite{stri}. What happens at  $p=1$  is interesting in itself and also part of this quest as  it can give  results for $p>1$ by interpolation.  As the ``geometric'' Riesz transform is form-valued, getting satisfactory $H^{1}$ boundedness statements for 
this operator requires the notion of Hardy spaces of differential 
forms. Our aim is, therefore, to develop an appropriate theory   of Hardy spaces on Riemannian manifolds (whether or not compact) and to apply this to the Riesz transform (for us, only the ``geometric''  Riesz transform matters so that we drop the ``s'' in transform). This will
indeed generalize the theory 
in \cite{loumcintosh} to a geometric context. In particular, their Hardy spaces will be 
our  spaces $H^1_d(\Lambda^kT^{\ast}\R^n)$ for $0\leq k\leq n$ (see Section \ref{molecgauss}). Also, our Hardy spaces are designed so that the Riesz transform is automatically  bounded on them. Specializing to specific situations allows us to recover results obtained by the third 
author alone \cite{russ} or with M. Marias \cite{mariasruss}. \par

\bigskip

We now describe precisely our setting. Let $M$ be a complete Riemannian manifold, $\rho$ the geodesic distance and 
$d\mu$ the Riemannian measure. Complete means that any two points can be joined by a geodesic, thus $M$ is connected.
For all $x\in M$ and all $r>0$, $B(x,r)$ stands for the open geodesic 
ball with center $x$ and radius $r$, and its measure will be denoted 
$V(x,r)$. 

For all $x\in M$, denote by $\Lambda T_x^{\ast}M$ the complex 
exterior algebra over the cotangent space $T_x^{\ast}M$. Let 
$\Lambda T^{\ast}M=\oplus_{0\leq k\leq \mbox{dim }M} 
\Lambda^k T^{\ast}M$ be the bundle over 
$M$ whose fibre at each $x\in M$ is given by $\Lambda T^{\ast}_xM$, and 
let $L^2(\Lambda T^{\ast}M)$ be the space of square integrable sections of $\Lambda T^{\ast}M$. Denote 
by $d$ the exterior differentiation. Recall that, for 
$0\leq k\leq \mbox{dim }M-1$, $d$ maps, for instance, 
$C^{\infty}_0(\Lambda^k T^{\ast}M)$ into $C^{\infty}_0(\Lambda^{k+1} 
T^{\ast}M)$ and that $d^2=0$. Denote also by $d^{\ast}$ the adjoint 
of $d$ on 
$L^2(\Lambda T^{\ast}M)$. Let $D=d+d^{\ast}$ be the Hodge-Dirac operator on 
$L^2(\Lambda T^{\ast}M)$, and 
$\Delta=D^2=dd^{\ast}+d^{\ast}d$ the (Hodge-de Rham) Laplacian. The 
$L^2$ Hodge decomposition, valid on any complete Riemannian manifold, 
states that
\[
L^2(\Lambda T^{\ast}M)=\overline{{\mathcal R}(d)}\oplus \overline{{\mathcal 
R}(d^{\ast})}\oplus 
{\mathcal N}(\Delta),
\]
where ${\mathcal R}(T)$ (resp. ${\mathcal N}(T)$) stands for the 
range (resp. the nullspace) of $T$, and the decomposition is 
orthogonal. See for instance \cite{derham}, Theorem 24, p. 165, and 
\cite{carronsurvey}.

In view of the previous discussions, we will start the approach of Hardy spaces via tent spaces. The first observation is that  this theory  can be developed in spaces of homogeneous type subject to an additional technical condition \cite{russtent}. For us, it only means that we impose on $M$ the doubling property:
 there exists $C>0$ 
such that, for all $x\in M$ and all $r>0$,
\begin{equation} \label{D}
V(x,2r)\leq CV(x,r).
\end{equation}
A straightforward consequence of (\ref{D}) is that there exist 
$C,\D>0$ such that, for all $x\in M$, all $r>0$ and all $\theta>1$,
\begin{equation} \label{Dbis}
V(x,\theta r)\leq C\theta^\D V(x,r).
\end{equation}
The hypothesis (\ref{D}) exactly means that $M$, equipped with its 
geodesic distance and its Riemannian measure, is a space of 
homogeneous type in the sense of 
Coifman and Weiss. 

There is a wide class of manifolds on which (\ref{D}) holds. First, 
it is true on Lie groups with polynomial volume growth (in particular 
on nilpotent Lie groups), and in this 
context the heat kernel on functions does satisfy Gaussian estimates, 
see \cite{salofflie}. In particular, (\ref{D}) is true if $M$ has 
nonnegative Ricci curvature thanks to the Bishop comparison theorem 
(see \cite{BC}). Recall also that (\ref{D}) remains valid if $M$
is quasi-isometric to a manifold with nonnegative Ricci curvature, or
is a cocompact covering manifold  whose deck transformation group has
polynomial
growth, \cite{CSC}. Contrary to the doubling property, the 
nonnegativity of the Ricci curvature is not stable under 
quasi-isometry.

The second observation is that the Euclidean  proofs using tent spaces 
and formul\ae\  such as \eqref{projection} use pointwise bounds on kernels of the Poisson semigroup (or of appropriate convolution operators).  Here, the only available operators are functions of $\Delta$ and this would require 
some knowledge, say,  on the kernel $p_t(x,y)$ of the heat semigroup. If one deals with the Laplace-Beltrami operator on 
functions, ``Gaussian'' pointwise estimates  may hold for the heat kernel $p_t$, but this depends 
 on further geometric assumptions on $M$. For instance, when  $M$ 
 non-compact, it is 
well-known (see \cite{saloff-coste}, \cite{grigoryan}) that  $M$ satisfies the doubling property  and a scaled $L^2$ Poincar\'e 
inequality on balls if and only if $p_t$ satisfies a ``Gaussian'' upper and 
lower estimate and 
is H\"older continuous. 
More precisely, there exist 
$C_1,c_1,C_2,c_2,\alpha>0$ such that, for all 
$x,x^{\prime},y,y^{\prime}\in M$ and all $t>0$,
\begin{equation} \label{G}
\begin{array}{l}
\displaystyle \frac 
{c_2}{V(x,\sqrt{t})}e^{-C_2\frac{\rho^2(x,y)}{t}}\leq p_t(x,y)\leq 
\frac {C_1}{V(x,\sqrt{t})}e^{-c_1\frac{\rho^2(x,y)}{t}},\\
\displaystyle \left\vert p_t(x,y)-p_t(x^{\prime},y)\right\vert \leq 
C\left(\frac{\rho(x,x^{\prime})}{\sqrt{t}}\right)^{\alpha},\ 
\left\vert p_t(x,y)-p_t(x,y^{\prime})\right\vert \leq 
C\left(\frac{\rho(y,y^{\prime})}{\sqrt{t}}\right)^{\alpha}.
\end{array}
\end{equation}
Note that such  a result concerns  the heat kernel on functions, {\it i.e.} on $0$-forms. For the heat kernel on 1-forms, the pointwise Gaussian domination holds for $|p_t(x,y)|$ if $M$ has non negative Ricci curvature from the  Weitzenb\"ock formula (see \cite{bakry} and also the recent work 
\cite{coulhonzhang} for more and the references therein).  
Very little seems to be known about estimates for the heat 
kernel on general forms.

Hence, for our theory to be applicable, we have to forbid  the use of Gaussian estimates 
similar to $(\ref{G})$. Fortunately, there is a weaker notion of Gaussian decay, 
which holds on any complete  Riemannian manifold, namely the notion of $L^2$ off-diagonal estimates, as introduced by Gaffney \cite{gaffney}. This notion has already proved to be a good substitute of Gaussian estimates for such questions as the Kato square root 
problem or $L^p$-bounds for Riesz transforms when dealing with elliptic operators (even in the Euclidean setting) for which Gaussian estimates do not hold 
(see \cite{auscher, ahlmt, akm} in the Euclidean setting, and 
\cite{acdh} in a complete Riemannian manifold). We show in the present work that a theory of Hardy spaces of differential forms can be developed under such a notion. 

\bigskip

The results of this work have been announced in the Note \cite{amr}. Let us state the main ones. We define in fact three classes of Hardy spaces of differential forms on 
manifolds satisfying (\ref{D}). These definitions require some preliminary material, and we remain  vague at this stage. The first class, 
denoted by $H^1(\Lambda T^{\ast}M)$,  is the one defined via tent spaces. Actually, using fully the theory of tent spaces, we also define 
$H^p(\Lambda T^{\ast}M)$ 
for all $1\leq p\leq +\infty$. The second class, 
$H^1_{mol}(\Lambda T^{\ast}M)$, is defined via  ``molecules'' (see above). Our third class, 
$H^1_{max}(\Lambda T^{\ast}M)$, is defined in terms of an appropriate maximal 
function 
associated to the Hodge-de Rham Laplacian. Within each class, the Hardy spaces are Banach spaces with norms depending on some parameters. We show they are identical spaces  with equivalence of norms. Eventually we prove that the three classes are the same. This can be summarized as follows: 

\begin{theo} \label{equality}
Assume (\ref{D}). Then, 
$H^1(\Lambda T^{\ast}M)=H^1_{mol}(\Lambda T^{\ast}M)=H^1_{max}(\Lambda T^{\ast}M)$.
\end{theo}

Let us mention that we do not use much of the differential structure to prove the first
equality.  As a matter of fact, it can be proved on a space of homogeneous type for  an operator satisfying $L^2$ off-diagonal bounds and $L^2$ quadratic estimates. We leave this point to further works (see also Remark 1.4 below).

As a corollary of Theorem \ref{equality}, we derive the following comparison between $H^p(\Lambda T^{\ast}M)$ and $L^p(\Lambda T^{\ast}M)$:
\begin{cor} \label{lp} Assume (\ref{D}).
\begin{itemize}
\item[$(a)$]
For all $1\leq p\leq 2$, $H^p(\Lambda T^{\ast}M)\subset L^p(\Lambda T^{\ast}M)$, and more precisely, $H^p(\Lambda T^{\ast}M)\subset\overline{{\mathcal R}(D)\cap L^p(\Lambda T^{\ast}M)}^{L^p(\Lambda T^{\ast}M)}$.
\item[$(b)$]
For $2\leq p<+\infty$, $\overline{{\mathcal R}(D)\cap L^p(\Lambda T^{\ast}M)}^{L^p(\Lambda T^{\ast}M)}\subset H^p(\Lambda T^{\ast}M)$.
\end{itemize}
\end{cor}

Of course, it may or may not be that equalities hold for some/all $p \in (1,\infty)\setminus \{2\}$. 

For our motivating operator, namely the Riesz transform   $D\Delta^{-1/2}$ on $M$,  we obtain a satisfactory answer. 

\begin{cor} \label{riesz}
Assume (\ref{D}). Then, for all $1\leq p\leq +\infty$, 
$D\Delta^{-1/2}$ is $H^p(\Lambda T^{\ast}M)$ 
bounded. Consequently, it is $H^1(\Lambda T^{\ast}M)-L^1(\Lambda T^{\ast}M)$ bounded.
\end{cor}

The plan of the paper is as follows. As a preliminary section (Section \ref{h2}), we focus on the case of $H^2(\Lambda T^{\ast}M)$ and define what we mean by Riesz transform, because this case just requires well-known facts of the Hodge-de Rham theory of $L^2(\Lambda T^{\ast}M)$, and this motivates the foregoing technical tools needed to define and study $H^p(\Lambda T^{\ast}M)$ spaces for $p\neq 2$. Section \ref{OD(N)} is devoted to 
the statement and the proof of the off-diagonal $L^2$ estimates for 
the Hodge-Dirac operator and the Hodge-de Rham Laplacian. In Section \ref{tent}, we present tent spaces on $M$ 
and 
establish the boundedness of some ``projectors'' on these spaces. 
Relying 
on this fact, we define Hardy spaces $H^p(\Lambda T^{\ast}M)$ for $1\leq p\leq 
+\infty$ 
in Section \ref{def} and state duality and interpolation 
results. We also establish the $H^p(\Lambda T^{\ast}M)$ boundedness of Riesz 
transforms 
(Corollary \ref{riesz}) and show more generally that there is a functional calculus on $H^p(\Lambda T^{\ast}M)$. 
Section \ref{mol} is devoted to the description of molecules and the 
identification of $H^1(\Lambda T^{\ast}M)$ with $H^1_{mol}(\Lambda T^{\ast}M)$. As a consequence, 
we obtain Corollary \ref{lp}, 
 which completes the proof of Corollary \ref{riesz}.
In Section \ref{max}, we prove the maximal 
characterization of 
$H^1(\Lambda T^{\ast}M)$, which ends the proof of Theorem \ref{equality}.

In Section \ref{appl}, we give 
further examples and applications of the previous results. Namely, 
specializing to the case of $0$-forms, we compare our Hardy spaces 
with the classical Hardy spaces for functions under suitable 
assumptions on $M$ (such as Poincar\'e inequalities), 
generalizing known results about the Riesz transform. 
 We also go further in the comparison of $H^p(\Lambda T^{\ast}M)$ with $L^p(\Lambda T^{\ast}M)$  
assuming  ``Gaussian'' estimates for the 
heat kernel, and recover well-known results about the $L^p$ boundedness for the Riesz transform. 

\medskip

\begin{rem}
During the preparation of this manuscript, we learnt that S. Hofmann and S. Mayboroda have been developing the theory of Hardy spaces associated with second order elliptic operators in divergence form in $\R^n$  \cite{HM}.  This is an alternative generalisation of the usual theory, which is associated with the Laplacian on $\R^N$. Although there is much in common, such as the use of off-diagonal estimates, the results are different, and the proofs have been obtained independently.
\end{rem}

{\bf Acknowledgments: }
This work was supported by the bilateral CNRS/ARC linkage (2003-2004): Espaces de Hardy de formes diff\'erentielles / Hardy spaces of differential forms. 
The research was conducted at the Centre for Mathematics and its Applications at the Australian National University, Canberra, and at the Laboratoire d'Analyse, 
Topologie et Probabilit\'es at the Facult\'e des Sciences et Techniques de Saint-J\'er\^ome, Universit\'e Paul C\'ezanne, Marseille. The authors would also 
like to thank Andreas Axelsson, Thierry Coulhon, Andrew Hassell, 
Zengjian Lou and Philippe Tchamitchian for interesting and helpful 
conversations.\par 
This work was presented by the second author at the Conference on 
Harmonic Analysis and Related Problems, 19-23 June 2006, Zaros 
(Crete), at the conference on ``Recent Advances in Nonlinear 
Partial Differential Equations, a Celebration of Norman Dancer's 60th 
Birthday'', 16-21 July 2006, University of New England (Australia) and 
at the 50th Annual Meeting of the Australian Mathematical Society, 
25-29 September 2006, Macquarie University (Australia).

{\bf Notation: }
If two quantities $A(f),B(f)$ depend on a function $f$ ranging over a 
certain
space $L$, $A(f)\sim B(f), \mbox{ for }f\in L,$
means that there exist $c,C>0$ such that $c\,A(f)\le B(f) \le C\,A(f), \ \forall\,f\in L.$

\SE{The $H^2(\Lambda T^{\ast}M)$ space and the Riesz transform} \label{h2}
Set $H^2(\Lambda T^{\ast}M)=\overline{{\mathcal R}(D)}= \overline{\{Du\in L^2(\Lambda T^{\ast}M);u\in L^2(\Lambda T^{\ast}M)\}}$ and note that
\[
L^2(\Lambda T^{\ast}M)=\overline{{\mathcal R}(D)}\oplus {\mathcal N}(D)=H^2(\Lambda T^{\ast}M)\oplus {\mathcal N}(D) .
\]
It is an essential fact for the sequel that $H^2(\Lambda T^{\ast}M)$ can be described in terms of tent spaces and appropriate quadratic functionals, which we describe now. If $\displaystyle \theta\in \left(0,\tfrac{\pi}2\right)$, set
\[
\begin{array}{lll}
\displaystyle \Sigma_{\theta^{+}} &= & \displaystyle \left\{z\in 
\C\setminus \left\{0\right\};\ \left\vert \mbox{arg }z\right\vert\leq 
\theta\right\} \cup \left\{0\right\},\\
\displaystyle \Sigma^{0}_{\theta^{+}} & = & \displaystyle \left\{z\in 
\C\setminus\left\{0\right\};\ \left\vert \mbox{arg }z\right\vert< \theta\right\},\\
\displaystyle \Sigma_{\theta} & = & \displaystyle 
\Sigma_{\theta^{+}}\cup \left(-\Sigma_{\theta^{+}}\right),\\
\displaystyle \Sigma_{\theta}^0 & = & \displaystyle 
\Sigma_{\theta^{+}}^0\cup \left(-\Sigma_{\theta^{+}}^0\right).
\end{array}
\]
and denote by 
$H^{\infty}(\Sigma_{\theta}^0)$ the algebra of bounded holomorphic 
functions on $\Sigma_{\theta}^0$.  Given $\sigma, \tau>0$, define
   $\Psi_{\sigma,\tau}(\Sigma_{\theta}^0)$ to be
the set of holomorphic functions $\psi\in H^{\infty}(\Sigma_{\theta}^0)$ which satisfy
\[
\left\vert \psi(z)\right\vert \leq C\inf\{\left\vert z\right\vert^\sigma, 
\left\vert z\right\vert^{-\tau}\}
\]
for some $C>0$ and all $z\in \Sigma_{\theta}^0$. Then let
$\Psi(\Sigma_{\theta}^0)=\cup_{\sigma,\tau>0}\Psi_{\sigma,\tau}(\Sigma_{\theta}^0)$. 

For example, if $\psi(z)= z^N(1\pm iz)^{-\alpha}$ for  integers $N,\alpha$ with $1\leq N< \alpha$ then $\psi\in\Psi_{N,\alpha - N}(\Sigma_{\theta}^0)$, if $\psi(z) = z^N(1+z^2)^{-\beta}$ for  integers $N,\beta$ with $1\leq N< 2\beta$ then $\psi\in\Psi_{N,2\beta - N}(\Sigma_{\theta}^0)$, and if $\psi(z)= z^N \exp(-z^2)$ for a non-negative integer $N$ then
$\psi\in \Psi_{N,\tau}(\Sigma_{\theta}^0)$ for all $\tau>0$.

Define $\displaystyle 
{\mathcal H}=L^2\left((0,+\infty),L^2(\Lambda T^{\ast}M),\frac{dt}t\right)$, equipped with 
the norm 
\[
\left\Vert F\right\Vert_{\mathcal H}=\left(\int_0^{+\infty}\int_M \left\vert 
F(x,t)\right\vert^2dx\frac{dt}t\right)^{1/2},
\]
with $\left\vert F(x,t)\right\vert^{2}=\langle 
F(x,t),F(x,t)\rangle_{x}$, where $\langle .,.\rangle_{x}$ stands for 
the inner complex product in $T^{\ast}_{x}M$, and we drop the 
subscript $x$ in the notation to simplify the exposition. Note also that, here and after, we write $dx,dy,\ldots$ instead of 
$d\mu(x),d\mu(y),\ldots$. If $F\in {\mathcal H}$ and $t>0$, denote by $F_t$ the map $x\mapsto 
F(x,t)$.

Given $\psi\in \Psi(\Sigma_{\theta}^0)$ for some $\theta>0$, set 
$\psi_{t}(z)=\psi(tz)$ for all $t>0$ and all $z\in 
\Sigma_{\theta}^{0}$ and define the operator $\mathcal Q_\psi : L^2(\Lambda T^{\ast}M) \to \mathcal H$ by  
\[
\left({\mathcal Q}_{\psi}h\right)_{t}=\psi_t(D)h\ ,\quad t>0\ .
\]
Since $D$ is a self-adjoint operator on $L^2(\Lambda T^{\ast}M)$, it follows from the spectral theorem that
$\mathcal Q_\psi$ is bounded, and indeed that 
\[
\left\Vert {\mathcal Q}_{\psi}f\right\Vert_{\mathcal H}\sim
\left\Vert f\right\Vert_2
\]
for all $f\in H^2(\Lambda T^{\ast}M)$. (Note that ${\mathcal Q}_{\psi}f=0$ for all $f\in\mathcal N(D)$.)
Also define the operator ${\mathcal S}_{\psi}: \mathcal H \to  L^2(\Lambda T^{\ast}M)$ by 
\[
{\mathcal S}_{\psi}H=\int_0^{+\infty} \psi_t(D)H_{t}\frac{dt}t
 =\lim_{\epsilon\to 0, N\to\infty}\int_\epsilon^N \psi_t(D)H_{t}\frac{dt}t\ \]
where the limit is in the $L^2(\Lambda T^{\ast}M)$ strong topology. This operator is also bounded, as ${\mathcal S}_{\psi}={\mathcal Q_{\overline\psi}}^*$  where 
$\overline\psi$ is defined by 
$\overline\psi(z) = \overline{\psi(\bar z)}$.

If $\widetilde \psi\in \Psi(\Sigma_{\theta}^0)$ is chosen to satisfy $\int_0^\infty\psi(\pm t)\widetilde\psi(\pm t)\frac{dt}t =1$ ({\it e.g.} by taking $\widetilde\psi(z)=\{\int_0^\infty|\psi(\pm t)|^2\frac{dt}t\}^{-1}\overline \psi(z)$ when $z\in \Sigma^{0}_{\theta^{\pm}}$), then the spectral theorem implies the 
following version of the Calder\'on reproducing theorem:
\[
{\mathcal S}_{\widetilde\psi}{\mathcal Q}_{\psi}f={\mathcal S}_{\psi}{\mathcal Q}_{\widetilde\psi}f = f
\]
for all $f\in \mathcal R(D)$ and hence for all $f\in H^2(\Lambda T^{\ast}M)$.
(Indeed ${\mathcal S}_{\psi}{\mathcal Q}_{\widetilde\psi}$ is the orthogonal projection of $L^2(\Lambda T^{\ast}M)$ onto $H^2(\Lambda T^{\ast}M)$.)
It follows that $\mathcal R({\mathcal S}_{\psi}) = H^2(\Lambda T^{\ast}M)$ and that 
\[
\left\Vert f\right\Vert_2\sim \inf\left\{\left\Vert H\right\Vert_{\mathcal H};\ f=S_{\psi}H\right\}
\]  for all $f\in H^2(\Lambda T^{\ast}M)$. 

\begin{rem} \label{care} With a little more care we could take $\widetilde\psi\in \Psi_{\sigma,\tau}(\Sigma_{\theta}^0)$ for any given $\sigma,\tau$. This fact will be used in Section \ref{firstdef}.
\end{rem} 

Thus, we have two descriptions of $H^2(\Lambda T^{\ast}M)$ in terms of quadratic functionals, involving the ${\mathcal H}$ space (which is nothing but the tent space $T^{2,2}(\Lambda T^{\ast}M)$ of Section \ref{tent} below) and independent of the choice of the function $\psi$. 

Let $\Delta=D^2$. Note that $\mathcal N(\Delta)= \mathcal N(D)$ and, hence,  $D$ and $\Delta$ are one-one operators on $H^2(\Lambda T^{\ast}M)$. Observe also that replacing $D$ by $\Delta$  and $\Psi(\Sigma_{\theta}^0)$ by $\Psi(\Sigma_{\theta^+}^0)$  would lead exactly to the similar descriptions of the Hardy space $H^2(\Lambda T^{\ast}M)$ in terms of functions of $\Delta$ only.

\bigskip

We define the Riesz transform on $M$ as  
the bounded operator $D\Delta^{-1/2} \colon H^2(\Lambda T^{\ast}M) \to H^2(\Lambda T^{\ast}M)$.

Set $H^2_d(\Lambda T^{\ast}M)=\overline{{\mathcal R}(d)}$ and $H^2_{d^{\ast}}(\Lambda T^{\ast}M)=\overline{{\mathcal R}(d^{\ast})}$, so that by the Hodge decomposition
\begin{equation} \label{Hodgereformul}
H^2(\Lambda T^{\ast}M)=H^2_d(\Lambda T^{\ast}M)\oplus H^2_{d^{\ast}}(\Lambda T^{\ast}M),
\end{equation}
and the sum is orthogonal. The orthogonal projections are given by 
$dD^{-1}$ and $d^{\ast}D^{-1}$. 

The Riesz transform $D\Delta^{-1/2}$ splits naturally as the sum 
of 
$d\Delta^{-1/2}$ and $d^*\Delta^{-1/2}$, which we call the Hodge-Riesz transforms.
As
$$
d\Delta^{-1/2} = (dD^{-1}) (D\Delta^{-1/2}) \quad \mathrm{and} \quad   d^*\Delta^{-1/2}= (d^*D^{-1}) (D\Delta^{-1/2}),
$$
they extend to bounded operators on $H^2(\Lambda T^{\ast}M)$. 
One further checks that $d\Delta^{-1/2}$ is bounded and invertible  from $H^2_{d^{\ast}}(\Lambda T^{\ast}M)$ to $H^2_d(\Lambda T^{\ast}M)$, that $d^{\ast}\Delta^{-1/2}$ is bounded and invertible  from $H^2_{d}(\Lambda T^{\ast}M)$ to $H^2_{d^{\ast}}(\Lambda T^{\ast}M)$,  and that they are inverse to one another.

\SE{Off-diagonal $L^2$-estimates for Hodge-Dirac and Hodge-Laplace 
operators} \label{OD(N)}

Throughout this section, $M$ is an arbitrary complete Riemannian manifold (we stress the fact that $M$ is not assumed to satisfy the doubling property (\ref{D})). We collect and prove all the 
off-diagonal 
$L^2$-estimates which will be used in the sequel for the Hodge-Dirac 
operator and the Hodge-de Rham Laplacian (and also for $d$ and $d^{\ast}$). 
We will make use of the following terminology: 
\begin{defi}
Let $A\subset \C$ be a non-empty set, $(T_z)_{z\in 
A}$ be a family of $L^2(\Lambda T^{\ast}M)$-bounded operators, 
$N\geq0$ and $C>0$. Say that $(T_z)_{z\in A}$ 
satisfies $OD_z(N)$ estimates with constant $C$ if, for all disjoint 
closed subsets 
$E,F\subset M$ and all $z\in 
A$,
\begin{equation} \label{OD}
\left\Vert M_{\chi_{F}}T_zM_{\chi_E}\right\Vert_{2,2}\leq C\inf\left(1,\left(\frac 
{\left\vert z\right\vert}{\rho(E,F)}\right)^{N}\right),
\end{equation}
where, for any $G\subset M$, $\chi_{G}$ denotes the characteristic 
function of $G$ and, for any bounded function $\eta$ on $M$, $M_{\eta}$ 
stands for the multiplication by $\eta$.
\end{defi}
In this definition and in the sequel, if $E$ and $F$ are any subsets 
of $M$, $\rho(E,F)$ is the infimum of $\rho(x,y)$ for all $x\in E$ 
and all $y\in F$. Moreover, if $T$ is a bounded linear operator from $L^p(\Lambda T^{\ast}M)$ to $L^q(\Lambda T^{\ast}M)$, its functional norm is denoted by $\left\Vert T\right\Vert_{q,p}$. 

\begin{rem} \label{ODrem} We remark that if $(T_z)_{z\in A}$ 
satisfies $OD_z(N)$ estimates, and  $0\leq N_1\leq N$, then $(T_z)_{z\in A}$ 
satisfies $OD_z(N_1)$ estimates. 
\end{rem}
The off-diagonal estimates to be used in the sequel 
will be presented in four lemmata. 
\begin{lem}\label{offdiag1}
Let $N$ and  $\alpha$  be nonnegative integers with $0\leq N\leq 
\alpha$ and $\displaystyle \mu\in 
\left(0,\frac{\pi}2\right)$. Then, for all integers $N^{\prime}\geq 
0$, 
$\left((zD)^{N}(I+izD)^{-\alpha}\right)_{z\in 
\Sigma_{\frac{\pi}2-\mu}}$  satisfies 
$OD_z(N^{\prime})$ estimates with  constants only depending on 
$\mu,N,N'$ and $\alpha$. 
\end{lem}
\begin{rem} 
Note that, with the same notations, if $\alpha\geq N+1$, $(zd(zD)^N(I+izD)^{-\alpha})_{z\in 
\Sigma_{\frac{\pi}2-\mu}}$ and $(zd^{\ast}(zD)^N(I+izD)^{-\alpha})_{z\in 
\Sigma_{\frac{\pi}2-\mu}}$ satisfy 
$OD_z(N^{\prime})$ estimates with  constants only depending on 
$\mu,N,N',\alpha$. However, these estimates will not be used in the sequel.
\end{rem}
\begin{lem} \label{offdiag2}
Let $k,N$ and $\alpha$ be nonnegative integers with $0\leq 
N\leq \alpha$ and $\displaystyle \mu\in 
\left(0,\frac{\pi}2\right)$. Then, for all 
$\displaystyle \tau\in \Sigma_{\frac{\pi}2-\mu}$, $\displaystyle 
\left((I+i\tau D)^{-k}(zD)^{N}(I+izD)^{-\alpha}\right)_{z\in 
\Sigma_{\frac{\pi}2-\mu}}$
satisfies $OD_z(N)$ estimates with constants only depending on 
$\mu,N,k,\alpha$ and $\beta$ (and, in particular, uniform in $\tau$). 
\end{lem}
Since the operator 
$D$ is self-adjoint in $L^2(\Lambda T^{\ast}M)$, one may 
define the $L^2$-bounded operator $f(D)$ for any $f\in 
H^{\infty}(\Sigma_{\theta}^0)$. If $f\in 
\Psi(\Sigma_{\theta}^0)$, then $f(D)$ can be computed with the Cauchy 
formula:
\begin{equation} \label{Cauchy}
f(D)=\frac 1{2i\pi} \int_{\gamma} (\zeta I-D)^{-1}f(\zeta)d\zeta,
\end{equation}
where $\gamma$ is made of two rays $re^{\pm i\beta}$, $r\geq 0$ and 
$\beta<\theta$, described counterclockwise (see \cite{macfunct}, 
\cite[Section 0.1]{asterisque}). Moreover, for every $f\in 
H^{\infty}(\Sigma_{\theta}^0)$ there is a uniformly bounded sequence 
of functions 
$f_n\in \Psi(\Sigma_{\theta}^0)$ which converges to $f$ uniformly on 
compact 
sets, and then $f(D)\psi(D) = \lim f_n(D)\psi(D)$ in the strong 
operator topology for all 
$\psi\in \Psi(\Sigma_{\theta}^0)$.

\begin{lem} \label{offdiag3} 
Let $N$ be a positive integer  and $\displaystyle \mu\in 
\left(0,\tfrac{\pi}2\right)$. 
\begin{itemize}
\item[$(a)$]
If $\left(g_{(t)}\right)_{t>0}$ is a uniformly bounded family of functions in $H^\infty(\Sigma_{\frac{\pi}2-\mu}^0)$ and $\alpha$ is an integer such that
$\alpha \geq N+1$, then 
$(g_{(t)}(D)(tD)^{N}(I\pm itD)^{-\alpha})_{t>0}$ satisfies $OD_t(N-1)$ 
estimates with constant bounded by $C\sup\limits_{t>0}\left\Vert 
g_{(t)}\right\Vert_{\infty}$.
\item[$(b)$] If $f\in H^\infty(\Sigma_{\frac{\pi}2-\mu}^0)$ and 
$\psi\in \Psi_{N,1}(\Sigma_{\pi/2-\mu}^0)$, then $(f(D)\psi_t(D))_{t>0}$
satisfies $OD_t(N-1)$ 
estimates with constants bounded by $C\left\Vert 
f\right\Vert_{\infty}$.
\end{itemize}
\end{lem}

In what follows, we set  $h_{a,b}(u)=\inf\left(u^{a},u^{-b}\right)$, 
where
 $a,b,u>0$. Recall that if $\psi\in \Psi(\Sigma_{\theta}^0)$ and $t>0$, then $\psi_t$ is defined by $\psi_t(z)=\psi(tz)$.
 
\begin{lem} \label{offdiag4} Let $\psi\in \Psi_{N_1,\alpha_1}(\Sigma_{\pi/2-\mu}^0)$
and $\widetilde \psi\in \Psi_{N_2,\alpha_2}(\Sigma_{\pi/2-\mu}^0)$ where 
$\alpha_1,\alpha_2,N_1,N_2$ are positive integers
and $\mu\in \left(0,\frac{\pi}2\right)$, and suppose that $a,b$ are nonnegative integers satisfying $a\leq\min\{N_1,\alpha_2-1\}$,
$b\leq\min\{N_2,\alpha_1-1\}$. Then, there exists $C>0$ 
such that, for 
all $f\in H^{\infty}(\Sigma^{0}_{\mu})$, there exists, for all 
$s,t>0$, an operator $T_{s,t}$ with the following properties:
\begin{itemize}
\item[$(i)$]
$\displaystyle 
\psi_s(D)f(D)\widetilde\psi_t(D)=h_{a,b}\left(\frac 
st\right) T_{s,t}$; 
\item[$(ii)$]
 $(T_{s,t})_{t\geq s}$ satisfies 
$OD_{t}(N_{2}+a-1)$ 
estimates uniformly in $s>0$;
\item[$(iii)$]
$(T_{s,t})_{s\geq t}$ satisfies 
$OD_{s}(N_{1}+b-1)$ 
estimates uniformly in $t>0$.
\end{itemize}
\end{lem} 

\noindent{\bf Proof of Lemma \ref{offdiag1}: } The proof is exactly 
as the one of Proposition 5.2 in \cite{akm}.
\hfill\fin\par

\bigskip

\noindent{\bf Proof of Lemma \ref{offdiag2}: } We use the notation
\[
\left[T,S\right]=TS-ST
\]
for the commutator of two operators $T$ and $S$. The proof is done by induction on $k$ and relies on a 
commutator argument, as in \cite{akm}, Proposition 5.2. \par
For $k=0$, the conclusion is given by Lemma \ref{offdiag1}. 
Let $k\geq 1$ and assume that, for all integers $0\leq N\leq \alpha$, 
$\left((I+i\tau D)^{-(k-1)}(zD)^{N}(I+izD)^{-\alpha}\right)_{z\in 
\Sigma_{\frac{\pi}2-\mu}}$ satisfies $OD_z(N)$ estimates uniformly in 
$\tau$. To establish that $\left((I+i\tau 
D)^{-k}(zD)^{N}(I+izD)^{-\alpha}\right)_{z\in 
\Sigma_{\frac{\pi}2-\mu}}$ satisfies $OD_z(N)$ estimates uniformly in 
$\tau$ whenever $0\leq N\leq \alpha$, we argue by induction on $N$. 
The case when $N=0$ is obvious. Let $1\leq N\leq \alpha$ and assume 
that 
$\displaystyle \left((I+i\tau 
D)^{-k}(zD)^{(N-1)}(I+izD)^{-\alpha}\right)_{z\in 
\Sigma_{\frac{\pi}2-\mu}}$ satisfies $OD_z(N-1)$ estimates uniformly 
in $\tau$. We intend to show that $\displaystyle \left((I+i\tau 
D)^{-k}(zD)^{N}(I+izD)^{-\alpha}\right)_{z\in 
\Sigma_{\frac{\pi}2-\mu}}$ satisfies $OD_z(N)$ estimates uniformly in 
$\tau$. 
Let $E,F$ be two disjoint closed subsets of $M$, $\chi$ the characteristic 
function of $E$ and $\eta$ a Lipschitz function on $M$ equal to $1$ on 
$F$, to $0$ on $E$ and satisfying
\[
\left\Vert \left\vert \nabla \eta\right\vert\right\Vert_{\infty}\leq 
C\rho(E,F)^{-1}, \rho(\mbox{supp }\eta,E)\sim \rho(E,F).
\]
Our conclusion reduces to proving that
\begin{equation} \label{commut2}
\left\Vert M_{\eta} \left((I+i\tau 
D)^{-k}(zD)^{N}(I+izD)^{-\alpha}\right)M_{\chi}\right\Vert_{2,2}\leq 
C\left(\frac {\left\vert z\right\vert}{\rho(E,F)}\right)^N
\end{equation}
where we recall that $M_{\eta}$ and $M_{\chi}$ denote the multiplication by $\eta$ 
and 
$\chi$ respectively. But, because of the supports of $\chi$ and 
$\eta$, the left hand side of (\ref{commut2}) is equal to the $\left\Vert . \right\Vert_{2,2}$ norm of 
\begin{equation}\label{sum}
\begin{array}{lll}
&&\left[M_{\eta}, \left((I+i\tau 
D)^{-k}(zD)^{N}(I+izD)^{-\alpha}\right)\right]M_{\chi}\\ 
& = & \displaystyle (I+i\tau D)^{-1}\left[M_{\eta}, \left((I+i\tau 
D)^{-(k-1)}(zD)^{N}(I+izD)^{-\alpha}\right)\right]M_{\chi}\\ 
& &+  \displaystyle \left[M_{\eta},(I+i\tau D)^{-1}\right] 
\left((I+i\tau 
D)^{-(k-1)}(zD)^{N}(I+izD)^{-\alpha}\right)M_{\chi}.
\end{array}
\end{equation}
By the induction assumption, the $\left\Vert .\right\Vert_{2,2}$ norm of the first term is bounded by
\[
\left\Vert \left[M_{\eta}, \left((I+i\tau 
D)^{-(k-1)}(zD)^{N}(I+izD)^{-\alpha}\right)\right]M_{\chi}\right\Vert_{2,2} 
\leq C\left(\frac {\left\vert z\right\vert}{\rho(E,F)}\right)^N.
\]
The second term in (\ref{sum}) is equal to
\[
(I+i\tau D)^{-1} z \left[D,M_{\eta}\right](i\tau 
D)(I+i\tau D)^{-k}  
(zD)^{N-1}(I+izD)^{-\alpha}M_{\chi}
\] 
and its $\left\Vert .\right\Vert_{2,2}$ norm is therefore bounded by
\[
\begin{array}{ll}
\displaystyle \left\Vert (I+i\tau D)^{-1} 
z\left[D,M_{\eta}\right] (I+i\tau 
D)^{-(k-1)} 
(zD)^{N-1}(I+izD)^{-\alpha}M_{\chi}\right\Vert_{2,2}
& + \\
\displaystyle \left\Vert (I+i\tau D)^{-1} 
z\left[D,M_{\eta}\right] (I+i\tau D)^{-k} (zD)^{N-1} 
(I+izD)^{-\alpha}M_{\chi}\right\Vert_{2,2}\\
\leq \displaystyle \left\vert z\right\vert \left\Vert \left\vert 
\nabla \eta\right\vert\right\Vert_{\infty} \left(\frac {\left\vert 
z\right\vert}{\rho(E,F)}\right)^{N-1}\\
\leq \displaystyle C \left(\frac {\left\vert 
z\right\vert}{\rho(E,F)}\right)^N,
\end{array}
\]
where the penultimate inequality follows from the induction 
assumptions and the formula
\begin{equation} \label{development}
D(\eta b)=\eta Db+d\eta\wedge b-d\eta\vee b,
\end{equation}
where
\[
\langle \alpha\vee \beta,\gamma\rangle:=\langle \beta,\alpha\wedge 
\gamma\rangle.
\]
This concludes the proof of (\ref{commut2}), and therefore of Lemma \ref{offdiag2}.
\hfill\fin\par

\bigskip

\noindent{\bf Proof of Lemma \ref{offdiag3}: }We begin with assertion $(a)$.
First note that for $f\in \Psi(\Sigma_{\pi/2-\mu}^0)$ and 
$0<r<R<\infty$, then $$\left|\frac 1{2i\pi}\int_{\zeta\in\gamma;r\leq |\zeta|\leq 
R} f(\zeta)\tfrac1\zeta\,d\zeta\right|\leq 2||f||_\infty\ .$$
To see this, apply Cauchy's theorem to change to an integral over 
four arcs.
This fact is used to handle the second last term in the following 
expression.
\begin{align*}
g_{(t)}(D) (tD)^{N}(I+itD)^{-\alpha}&=\frac 1{2i\pi} \int_\gamma 
g_{(t)}(\zeta) (\zeta I-D)^{-1}(tD)^{N}(I+itD)^{-\alpha}\,d\zeta\\
&=\frac t{2i\pi} 
\int_{\zeta\in\gamma;|\zeta|<1/t}g_{(t)}(\zeta) \tfrac1\zeta D( 
I-\tfrac1\zeta D)^{-1}(tD)^{N-1}(I+itD)^{-\alpha}\,d\zeta\\
&\quad+\lim_{R\to\infty}\frac 
1{2i\pi}\int_{\zeta\in\gamma;1/t\leq|\zeta|\leq 
R}g_{(t)}(\zeta) \tfrac1\zeta \,(tD)^{N}(I+itD)^{-\alpha}  \, d\zeta\\
&\quad+\lim_{R\to\infty}\frac 
1{2i\pi}\int_{\zeta\in\gamma;1/t\leq|\zeta|\leq R} g_{(t)}(\zeta) \tfrac 
1{t\zeta^2} ( I-\tfrac1\zeta 
D)^{-1}(tD)^{N+1}(I+itD)^{-\alpha}\,d\zeta
\end{align*}
Apply Lemmas \ref{offdiag1} and \ref{offdiag2} to see that each term 
satisfies $OD_t(N-1)$ estimates. 
A limiting argument gives the result for a family $(g_{(t)})_{t>0}$ uniformly bounded in $H^\infty(\Sigma_{\pi/2-\mu}^0)$. 

To prove assertion $(b)$ in Lemma \ref{offdiag3}, apply assertion $(a)$ with $g_t(z)=f(z)\psi_t(z)(tz)^{-N}(1+itz)^{(N+1)}$.
 \hfill\fin\par
 
 \bigskip
 
\noindent{\bf Proof of Lemma \ref{offdiag4}: }If $s\leq t$, write
\[
\psi_s(D)f(D)\widetilde\psi_t(D)=\left(\frac 
st\right)^{a}(sD)^{-a}\psi_s(D)f(D)(tD)^a\widetilde\psi_t(D)=
\left(\frac st\right)^{a}T_{s,t}
\]
where 
\[
T_{s,t}=f_{(s)}(D)\widetilde{\widetilde\psi}_t(D)\ .
\]
with $f_{(s)}(z)=(sz)^{-a}\psi(sz)f(z)$ and $\widetilde{\widetilde\psi}(z)=z^a\widetilde\psi(z)$.
Now $f_{(s)}\in H^\infty(\Sigma_{\pi/2-\mu}^0)$ with 
$\left\Vert f_{(s)}\right\Vert_{\infty}\leq C_1\left\Vert 
f\right\Vert_{\infty} $, and $\widetilde{\widetilde\psi}\in \Psi_{N_2+a,\alpha_2-a}$, so
Lemma \ref{offdiag3} ensures that $T_{s,t}$ satisfies 
$OD_t(N_2+a-1)$ estimates with a constant not exceeding 
$C\left\Vert 
f\right\Vert_{\infty} $. The part $t\leq s$ is proved in a similar 
way. \hfill\fin

We remark that, since $\Delta =D^{2}$, Lemmata \ref{offdiag1}, 
\ref{offdiag2}, \ref{offdiag3} and \ref{offdiag4}
imply similar off-diagonal estimates when $(tD)^N(I+itD)^{-\alpha}$ is replaced by 
$(t^{2}\Delta)^{N}(I+t^{2}\Delta)^{-\alpha}$ for appropriate $N$ 
and $\alpha$. Furthermore, we can strengthen these to
 ``Gaffney'' type estimates for 
the heat semigroup. 
\begin{lem} \label{gaffneyheat}
For all $N\geq 0$, there exists $C,\alpha>0$ such that, for all 
disjoint closed subsets 
$E,F\subset M$ and all $t>0$,
\[
 \left\Vert 
M_{\chi_{F}}(t^{2}\Delta)^N e^{-t^{2}\Delta}M_{\chi_E}\right\Vert_{2,2} + 
\left\Vert 
M_{\chi_{F}}tD(t^{2}\Delta)^Ne^{-t^{2}\Delta}M_{\chi_E}\right\Vert_{2,2} \leq 
Ce^{-\alpha 
\frac{\rho^{2}(E,F)}{t^{2}}}.
\]
In particular, $((t^{2}\Delta)^N e^{-t^{2}\Delta})_{t>0}$ and 
$(tD(t^{2}\Delta)^Ne^{-t^{2}\Delta})_{t>0}$ satisfy 
$OD_{t}(N^{\prime})$ estimates for any integer $N^{\prime}\geq 0$.
\end{lem}
\noindent{\bf Proof: }The proof of the estimate for the first term is 
analogous to 
\cite{davies1} (this kind of estimate originated in Gaffney's work 
\cite{gaffney}) and \cite{davies2}, Lemma 7, whereas the second term can be estimated by the same method as in \cite{acdh}, estimate 
(3.1) p. 930. \hfill\fin

Observe that the same argument yields
\[
\left\Vert 
M_{\chi_{F}}tde^{-t^{2}\Delta}M_{\chi_E}\right\Vert_{2,2} + \left\Vert 
M_{\chi_{F}}td^{\ast}e^{-t^{2}\Delta}M_{\chi_E}\right\Vert_{2,2}\leq Ce^{-\alpha 
\frac{\rho^{2}(E,F)}{t^{2}}}.
\]
\SE{Tent spaces on $M$} \label{tent}
\subsection{Definition, atomic decomposition and duality for tent 
spaces} \label{deftent}
We first present tent spaces on $M$, following 
\cite{coifmanmeyerstein}. For all $x\in M$ and $\alpha>0$, the cone 
of aperture $\alpha$ and vertex $x$ is the set
\[
\Gamma_{\alpha}(x)=\left\{(y,t)\in M\times \left(0,+\infty\right);\ 
y\in B(x,\alpha t)\right\}.
\]
When $\alpha=1$, $\Gamma_{\alpha}(x)$ will simply be denoted by 
$\Gamma(x)$. For any closed set $F\subset M$, let ${\mathcal R}(F)$ 
be the union of all cones with aperture $1$ and vertices in $F$. 
Finally, if $O\subset M$ is an open set and $F=M\setminus O$, the 
tent 
over $O$, denoted by $T(O)$, is the complement of ${\mathcal R}(F)$ 
in $M\times \left(0,+\infty\right)$.  

Let $F=(F_t)_{t>0}$ be a family of measurable sections of $\Lambda T^{\ast}M$. Write $F(y,t):=F_t(y)$ for all $y\in M$ and all $t>0$ and assume that $F$ is measurable on $M\times (0,+\infty)$. Define then, for all $x\in M$,
\[
\S F(x)=\left(\iint_{\Gamma(x)} \left\vert F(y,t)\right\vert^2 
\frac{dy}{V(x,t)} \frac {dt}t\right)^{1/2},
\]
and, if $1\leq p<+\infty$, say that $F\in T^{p,2}(\Lambda T^{\ast}M)$ if
\[
\left\Vert F\right\Vert_{T^{p,2}(\Lambda T^{\ast}M)}:=\left\Vert \S 
F\right\Vert_{L^p(M)}<+\infty.
\]
\begin{rem} \label{aperture}
Assume that (\ref{D}) holds. If $\alpha>0$ and if we define, for all $x\in M$, 
\[
\S_{\alpha}F(x)=\left(\iint_{\Gamma_{\alpha}(x)} \left\vert F(y,t)\right\vert^2 
\frac{dy}{V(x,t)} \frac {dt}t\right)^{1/2},
\] 
then $\left\Vert F\right\Vert_{T^{p,2}(\Lambda T^{\ast}M)}\sim \left\Vert \S_{\alpha}F\right\Vert_{L^p(M)}$ for all $1\leq p<+\infty$ (see \cite{coifmanmeyerstein}). 
\end{rem}
In order to ensure duality results for tent spaces, we do not define 
$T^{\infty,2}(\Lambda T^{\ast}M)$ in the same way. For any family $(F_t)_{t>0}$ of measurable sections of 
$\Lambda T^{\ast}M$ and all $x\in 
M$, define
\[
\c F(x)=\sup_{B\ni x} \left(\frac 1{V(B)} \iint_{T(B)} 
\left\vert F(y,t)\right\vert^{2} dy\frac{dt}t\right)^{1/2},
\]
where the supremum is taken over all open balls $B$ containing $x$, 
and 
say that $F\in T^{\infty,2}(\Lambda T^{\ast}M)$ if $\left\Vert 
F\right\Vert_{T^{\infty,2}(\Lambda T^{\ast}M)}:=\left\Vert \c 
F\right\Vert_{\infty}<+\infty$.

We first state a density result for tent spaces which does not require  (\ref{D})
\begin{pro} \label{density} When $1\leq p<\infty$, then $T^{p,2}(\Lambda T^{\ast}M)\cap T^{2,2}(\Lambda T^{\ast}M)$ is dense in $T^{p,2}(\Lambda T^{\ast}M)$.
\end{pro}

\noindent{\bf Proof: }Set
\[
{\mathcal E}=\left\{F\in T^{2,2}(\Lambda T^{\ast}M);\ F\mbox{ is bounded and has compact support in }M\times \left(0,+\infty\right)\right\},
\]
which is obviously contained in $T^{p,2}(\Lambda T^{\ast}M)\cap T^{2,2}(\Lambda T^{\ast}M)$. Fix any point $x_0$ in $M$ and, for all $n\geq 1$, define $\chi_n=\chi_{B(x_0,n)\times \left(\frac 1n,n\right)}$. Then, it is easy to check that, for all $F\in T^{p,2}(\Lambda T^{\ast}M)$, if
\[
F_n=\chi_n\chi_{\left\{(x,t)\in M\times (0,+\infty);\ \left\vert F(x,t)\right\vert<n\right\}}F
\]
for all $n\geq 1$, then $F_n\in {\mathcal E}$ and $F_n\rightarrow F$ in $T^{p,2}(\Lambda T^{\ast}M)$.\hfill\fin
\begin{rem} \label{truncation}
The same argument shows that, if $F\in T^{p,2}(\Lambda T^{\ast}M)\cap T^{2,2}(\Lambda T^{\ast}M)$, then $\chi_nF\rightarrow F$ both in $T^{p,2}(\Lambda T^{\ast}M)$ and in $T^{2,2}(\Lambda T^{\ast}M)$.
\end{rem}

If we assume furthermore property (\ref{D}), duality, atomic decomposition and 
interpolation results hold for tent spaces as in the Euclidean case. The proofs are analogous 
to the corresponding ones in \cite{coifmanmeyerstein}, and we will 
therefore not write them down (see however \cite{russtent} for the atomic decomposition for tent spaces on spaces of homogeneous type). Let us 
just mention that, apart from property (\ref{D}), these proofs rely on the existence of 
$\alpha>0$ 
such that, for all $r>0$ and all $x,y\in M$ satisfying 
$\rho(x,y)<r$,
\[
\mu(B(x,r)\cap B(y,r))\geq \alpha V(x,r).
\]
This last assertion follows from the definition of the geodesic 
distance on $M$, the completeness of $M$ and the doubling property. 

\medskip

The duality for tent spaces is as follows:
\begin{theo} \label{duality}
Assume  (\ref{D}). Then:
\begin{itemize}
\item[$(a)$]
There exists $C>0$ such that, for all $F\in T^{1,2}(\Lambda T^{\ast}M)$ and all 
$G\in 
T^{\infty,2}(\Lambda T^{\ast}M)$,
\[
\iint_{M\times \left(0,+\infty\right)} \left\vert F(x,t)\right\vert 
\left\vert G(x,t)\right\vert dx\frac{dt}t \leq C\int_{M}\S F(x)\c 
G(x)dx.
\]
\item[$(b)$]
The pairing $\langle F,G\rangle\mapsto \iint_{M\times 
\left(0,+\infty\right)} \langle F(x,t),G(x,t)\rangle dx\frac{dt}t$ realizes 
$T^{\infty,2}(\Lambda T^{\ast}M)$ as equivalent with the dual of $T^{1,2}(\Lambda T^{\ast}M)$ 
and 
$T^{p^{\prime},2}(\Lambda T^{\ast}M)$ as equivalent with the dual of 
$T^{p,2}(\Lambda T^{\ast}M)$ if 
$1<p<+\infty$ and $1/p+1/p^{\prime}=1$.
\end{itemize}
\end{theo}
In assertion $(b)$ and in the sequel, $\langle .,.\rangle$ denotes the complex inner product in $\Lambda T^{\ast}M$. \par

The $T^{p,2}(\Lambda T^{\ast}M)$ spaces interpolate by the complex interpolation method:
\begin{theo} \label{interpolation}
Assume (\ref{D}). Let $1\leq p_{0}<p<p_{1}\leq +\infty$ and $\theta\in 
\left(0,1\right)$ such that $1/p=(1-\theta)/p_{0}+\theta/p_{1}$. 
Then 
$\left[T^{p_{0},2}(\Lambda T^{\ast}M),T^{p_{1},2}(\Lambda T^{\ast}M)\right]_{\theta}=T^{p,2}(\Lambda T^{\ast}M)$. 
\end{theo}

The $T^{1,2}(\Lambda T^{\ast}M)$ space admits an atomic decomposition. An atom is 
a 
function $A\in L^{2}\left((0,+\infty),L^2(\Lambda T^{\ast}M),dt/t\right)$ supported 
in $T(B)$ for some ball 
$B\subset M$ and satisfying
\[
\iint_{T(B)} \left\vert A(x,t)\right\vert^{2} dx\frac{dt}t\leq \frac 
1{V(B)}.
\]
An atom belongs to $T^{1,2}(\Lambda T^{\ast}M)$ with a 
norm 
controlled by a constant only depending on $M$. It turns out that 
every $F\in T^{1,2}(\Lambda T^{\ast}M)$ has an atomic decomposition (see \cite{russtent}):
\begin{theo} \label{atomic}
Assume  (\ref{D}). There exists $C>0$ such that every $F\in T^{1,2}(\Lambda T^{\ast}M)$ can be 
written 
as $F=\sum_{j}\lambda_jA_{j}$, where the $A_{j}$'s are atoms and 
$\sum_{j\geq 0} \left\vert \lambda_{j}\right\vert \leq C\left\Vert 
F\right\Vert_{T^{1,2}(\Lambda T^{\ast}M)}$.
\end{theo}
\begin{rem}\label{box}
It is plain to see that, in the definition of an atom, up to changing 
the constant in Theorem \ref{atomic}, the tent $T(B)$ over the ball 
$B$ can be replaced by the Carleson box
\[
{\mathcal B}(B)=B\times \left[0,r(B)\right]
\]
where $r(B)$ is the radius of $B$.
\end{rem}
We end up this section by a technical lemma for later use:
\begin{lem} \label{t1t2conv}
Assume (\ref{D}).
\begin{itemize}
\item[$(a)$]
Let $1\leq p<+\infty$. If $(H_n)_{n\geq 1}$ is any sequence in $T^{p,2}(\Lambda T^{\ast}M)$ which converges to $H\in T^{p,2}(\Lambda T^{\ast}M)$, there exists an increasing map $\varphi:\N^{\ast}\rightarrow \N^{\ast}$ such that $H_{\varphi(n)}(y,t)\rightarrow H(y,t)$ for almost every $(y,t)\in M\times (0,+\infty)$. 
\item[$(b)$] 
Let $(H_n)_{n\geq 1}$ be a sequence in $T^{1,2}(\Lambda T^{\ast}M)\cap T^{2,2}(\Lambda T^{\ast}M)$ which converges to $H$ in $T^{1,2}(\Lambda T^{\ast}M)$ and to $G$ in $T^{2,2}(\Lambda T^{\ast}M)$. Then, $H=G$.
\end{itemize}
\end{lem}
\noindent{\bf Proof: }For assertion $(a)$, since, for all $j\geq 1$, ${\mathcal S}_j(H_n-H)\rightarrow 0$ in $L^p(\Lambda T^{\ast}M)$ (see Remark \ref{aperture}), a diagonal argument shows that, up to a subsequence, ${\mathcal S}_j(H_n-H)(x)\rightarrow 0$ for all $j\geq 1$ and almost every $x\in M$. Fix then $x\in M$ such that $\S_j(H_n-H)(x)\rightarrow 0$ for all $j\geq 1$. Thanks to a diagonal argument again, one has $\left\vert (H_n-H)(y,t)\right\vert\rightarrow 0$ for almost every $(y,t)\in \Gamma_j(x)$, up to a subsequence, which gives the conclusion. Assertion $(b)$ is an immediate consequence of assertion $(a)$. \hfill\fin 
\subsection{The main estimate} \label{mainestim}
\noindent Recall from Section \ref{h2} that $\displaystyle 
{\mathcal H}=L^2\left((0,+\infty),L^2(\Lambda T^{\ast}M),\frac{dt}t\right)$, equipped with 
the norm 
\[
\left\Vert F\right\Vert_{\mathcal H}=\left(\int_0^{+\infty}\int_M \left\vert 
F(x,t)\right\vert^2dx\frac{dt}t\right)^{1/2}.
\]
It is easy to see that $\mathcal H = T^{2,2}(\Lambda T^{\ast}M)$ with equivalent norms.

The main result of the present section is the following theorem, 
which 
will play a crucial role in our definition of Hardy spaces via tent spaces 
on 
$M$, and in establishing the functional calculus for Hardy spaces.

\begin{theo} \label{projector} Assume that $M$ is a complete 
connected Riemannian manifold which satisfies  the doubling property  (\ref{D}). Define $\kappa$ as in (\ref{Dbis}), and let $\beta=\left[\frac\kappa 2\right]+1$ (the smallest integer larger than $\frac\kappa 2$)
and $\theta\in \left(0,\frac{\pi}2\right)$. For given $\psi,\widetilde\psi\in \Psi(\Sigma_{\theta}^0)$ and $f\in H^{\infty}(\Sigma_{\theta}^0)$,
define the bounded operator $Q_f:\mathcal H \to \mathcal H$ to be $Q_f = \mathcal Q_\psi f(D)\mathcal S_{\widetilde \psi}$, {\it i.e.},
\[
Q_f(F)_s=\int_0^{+\infty} 
\psi_s(D)f(D)\widetilde{\psi}_t(D)F_t\frac{dt}t
\]
for all $F\in {\mathcal H}$ and all $s>0$. Suppose either \newline
\noindent(a) $1\leq p<2$ and $\psi\in \Psi_{1,\beta+1}(\Sigma_{\theta}^0)$,
 $\widetilde \psi\in \Psi_{\beta,2}(\Sigma_{\theta}^0)$; or \newline
\noindent(b) $2<p\leq\infty$, 
and $ \psi\in \Psi_{\beta,2}(\Sigma_{\theta}^0)$,
$\widetilde\psi\in \Psi_{1,\beta+1}(\Sigma_{\theta}^0)$.

Then  $Q_f$ extends to a 
$T^{p,2}(\Lambda T^{\ast}M)$-bounded map, and, for all 
$F\in 
T^{p,2}(\Lambda T^{\ast}M)$,
\begin{equation} \label{eqprojector}
\left\Vert Q_f(F)\right\Vert_{T^{p,2}(\Lambda T^{\ast}M)}\leq C_p\left\Vert 
f\right\Vert_{\infty} \left\Vert F\right\Vert_{T^{p,2}(\Lambda T^{\ast}M)},
\end{equation}
where $C_{p}>0$ only depends on the constant in (\ref{D}), $\kappa$, $\theta$, 
$p$, $\psi$ and $\widetilde\psi$. 
\end{theo}

\begin{rem}\label{interp} In the case when $\int_0^\infty\psi(\pm t)\widetilde\psi(\pm t)\frac{dt}t =1$ and hence
${\mathcal S}_{\widetilde\psi}{\mathcal Q}_{\psi}h = h$ for all $h\in\mathcal R(D)$, then $Q_fQ_g=Q_{fg}$. 
In particular $\mathcal P_{\{\psi,\widetilde \psi\}}:= Q_1$ is a bounded projection on $T^{p,2}(\Lambda T^{\ast}M)$ when $1\leq p \leq \infty$.  In order that these operators be the same for all $p$, choose $\psi=\widetilde\psi\in \Psi_{1,\beta+1}(\Sigma_{\theta}^0)\cap
\Psi_{\beta,2}(\Sigma_{\theta}^0)$ and set $\mathcal P_{\psi}:=\mathcal P_{\{\psi,\widetilde \psi\}}={\mathcal Q}_{\psi}{\mathcal S}_{\psi}$.
In this case we see that the spaces
$\mathcal P_\psi T^{p,2}(\Lambda T^{\ast}M)$ interpolate by the complex method for $1\leq p\leq\infty$.
\end{rem}

\noindent{\bf Proof of Theorem \ref{projector}: } This proof will be divided in several steps.

\medskip

\noindent{\bf Step 1: }The boundedness of $Q_f$ in $T^{2,2}(\Lambda T^{\ast}M)=\mathcal H$ follows immediately from the results in Section \ref{h2}.\par

\medskip

\noindent{\bf Step 2: An inequality for $T^{1,2}(\Lambda T^{\ast}M)$ atoms. }We now assume that $\psi\in \Psi_{1,\beta+1}(\Sigma_{\theta}^0)$ and
 $\widetilde \psi\in \Psi_{\beta,2}(\Sigma_{\theta}^0)$. Let us prove that, for any atom $A\in 
T^{1,2}(\Lambda T^{\ast}M)$,
\begin{equation} \label{atom}
\left\Vert Q_f(A)\right\Vert_{T^{1,2}(\Lambda T^{\ast}M)}\leq C.
\end{equation}
Let $A$ be an atom in $T^{1,2}(\Lambda T^{\ast}M)$. There exists a ball $B\subset 
M$ such that $A$ is supported in $T(B)$ and
\[
\iint_{T(B)} \left\vert A(x,t)\right\vert^2 dx\frac{dt}t\leq 
V^{-1}(B).
\]
Set $\widetilde{A}=Q_f(A)$, $\widetilde{A_1}=\widetilde{A}\chi_{T(4B)}$ 
and, for all $k\geq 2$, $\widetilde{A_k}=\widetilde{A} 
\chi_{T(2^{k+1}B)\setminus T(2^kB)}$, so that 
$\widetilde{A}=\sum\limits_{k\geq 1} \widetilde{A_k}$ (actually, we 
should truncate $\widetilde{A}$ by imposing $\delta\leq s\leq R$, 
obtain bounds independent of $\delta$ and $R$, and then let 
$\delta$ go to $0$ and $R$ to $+\infty$; we ignore this point and 
argue directly without this truncation, to simplify the notation). 

We need  to show that, for some $\varepsilon>0$ and $C>0$ 
independent of $k,A$ and $f$, 
$\frac{2^{k\varepsilon}}{C}\widetilde{A_k}$ is a 
$T^{1,2}$ atom, which will prove that $\widetilde{A}\in 
T^{1,2}(\Lambda T^{\ast}M)$ 
with a controlled norm, by the atomic decomposition of 
$T^{1,2}(\Lambda T^{\ast}M)$. Since 
$\widetilde{A_k}$ is supported in $T\left(2^{k+1}B\right)$, it is 
enough to check that, for all $k\geq 1$,
\begin{equation} \label{atomest}
\iint \left\vert \widetilde{A_k}(x,s)\right\vert^2 dx\frac{ds}s\leq 
\frac{C^2}{V(2^{k+1}B)}2^{-2k\varepsilon}.
\end{equation}
For $k=1$, using the $T^{2,2}(\Lambda T^{\ast}M)$ boundedness of $Q_f$, the fact 
that 
$\left\Vert A\right\Vert_{\mathcal H}\leq V(B)^{-1/2}$ and the doubling 
property, one obtains
\[
\left\Vert \widetilde{A_1}\right\Vert_{\mathcal H}\leq C\left\Vert 
A\right\Vert_{{\mathcal H}}\leq CV^{-1/2}(B)\leq C^{\prime}V^{-1/2}(4B).
\]  
Fix now $k\geq 2$, and suppose $0<\delta < \beta-\D/2$. Applying Lemma \ref{offdiag4} with $a=1, b=\beta, N_1=1, N_2= \beta,  
\alpha_1= \beta + 1, \alpha _2= 2$ and the fact that  $A$ is supported in $T(B)$, write
\[
\widetilde{A}_{s}= Q(A)_s= \int_0^{+\infty} 
\psi_s(D)f(D)\widetilde{\psi}_t(D)A_t\frac{dt}t= \int_0^r h_{1,\beta}\left(\frac 
st\right)T_{s,t}A_t\frac{dt}t,
\]
where $r$ is the radius of $B$. The Cauchy-Schwarz inequality yields
\[
\begin{array}{lll}
\displaystyle\left\vert \widetilde{A_k}_{s}\right\vert^2 &\leq & \displaystyle 
\left(\int_0^{+\infty} h_{1,\delta}\left(\frac 
st\right)\frac{dt}t\right) \left(\int_0^r h_{1,2\beta-\delta}\left(\frac 
st\right)\left\vert T_{s,t}A_t\right\vert^2\frac{dt}t\right)\\
& \leq & \displaystyle C \int_0^r h_{1,2\beta-\delta}\left(\frac 
st\right)\left\vert T_{s,t}A_t\right\vert^2\frac{dt}t.
\end{array}
\]
Since $\widetilde{A_k}$ is supported in 
$T\left(2^{k+1}B\right)\setminus T\left(2^kB\right)$, one may assume 
that $0<s<2^{k+1}r$. Moreover, if $s<2^{k-1}r$ and 
if $(x,s)$ belongs to $T(2^{k+1}B)\setminus T(2^kB)$, then $x\in 
2^{k+1}B\setminus 2^{k-1}B$ , so that
\begin{equation} \label{widetildeak}
\begin{array}{lll}
\displaystyle \iint \left\vert \widetilde{A_k}(x,s)\right\vert^2 
dx\frac{ds}s & \leq & \displaystyle C \int_0^{2^{k-1}r} \int_0^r 
h_{1,2\beta-\delta}
\left(\frac st\right)\left\Vert 
\chi_{2^{k+1}B\setminus 
2^{k-1}B}T_{s,t}A_t\right\Vert_{L^2(\Lambda T^{\ast}M)}^2\frac{dt}t\frac{ds}s \\
& + & \displaystyle C \int_{2^{k-1}r}^{2^{k+1}r} \int_0^r 
h_{1,2\beta-\delta}\left(\frac st\right)\left\Vert 
T_{s,t}A_t\right\Vert_{L^2(\Lambda T^{\ast}M)}^2\frac{dt}t\frac{ds}s\ .
\end{array}
\end{equation} 
Thanks to (\ref{Dbis}), the last integral in (\ref{widetildeak}) is bounded by
\[
\begin{array}{lll}
\displaystyle C\int_{2^{k-1}r}^{2^{k+1}r}\int_0^r \left(\frac 
ts\right)^{2\beta-\delta} \left\Vert 
A_t\right\Vert_{L^2(\Lambda T^{\ast}M)}^2\frac{dt}t\frac{ds}s & \leq & \displaystyle C\int_0^r \left(\frac t{2^kr}\right)^{2\beta-\delta} 
\left\Vert A_t\right\Vert_{L^2(\Lambda T^{\ast}M)}^2\frac{dt}t\\
& \leq & \displaystyle C2^{-k(2\beta-\delta)}V^{-1}(B)\\
& \leq & \displaystyle C2^{-k(2\beta-\delta-\D)}V^{-1}(2^{k+1}B),
\end{array}
\]
where we now need the fact that $2\beta-\delta-\D>0$. Moreover, Lemma \ref{offdiag4} yields that $(T_{s,t})_{s\geq t}$ satisfies $OD_s(\beta)$ estimates, and hence $OD_s(\beta-\delta)$ estimates by Remark \ref{ODrem}. The first integral on 
the right hand side of (\ref{widetildeak}) is therefore dominated by
\[\begin{array}{ll}
\displaystyle C\int_0^r \left\Vert A_t\right\Vert_{L^2(\Lambda T^{\ast}M)}^2 
\left(\int_0^t \left(\frac st\right) \left(\frac 
t{2^kr}\right)^{2\beta}\frac{ds}s+\int_t^{2^{k-1}r} 
\left(\frac ts\right)^{2\beta-\delta} \left(\frac 
s{2^kr}\right)^{2\beta-2\delta}\frac{ds}s\right)\frac{dt}t & \leq \\
\displaystyle C\int_0^r \left\Vert A_t\right\Vert_{L^2(\Lambda T^{\ast}M)}^2 
\left(\left(\frac t{2^kr}\right)^{2\beta}+\left(\frac t 
{2^kr}\right)^{2\beta-2\delta} 
\right)\frac{dt}t 
 & \leq \\
\displaystyle C2^{-k(2\beta-2\delta)} \int_0^r \left\Vert 
A_t\right\Vert_{L^2(\Lambda T^{\ast}M)}^2 \frac{dt}t \leq C2^{-k(2\beta-2\delta)} V^{-1}(B) \leq 
C2^{-k(2\beta-2\delta-\D)}V^{-1}(2^{k+1}B)\ ,
\end{array}
\]
using Lemma \ref{offdiag4}.
We now need the fact that $2\beta-2\delta-{\D}>0$ to complete 
the proof of (\ref{atom}).

\medskip

\noindent{\bf Step 3: conclusion of the proof when $p=1$. } Consider again $\psi$ and $\widetilde{\psi}$ as in assertion $(a)$ of Theorem \ref{projector}. Observe first that the extension of $Q_f$ to a $T^{1,2}(\Lambda T^{\ast}M)$-bounded operator does not follow at once from (\ref{atom}). Indeed, up to this point, $Q_f$ is only defined on $T^{2,2}(\Lambda T^{\ast}M)$, and our task is to define it properly on $T^{1,2}(\Lambda T^{\ast}M)$. One way to do this could be to observe that, by Theorem \ref{atomic}, any element $F\in T^{1,2}(\Lambda T^{\ast}M)$ has an atomic decomposition $F=\sum_j \lambda_jA_j$, and to define $Q_f(F)=\sum_j \lambda_j Q_f(A_j)$ (which converges in $T^{1,2}(\Lambda T^{\ast}M)$), but we should then check that this definition does not depend on the decomposition of $F$ (which is not unique). Here, we argue differently. Since, by Proposition \ref{density}, $T^{1,2}(\Lambda T^{\ast}M)\cap T^{2,2}(\Lambda T^{\ast}M)$ is dense in $T^{1,2}(\Lambda T^{\ast}M)$, it is enough to show that there exists $C>0$ such that, for all $F\in T^{1,2}(\Lambda T^{\ast}M)\cap T^{2,2}(\Lambda T^{\ast}M)$, (\ref{eqprojector}
 ) holds for $F$ with $p=1$. \par
Consider such an $F$, and write $F=\sum_j \lambda_jA_j$ where $\sum_j \left\vert \lambda_j\right\vert\sim \left\Vert F\right\Vert_{T^{1,2}(\Lambda T^{\ast}M)}$ and, for each $j\geq 1$, $A_j$ is a $T^{1,2}(\Lambda T^{\ast}M)$ atom supported in $B_j\times [0,r_j]$ ($r_j$ denotes the radius of $B_j$). By Remark \ref{truncation}, if $x_0\in M$ and $\chi_n=\chi_{B(x_0,n)\times \left(\frac 1n,n\right)}$, then $F_n:=\chi_nF$ converges to $F$ both in $T^{1,2}(\Lambda T^{\ast}M)$ and in $T^{2,2}(\Lambda T^{\ast}M)$. For all $n\geq 1$, $F_n$ has an atomic decomposition in $T^{1,2}(\Lambda T^{\ast}M)$:
\begin{equation} \label{fnatom}
F_n=\sum_j \lambda_j(\chi_nA_j),
\end{equation}
where, for each $j$, $\chi_nA_j$ is a $T^{1,2}(\Lambda T^{\ast}M)$ atom. In particular, The series in (\ref{fnatom}) clearly converges in $T^{1,2}(\Lambda T^{\ast}M)$, but we claim that it also converges in $T^{2,2}(\Lambda T^{\ast}M)$. This relies on the following observation:
\begin{fact} \label{fac1}
For all $n\geq 1$, there exists $\kappa_n>0$ such that, for all $j\geq 1$, if $V(B_j)\leq \kappa_n$, then $\left(B(x_0,n)\times \left(\frac 1n,n\right)\right)\cap \left(B_j\times \left[0,r_j\right]\right)=\emptyset$.
\end{fact}
\noindent{\bf Proof of the fact: }We claim that
\[
\kappa_n=\frac{V(x_0,n)}{C(1+4n^2)^{\D}},
\]
where $C$ and $\D$ appear in (\ref{Dbis}), does the job. Indeed, assume now that $V(B_j)\leq \kappa_n$. If $B(x_0,n)\cap B_j=\emptyset$, there is nothing to do. Otherwise, let $y\in B(x_0,n)\cap B_j$, and write $B_j=B(x_j,r_j)$. The doubling property yields
\[
\begin{array}{lll}
V(x_0,n)& \leq & V(x_j,n+d(x_0,y)+d(y,x_j))\\
& \leq & V(x_j,2n+r_j)\\
& \leq & \displaystyle CV(B_j)\left(1+\frac{2n}{r_j}\right)^{\D}.
\end{array}
\]
Since $V(B_j)\leq \kappa_n$, it follows at once that $r_j\leq \frac{1}{2n}$, which obviously implies the desired conclusion. \hfill\fin\par
This fact easily implies that the series in (\ref{fnatom}) converges in $T^{2,2}(\Lambda T^{\ast}M)$. Indeed, we can drop in this series all the $j$'s such that $V(B_j)\leq \kappa_n$, and, if $V(B_j)>\kappa_n$,
\[
\left\Vert \chi_nA_j\right\Vert_{T^{2,2}(\Lambda T^{\ast}M)}\leq V(B_j)^{-1/2}\leq \kappa_n^{-1/2},
\]
which proves the convergence (remember that $\sum\left\vert \lambda_j\right\vert<+\infty$). \par
As a consequence, 
\begin{equation} \label{qffn}
Q_f(F_n)=\sum_j \lambda_j Q_f(\chi_nA_j),
\end{equation}
and this series converges in $T^{2,2}(\Lambda T^{\ast}M)$. But, since $\left\Vert Q_f(\chi_nA_j)\right\Vert_{T^{1,2}(\Lambda T^{\ast}M)}\leq C$ for all $j\geq 1$, the series in the right-hand side of (\ref{qffn}) also converges in $T^{1,2}(\Lambda T^{\ast}M)$ to some $G\in T^{1,2}(\Lambda T^{\ast}M)$, and, according to Lemma \ref{t1t2conv}, $G=Q_f(F_n)$. Therefore,
\begin{equation} \label{projectorfn}
\left\Vert Q_f(F_n)\right\Vert_{T^{1,2}(\Lambda T^{\ast}M)}\leq \sum_j \left\vert \lambda_j\right\vert \left\Vert Q_f(\chi_nA_j)\right\Vert_{T^{1,2}(\Lambda T^{\ast}M)}\leq C\left\Vert F\right\Vert_{T^{1,2}(\Lambda T^{\ast}M)}.
\end{equation}
Let us now prove that $(Q_f(F_n ))_{n\geq 1}$ is a Cauchy sequence in $T^{1,2}(\Lambda T^{\ast}M)$. From (\ref{qffn}), one has
\[
Q_f(F_n-F_m)=\sum_j \lambda_j Q_f((\chi_n-\chi_m)A_j),
\]
where the series converges in $T^{1,2}(\Lambda T^{\ast}M)$. Thus,
\[
\left\Vert Q_f(F_n-F_m)\right\Vert_{T^{1,2}(\Lambda T^{\ast}M)}\leq C\sum_j \left\vert \lambda_j\right\vert  \left\Vert Q_f((\chi_n-\chi_m)A_j)\right\Vert_{T^{1,2}(\Lambda T^{\ast}M)}.
\]
Fix now $\varepsilon>0$. There exists $J\geq 2$ such that $\sum\limits_{j\geq J} \left\vert \lambda_j\right\vert<\varepsilon$. Moreover, for each $1\leq j\leq J-1$, there exists $N_j\geq 1$ such that, for all $n\geq N_j$, $\chi_nA_j=A_j$. Therefore, there exists $N\geq 1$ such that, for all $n,m\geq N$ and all $1\leq j\leq J-1$, $(\chi_n-\chi_m)A_j=0$. As a consequence,
\[
\left\Vert Q_f(F_n-F_m)\right\Vert_{T^{1,2}(\Lambda T^{\ast}M)}\leq C\varepsilon
\]
for all $n,m\geq N$. Since $(Q_f(F_n ))_{n\geq 1}$ is a Cauchy sequence in $T^{1,2}(\Lambda T^{\ast}M)$, $Q_f(F_n)\rightarrow U$ in $T^{1,2}(\Lambda T^{\ast}M)$ for some $U\in T^{1,2}(\Lambda T^{\ast}M)$. Moreover, since $F_n\rightarrow F$ in $T^{2,2}(\Lambda T^{\ast}M)$, $Q_f(F_n)\rightarrow Q_f(F)$ in $T^{2,2}(\Lambda T^{\ast}M)$, and a new application of Lemma \ref{t1t2conv} yields $Q_f(F)=U$. It follows that $Q_f(F_n)\rightarrow Q_f(F)$ in $T^{1,2}(\Lambda T^{\ast}M)$. Letting $n$ go to $+\infty$ in (\ref{projectorfn}) gives the desired result, which ends up the proof of the case $p=1$.

\medskip

\noindent{\bf Step 4: End of the proof. } Using the interpolation results for tent spaces (Theorem 
\ref{interpolation}), we obtain the conclusion of Theorem 
\ref{projector} for all $1\leq p\leq 2$. Finally, the duality for 
tent spaces (Theorem \ref{duality}) also yields this conclusion for 
$p\geq 
2$ (note that the assumptions on $\psi$ and $\widetilde{\psi}$ have been switched), which ends the 
proof of Theorem \ref{projector}. \hfill\fin 

\bigskip

\textbf{From now on, we constantly assume the doubling property \eqref{D}.}

\SE{Definition of Hardy spaces and first results} \label{def}

\subsection{Definition and first properties of Hardy spaces} 
\label{firstdef}

We are now able to give the definition of the $H^p(\Lambda T^{\ast}M)$ space for all 
$1\leq p\leq +\infty$, $p\neq 2$, by means of quadratic functionals (as was done for $H^2(\Lambda T^{\ast}M)$ in Section \ref{h2}) and tent spaces. Theorem \ref{projector} tells us that we have to distinguish between the cases $1\leq p<2$ and $2<p\leq +\infty$. 

Given $\psi\in \Psi(\Sigma_{\theta}^0)$ for some $\theta>0$, set 
$\psi_{t}(z)=\psi(tz)$ for all $t>0$ and all $z\in 
\Sigma_{\theta}^{0}$. Recall that the operator ${\mathcal 
S}_\psi:T^{2,2}(\Lambda T^{\ast}M)\longrightarrow 
L^2(\Lambda T^{\ast}M)$ is defined by 
\[
{\mathcal S}_{\psi}H=\int_0^{+\infty} \psi_t(D)H_{t}\frac{dt}t\]
and ${\mathcal Q}_{\psi}:L^2(\Lambda T^{\ast}M)\longrightarrow 
T^{2,2}(\Lambda T^{\ast}M)$ by
\[
\left({\mathcal Q}_{\psi}h\right)_{t}=\psi_t(D)h
\]
for all $h\in L^2(\Lambda T^{\ast}M)$ and all $t>0$.

\bigskip

\begin{defi} \label{hardyleq2} For each $\psi \in \Psi(\Sigma_{\theta}^0)$, define
$E^p_{D,\psi}(\Lambda T^{\ast}M)={\mathcal S}_{\psi}(T^{p,2}(\Lambda T^{\ast}M)\cap 
T^{2,2}(\Lambda T^{\ast}M))$ with semi-norm
\[
\left\Vert h\right\Vert_{H^p_{D,\psi}(\Lambda T^{\ast}M)}=\inf\{\left\Vert 
H\right\Vert_{T^{p,2}(\Lambda T^{\ast}M)}; H\in T^{p,2}(\Lambda T^{\ast}M)\cap 
T^{2,2}(\Lambda T^{\ast}M), {\mathcal S}_{\psi}H=h\}\ .
\]
\end{defi}
{\bf The case when $1\leq p<2$: } Recall that $\beta=\left[\frac\kappa 2\right]+1$ (the smallest integer larger than $\frac\kappa 2$). It turns out that, provided $\psi\in\Psi_{\beta,2}(\Sigma_{\theta}^0)$, $E^p_{D,\psi}(\Lambda T^{\ast}M)$ is actually independent from the choice of $\psi$, and can be described by means of the operators ${\mathcal Q}_{\widetilde{\widetilde{\psi}}}$ if $\widetilde{\widetilde{\psi}}\in \Psi_{1,\beta+1}(\Sigma_{\theta}^0)$ (see Section \ref{h2}) : 
\begin{lem} \label{equivalent} If $\psi,\widetilde\psi\in\Psi_{\beta,2}(\Sigma_{\theta}^0)$ and $\widetilde{\widetilde\psi}\in\Psi_{1,\beta+1}(\Sigma_{\theta}^0)$,
 then $$E^p_{D,\psi}(\Lambda T^{\ast}M)=E^p_{D,\widetilde\psi}(\Lambda T^{\ast}M)=\{h\in H^2(\Lambda T^{\ast}M)\ ;\ \|\mathcal Q_{\widetilde{\widetilde\psi}}h\|_{T^{p,2}(\Lambda T^{\ast}M)}<\infty\}$$ with norm $$\left\Vert h\right\Vert_{H^p_{D,\psi}(\Lambda T^{\ast}M)}\sim\left\Vert h\right\Vert_{H^p_{D,\widetilde\psi}(\Lambda T^{\ast}M)}\sim
 \|\mathcal Q_{\widetilde{\widetilde\psi}}h\|_{T^{p,2}(\Lambda T^{\ast}M)}\ . $$\end{lem}

\noindent{\bf Proof: } Fix $\psi,\widetilde{\psi}$ and $\widetilde{\widetilde{\psi}}$ as in Lemma \ref{equivalent}. Observe first that, if $\varphi\in 
\Psi_{1,\beta+1}(\Sigma_{\theta}^0)$, then Theorem \ref{projector} tells us that 
$\mathcal Q_{\varphi} \mathcal S_{\psi}$ extends to a bounded operator in $T^{p,2}(\Lambda T^{\ast}M)$.\par

\noindent{\bf $(a)$} First, let $h\in E^p_{D,\psi}(\Lambda T^{\ast}M)$ and $\varphi\in 
\Psi_{1,\beta+1}(\Sigma_{\theta}^0)$. There exists $H\in T^{p,2}(\Lambda T^{\ast}M)$ 
with $\left\Vert H\right\Vert_{T^{p,2}(\Lambda T^{\ast}M)}\leq 2\left\Vert 
h\right\Vert_{H^p_{D,\psi}(\Lambda T^{\ast}M)}$ such that  $h={\mathcal S}_{\psi}H$, so that $h\in 
H^{2}(\Lambda T^{\ast}M)$ and, because of our observation, 
$\left\Vert {\mathcal Q}_{\varphi}h\right\Vert_{T^{p,2}(\Lambda T^{\ast}M)}\leq 
C\left\Vert H\right\Vert_{T^{p,2}(\Lambda T^{\ast}M)}\leq C\left\Vert 
h\right\Vert_{H^p_{D,\psi}(\Lambda T^{\ast}M)}$. In particular, $E^p_{D,\psi}(\Lambda T^{\ast}M)\subset \{h\in H^2(\Lambda T^{\ast}M)\ ;\ \|\mathcal Q_{\widetilde{\widetilde\psi}}h\|_{T^{p,2}(\Lambda T^{\ast}M)}<\infty\}$.\par
\noindent{\bf $(b)$} Assume now that $h\in H^2(\Lambda T^{\ast}M)$ and ${\mathcal 
Q}_{\widetilde{\widetilde\psi}}h\in T^{p,2}(\Lambda T^{\ast}M)$. We claim that there 
exists $\zeta\in \Psi_{\beta,2}(\Sigma_{\theta}^0)$ such 
that $h\in E^p_{D,\zeta}(\Lambda T^{\ast}M)$. Indeed, by Remark 
\ref{care}, there exists $\zeta\in \Psi_{\beta,2}(\Sigma_{\theta}^0)$ 
such that ${\mathcal S}_{\zeta}{\mathcal 
Q}_{\widetilde{\widetilde{\psi}}}=\mbox{Id}$ on $H^{2}(\Lambda T^{\ast}M)$. Therefore, if 
$H=Q_{\widetilde{\widetilde\psi}}h$, one has $h={\mathcal 
S}_{\zeta}H$, which shows that $\left\Vert 
h\right\Vert_{H^p_{D,\zeta}(\Lambda T^{\ast}M)}\leq \left\Vert 
Q_{\widetilde{\widetilde\psi}}h\right\Vert_{T^{p,2}(\Lambda T^{\ast}M)}$.\par
\noindent{\bf $(c)$} We check now that, if $\zeta$ is as in step $(b)$ and $h\in E^p_{D,\zeta}(\Lambda T^{\ast}M)$, 
then $h\in E^p_{D,\psi}(\Lambda T^{\ast}M)$. Indeed, thanks to Remark \ref{care} again, 
there exists $\varphi\in \Psi_{1,\beta+1}(\Sigma_{\theta}^0)$ such 
that ${\mathcal S}_{\psi}{\mathcal Q}_{\varphi}=\mbox{Id}$ on 
$H^2(\Lambda T^{\ast}M)$. According to $(a)$, $\left\Vert {\mathcal Q}_{\varphi}h\right\Vert_{T^{p,2}(\Lambda T^{\ast}M)}\leq C\left\Vert h\right\Vert_{H^p_{D,\zeta}(\Lambda T^{\ast}M)}$. Since $h\in H^2(\Lambda T^{\ast}M)$, one has $h={\mathcal S}_{\psi} {\mathcal Q}_{\varphi}h$, which shows that $\left\Vert h\right\Vert_{H^p_{D,\psi}(\Lambda T^{\ast}M)}\leq C\left\Vert h\right\Vert_{H^p_{D,\zeta}(\Lambda T^{\ast}M)}$. 

\noindent {\bf $(d)$} It remains to be shown that $\left\Vert h\right\Vert_{H^p_{D,\psi}(\Lambda T^{\ast}M)}$ is a norm rather than a seminorm on $E^p_{D,\psi}(\Lambda T^{\ast}M)$. Let $h\in
E^p_{D,\psi}(\Lambda T^{\ast}M)$ with $\left\Vert h\right\Vert_{H^p_{D,\psi}(\Lambda T^{\ast}M)}=\|\mathcal Q_{\widetilde{\widetilde\psi}}h\|_{T^{p,2}(\Lambda T^{\ast}M)}=0$.
Then $h\in H^2(\Lambda T^{\ast}M)\cap \mathcal N(\mathcal Q_{\widetilde{\widetilde\psi}})=\overline{\mathcal R(D)}\cap\mathcal N(D)=\{0\}$. {\it i.e.} $h=0$ as required.
\hfill\fin

\begin{rem} It follows that these spaces and maps are independent of $\theta\in (0,\tfrac\pi2)$ too.
\end{rem}

\begin{pro} With the notation of Lemma \ref{equivalent}, $\{h\in \mathcal R(D)\ ;\ \|\mathcal Q_{\widetilde{\widetilde\psi}}h\|_{T^{p,2}(\Lambda T^{\ast}M)}<\infty\}$ is dense in $E^p_{D,\psi}(\Lambda T^{\ast}M)$ for all $\widetilde{\widetilde\psi}\in\Psi_{1,\beta+1}(\Sigma_{\theta}^0)$.
\end{pro}

\noindent{\bf Proof: } Fix $\widetilde{\widetilde\psi}\in\Psi_{1,\beta+1}(\Sigma_{\theta}^0)$ and choose $\widetilde{\psi} \in \Psi_{\beta,2}(\Sigma_{\theta}^0)$ such that $\mathcal S_{\widetilde\psi}\mathcal Q_{\widetilde{\widetilde\psi}}h = h$ for all $h\in H^2(\Lambda T^{\ast}M)$. For a given $h\in E^p_{D,\psi}(\Lambda T^{\ast}M)$, set  $H= \mathcal Q_{\widetilde{\widetilde\psi}}h\in T^{p,2}(\Lambda T^{\ast}M)\cap T^{2,2}(\Lambda T^{\ast}M)$, and define, for each natural number $N$,  
$H_N\in T^{p,2}(\Lambda T^{\ast}M)\cap T^{2,2}(\Lambda T^{\ast}M)$ by $H_N(x,t)=H(x,t)\chi_{[\frac1N,N]}(t)$. It is not difficult to show that $H_N\to H$ in $T^{p,2}(\Lambda T^{\ast}M)$, and so $h_N:= \mathcal S_{\widetilde\psi} H_N \to h$ in 
$E^p_{D,\psi}(\Lambda T^{\ast}M)=E^p_{D,\widetilde{\psi}}(\Lambda T^{\ast}M)$. 

It remains to be shown that $h_N\in\mathcal R(D)$. This holds because
\begin{equation*}
h_N=\int_{\frac1N }^{N} \widetilde\psi_t(D)H_{t}\frac{dt}t=D\int_{\frac1N} ^{N} \phi(tD)H_{t}\ dt
\end{equation*}
where $\phi\in H^\infty(S^0_\mu)$ is defined by $\phi(z)=\frac1z\widetilde\psi(z)$.
\hfill\fin
\bigskip

We are now in a position to define the {\it Hardy spaces} associated with $D$.

\begin{defi} Suppose $1\leq p<2$. Define $H^p_D(\Lambda T^{\ast}M)$ to be the completion of $E^p_{D,\psi}(\Lambda T^{\ast}M)$ under any of the equivalent norms $\left\Vert h\right\Vert_{H^p_{D,\psi}(\Lambda T^{\ast}M)}$ with $\psi\in\Psi_{\beta,2}(\Sigma^0_\theta)$, which we write as just $\left\Vert h\right\Vert_{H^p_{D}(\Lambda T^{\ast}M)}$.
\end{defi}

In particular, $$H^p_D(\Lambda T^{\ast}M)=\overline{\{h\in \mathcal R(D)\ ;\ \|\mathcal Q_{\psi}h\|_{T^{p,2}(\Lambda T^{\ast}M)}<\infty\}}$$ under the norm 
 $\|\mathcal Q_\psi h\|_{T^{p,2}(\Lambda T^{\ast}M)}$ for any $\psi\in \Psi_{1,\beta+1}(\Sigma^0_\theta)$. For example,
\[
\begin{array}{lll}
\displaystyle \left\Vert h\right\Vert_{H^p_{D}(\Lambda T^{\ast}M)} & \sim & \displaystyle  \left\Vert tD e^{-t\sqrt{\Delta}}h\right\Vert_{T^{p,2}(\Lambda T^{\ast}M)}\\
&  \sim & \displaystyle \left\Vert t^2\Delta e^{-t^2\Delta}h\right\Vert_{T^{p,2}(\Lambda T^{\ast}M)}\\
& \sim &\displaystyle \left\Vert tD (I+t^2\Delta)^{-N} h\right\Vert_{T^{p,2}(\Lambda T^{\ast}M)}
\end{array}
\]
where $N\geq \tfrac \beta2+1$.
\bigskip

\noindent {\bf The case when $2<p<\infty$}: The same procedure works, but with the 
roles of $\Psi_{\beta,2}(\Sigma_{\theta}^0)$ and 
$\Psi_{1,\beta+1}(\Sigma_{\theta}^0)$ interchanged.

\begin{defi} Suppose $2< p<\infty$. Define $H^p_D(\Lambda T^{\ast}M)$ to be the completion of $E^p_{D,\psi}(\Lambda T^{\ast}M)$ under any of the equivalent norms $\left\Vert h\right\Vert_{H^p_{D,\psi}(\Lambda T^{\ast}M)}$ with $\psi\in\Psi_{1,\beta+1}(\Sigma^0_\theta)$, which we write as just $\left\Vert h\right\Vert_{H^p_{D}(\Lambda T^{\ast}M)}$.
\end{defi}

In particular, $$H^p_D(\Lambda T^{\ast}M)=\overline{\{h\in \mathcal R(D)\ ;\ \|\mathcal Q_{\psi}h\|_{T^{p,2}(\Lambda T^{\ast}M)}<\infty\}}$$ under the norm  $\|\mathcal Q_\psi h\|_{T^{p,2}(\Lambda T^{\ast}M)}$ for any $\psi\in \Psi_{\beta,2}(\Sigma^0_\theta)$. For example,
$$
\begin{array}{lll}
\displaystyle \left\Vert h\right\Vert_{H^p_{D}(\Lambda T^{\ast}M)}& \sim &
\displaystyle \left\Vert (tD)^{\beta} e^{-t\sqrt{\Delta}}h\right\Vert_{T^{p,2}(\Lambda T^{\ast}M)}\\
& \sim&
\displaystyle \left\Vert (t^2\Delta)^{M} e^{-t^2\Delta}h\right\Vert_{T^{p,2}(\Lambda T^{\ast}M)}\\
& \sim&
\displaystyle \left\Vert (tD)^\beta (I+t^2\Delta)^{-N} h\right\Vert_{T^{p,2}(\Lambda T^{\ast}M)}
\end{array}
$$
where $M\geq \tfrac \beta2$ and $N\geq \tfrac \beta2+1$.
\bigskip

Suppose that the function $\psi$ used in any of the above norms is an 
even function. Then $\psi_t(D)=\widetilde\psi(t^2\Delta)$, where 
$\widetilde\psi\in\Psi(\Sigma^0_{2\theta^{+}})$. We thus see that we have defined Hardy spaces 
$H^p_\Delta(\Lambda T^{\ast}M)$ corresponding to the Laplacian $\Delta$, and that they 
are the same as the spaces $H^p_D(\Lambda T^{\ast}M)$.  From now on, for all $1\leq p<+\infty$, the 
$H^p_{D}(\Lambda T^{\ast}M)$ space, which 
coincides with $H^p_{\Delta}(\Lambda T^{\ast}M)$, will be denoted by $H^p(\Lambda T^{\ast}M)$.

\medskip
We define $H^{\infty}(\Lambda T^{\ast}M)$ in a different way. This definition relies on the following lemma:
\begin{lem} \label{dualh1}
Let $\psi\in \Psi_{1,\beta+1}(\Sigma_{\theta}^0)$.
\begin{itemize}
\item[$(a)$]
Let $G\in T^{\infty,2}(\Lambda T^{\ast}M)$. Then the map $T_G$, initially defined on $E^1_D(\Lambda T^{\ast}M)$ by $T_G(f)= \iint \langle \left({\mathcal Q}_{\psi}f\right)_t(x),G(x,t)\rangle dx\frac{dt}t$, extends in a unique way to a bounded linear functional on $H^1(\Lambda T^{\ast}M)$, denoted again by $T_G$.
\item[$(b)$]
Conversely, if $U$ is a bounded linear functional on $H^1(\Lambda T^{\ast}M)$, there exists $G\in T^{\infty,2}(\Lambda T^{\ast}M)$ such that $U=T_G$.
\end{itemize}
\end{lem}
The proof is an immediate consequence of assertion $(b)$ in Theorem \ref{duality} and the definition of $E^1_{D,\psi}(\Lambda T^{\ast}M)$. We define $H^{\infty}(\Lambda T^{\ast}M)$ as the dual space of $H^1(\Lambda T^{\ast}M)$, equipped with the usual dual norm. Observe that, by Lemma \ref{dualh1}, one has
\[
\left\Vert U\right\Vert_{H^{\infty}(\Lambda T^{\ast}M)}\sim \inf\left\{\left\Vert G\right\Vert_{T^{\infty,2}(\Lambda T^{\ast}M)};\ U=T_G\right\}.
\]

Theorems \ref{duality} and \ref{interpolation} yield duality and 
interpolation results for Hardy spaces:

\begin{theo} \label{dualityhardy}
The pairing $\langle g,h\rangle\mapsto \int_{M} \langle g(x),h(x)\rangle dx$ realizes 
$H^{p^{\prime}}(\Lambda T^{\ast}M)$ as equivalent with the dual of 
$H^{p}(\Lambda T^{\ast}M)$ if 
$1<p<+\infty$ and $1/p+1/p^{\prime}=1$. Moreover, by definition, the dual of $H^1(\Lambda T^{\ast}M)$ is $H^{\infty}(\Lambda T^{\ast}M)$.
\end{theo}    
\begin{theo} \label{interpohardy}
Let $1\leq p_{0}<p<p_{1}\leq +\infty$ and $\theta\in 
\left(0,1\right)$ such that $1/p=(1-\theta)/p_{0}+\theta/p_{1}$. 
Then 
$\left[H^{p_{0}}(\Lambda T^{\ast}M),H^{p_{1}}(\Lambda T^{\ast}M)\right]_{\theta}=H^{p}(\Lambda T^{\ast}M)$. 
\end{theo}
\noindent{\bf Proof: } The spaces
$\mathcal P_\psi T^{p,2}(\Lambda T^{\ast}M)$ defined in Remark \ref{interp} interpolate by the complex method for $1\leq p\leq\infty$, where we have taken 
$\psi\in \Psi_{1,\beta+1}(\Sigma_{\theta}^0)\cap
\Psi_{\beta,2}(\Sigma_{\theta}^0)$ with $\int_0^\infty\psi(\pm t)\psi(\pm t)\frac{dt}t =1$ and defined the projection
$\mathcal P_{\psi}:={\mathcal Q}_{\psi}{\mathcal S}_{\psi}$. It is straightforward to see that the map ${\mathcal Q}_{\psi}$ extends to an isomorphism from 
$\mathcal P_\psi T^{p,2}(\Lambda T^{\ast}M)$ to $H^p(\Lambda T^{\ast}M)$ with inverse ${\mathcal S}_{\psi}$ for each $p$, and that these maps coincide for different values of $p$. The result follows.
\hfill\fin

\begin{rem} \label{BMO}
Since $H^{\infty}(\Lambda T^{\ast}M)$ is the dual space of $H^1(\Lambda T^{\ast}M)$, it turns out that $H^{\infty}(\Lambda T^{\ast}M)$ is actually a $BMO$-type space. Recall that $BMO(\R^n)$ is the dual space of $H^1(\R^n)$ (\cite{feffermanstein}) and that similar duality results have been established for other kinds of Hardy spaces, in particular in \cite{loumcintosh} for Hardy spaces of exact differential forms in $\R^n$. To keep homogeneous notations and simplify our previous and foregoing statements about Hardy spaces, we write $H^{\infty}(\Lambda T^{\ast}M)$ instead of $BMO(\Lambda T^{\ast}M).$ \end{rem}

\subsection{Riesz transform and Functional calculus}

We are now ready to prove the first part of Corollary \ref{riesz}, namely 

\begin{theo} \label{rieszbound}
For all $1\leq p\leq +\infty$, the Riesz transform $D\Delta^{-1/2}$, initially defined on ${\mathcal R}(\Delta)$, extends to a $H^p(\Lambda T^{\ast}M)$-bounded operator. More precisely, one has $\left\Vert 
D\Delta^{-1/2}h\right\Vert_{H^p(\Lambda T^{\ast}M)}\sim \left\Vert 
h\right\Vert_{H^p(\Lambda T^{\ast}M)}$.
\end{theo}
\noindent{\bf Proof: } The case when $p=2$ is in Section \ref{h2}. Consider now the case when $1\leq p<+\infty$ and $p\neq 2$. Choose $\psi\in \Psi_{1,\beta+1}(\Sigma_{\theta}^0)$ 
when $1\leq p<2$, and $\psi\in \Psi_{\beta,2}(\Sigma_{\theta}^0)$ 
when $2<p<\infty$. 
In either case the holomorphic function $\widetilde\psi$ defined by $\widetilde\psi(z) 
=\text{sgn(\text{Re}\,z)}\psi(z)$ belongs to the same space, and 
moreover 
$\widetilde\psi(D)=D\Delta^{-1/2}\psi(D)$. Hence, by Lemma \ref{equivalent},
\begin{equation*} \left\Vert 
D\Delta^{-1/2}h\right\Vert_{H^p(\Lambda T^{\ast}M)}\sim 
\left\Vert 
D\Delta^{-1/2}\psi(D)h\right\Vert_{T^{p,2}(\Lambda T^{\ast}M)}=\left\Vert\widetilde\psi(D)h\right\Vert_{T^{p,2}(\Lambda T^{\ast}M)}\sim 
\left\Vert h\right\Vert_{H^p(\Lambda T^{\ast}M)}
\end{equation*}
for all $h\in H^p(\Lambda T^{\ast}M)$. The case $p=+\infty$ follows from the case $p=1$ by duality. \hfill\fin 

\medskip
Similar estimates actually give the following more general result on the holomorphic functional calculus of $D$:
\begin{theo} \label{hardycalculus}
For all $1\leq p\leq +\infty$, $f(D)$ is $H^p(\Lambda T^{\ast}M)$-bounded for all 
$f\in H^\infty(\Sigma^0_{\theta})$ with
$\left\Vert f(D)h\right\Vert_{H^p(\Lambda T^{\ast}M)}\leq 
C\left\Vert f\right\Vert_{\infty}\left\Vert 
h\right\Vert_{H^p(\Lambda T^{\ast}M)}$\ .
\end{theo}

When $1\leq p<\infty$, this estimate follows from Theorem \ref{projector} and the definitions of $H^p(\Lambda T^{\ast}M)$.  When $p=\infty$, use duality.

\bigskip

Let us finish this section by discussing the boundedness of the Hodge-Riesz transforms.  Let $1\leq p\leq +\infty$, and denote by $n$ the dimension of $M$. First, 
the splitting $\Lambda T^{\ast}M=\oplus_{0\leq k\leq n} 
\Lambda^k T^{\ast}M$ allows us to define naturally $H^p(\Lambda^k 
T^{\ast}M)$ for all $0\leq k\leq n$ (first for $1\leq p<+\infty$, then for $p=+\infty$ by duality), and one has, if 
$f=(f_{0},\ldots,f_{n})\in \Lambda T^{\ast}M$,
\begin{equation} \label{splitting}
\left\Vert f\right\Vert_{H^p(\Lambda T^{\ast}M)}\sim\sum_{k=0}^n \left\Vert 
f_{k}\right\Vert_{H^p(\Lambda^kT^{\ast}M)}.
\end{equation} 
To see this when $1\leq p<\infty$, recall that  $H^p(\Lambda^k 
T^{\ast}M) = H^p_\Delta(\Lambda^k 
T^{\ast}M)$, and note that $\Delta$ preserves the decomposition into $k$-forms. Specializing Theorem \ref{rieszbound} to $k$ forms implies that, for all $0\leq k\leq n-1$, 
$d\Delta^{-1/2}$ is 
$H^p(\Lambda^kT^{\ast}M)-H^p(\Lambda^{k+1}T^{\ast}M)$ bounded, and 
that, for all $1\leq k\leq n$, $d^{\ast}\Delta^{-1/2}$ is 
$H^p(\Lambda^kT^{\ast}M)-H^p(\Lambda^{k-1}T^{\ast}M)$ bounded. 
Using (\ref{splitting})  we have obtained:

\begin{theo} \label{rieszboundbis}
For all $1\leq p\leq +\infty$, $d\Delta^{-1/2}$ and 
$d^{\ast}\Delta^{-1/2}$ are both $H^p(\Lambda T^{\ast}M)$ bounded.
\end{theo}
 
\subsection{The Hodge decomposition for $H^p(\Lambda T^{\ast}M)$} \label{hodgedecompo}

We can define other Hardy spaces, associated 
to the operators $d$ and $d^{\ast}$, which leads us to a Hodge decomposition for $H^p(\Lambda T^{\ast}M)$. 
Recall from Section \ref{h2} that 
\[
H^2(\Lambda T^{\ast}M)=H^2_d(\Lambda T^{\ast}M)\oplus H^2_{d^{\ast}}(\Lambda T^{\ast}M),
\]
where $H^2_d(\Lambda T^{\ast}M)=\overline{{\mathcal R}(d)}$ and $H^2_{d^{\ast}}(\Lambda T^{\ast}M)=\overline{{\mathcal R}(d^{\ast})}$ and that the orthogonal projections are given by $dD^{-1}$ and $d^{\ast}D^{-1}$.

For $1\leq p<+\infty$ and $p\neq 2$, set
\[
H^p_d(\Lambda T^{\ast}M)=\overline{{\mathcal R}(d)\cap H^p(\Lambda T^{\ast}M)},\ H^p_{d^{\ast}}(\Lambda T^{\ast}M)=\overline{{\mathcal R}(d^{\ast})\cap H^p(\Lambda T^{\ast}M)}
\]
where the closure is taken in the $H^p(\Lambda T^{\ast}M)$ topology. We have the following Hodge decomposition for $H^p(\Lambda T^{\ast}M)$:
\begin{theo} \label{Hodge}
For all $1\leq p<+\infty$, one has $H^p(\Lambda T^{\ast}M)=H^p_d(\Lambda T^{\ast}M)\oplus H^p_{d^{\ast}}(\Lambda T^{\ast}M)$, and the sum is topological.
\end{theo}
\noindent{\bf Proof: } The orthogonal projection $dD^{-1}$ from $H^2(\Lambda T^{\ast}M)$ to $H^2_d(\Lambda T^{\ast}M)$ defines a bounded operator from $H^p(\Lambda T^{\ast}M)$ to $H^p_d(\Lambda T^{\ast}M)$. Indeed, $dD^{-1}=d\Delta^{-1/2}D\Delta^{-1/2}$, $D\Delta^{-1/2}$ is $H^p(\Lambda T^{\ast}M)$ bounded by Theorem \ref{rieszbound} and $d\Delta^{-1/2}$ is $H^p(\Lambda T^{\ast}M)-H^p_d(\Lambda T^{\ast}M)$ bounded by Theorem \ref{rieszboundbis}. Similarly, $d^{\ast}D^{-1}$ is a bounded operator from $H^p(\Lambda T^{\ast}M)$ to $H^p_{d^{\ast}}(\Lambda T^{\ast}M)$. Since, for all $f\in H^p(\Lambda T^{\ast}M)$, one has $f=dD^{-1}f+d^{\ast}D^{-1}f$, the theorem is proved. \hfill\fin

Note that, for all $1\leq p<+\infty$, $H^p_d(\Lambda T^{\ast}M)$ and $H^p_{d^{\ast}}(\Lambda T^{\ast}M)$ can also be described 
by means of tent spaces
in the same way as $H^p(\Lambda T^{\ast}M)$. More precisely, if $\psi\in 
\Psi_{1,\tau}(\Sigma^{0}_{\theta})$ for some $\tau, \theta>0$, define $\phi\in H^\infty(\Sigma^{0}_{\theta})$ by $\phi(z)=\frac1z\psi(z)$ and then define,  for $H\in 
T^{2,2}(\Lambda T^{\ast}M)$,
\[
{\mathcal S}_{d,\psi}(H)=\int_{0}^{+\infty} 
td\phi_{t}(D)H_{t}\frac{dt}t\mbox{ and }{\mathcal S}_{d^{\ast},\psi}(H)=\int_{0}^{+\infty} 
td^{\ast}\phi_{t}(D)H_{t}\frac{dt}t,
\]
and, for all $h\in L^{2}(\Lambda T^{\ast}M)$, define
\[
(\mathcal Q_{d,\psi}h)_{t}=td\phi_{t}(D)h\mbox{ and }(\mathcal 
Q_{d^{\ast},\psi}h)_{t}=td^{\ast}\phi_{t}(D)h.
\]
Then, for $1\leq p< +\infty$, replacing ${\mathcal S}_{\widetilde{\psi}}$ by 
${\mathcal S}_{d,\widetilde{\psi}}$ (resp. by ${\mathcal S}_{d^{\ast},{\widetilde{\psi}}}$) and 
${\mathcal Q}_{\psi}$ by ${\mathcal Q}_{d^{\ast},\psi}$ (resp. ${\mathcal 
Q}_{d,\psi}$) in Section \ref{firstdef}, one obtains a characterization of $H^p_{d}(\Lambda T^{\ast}M)$ (resp. 
$H^p_{d^{\ast}}(\Lambda T^{\ast}M)$) by means of ${\mathcal S}_{d,\widetilde{\psi}}$ and ${\mathcal 
Q}_{d^{\ast},\psi}$ (resp. ${\mathcal S}_{d^{\ast},\widetilde{\psi}}$ and ${\mathcal 
Q}_{d,\psi}$), provided that $\psi\in \Psi_{1,\beta+1}(\Sigma^0_{\theta})$ and $\widetilde{\psi}\in \Psi_{\beta,2}(\Sigma^0_{\theta})$ if $1\leq p< 2$, and $\psi\in \Psi_{\beta,2}(\Sigma^0_{\theta})$ and $\widetilde{\psi}\in \Psi_{1,\beta+1}(\Sigma^0_{\theta})$ if $2< p< \infty$. 
\bigskip

It is plain to observe that, if $1<p<+\infty$, the dual of $H^p_d(\Lambda T^{\ast}M)$ is isomorphic to $H^{p^{\prime}}_d(\Lambda T^{\ast}M)$, where $1/p+1/p^{\prime}=1$. We define $H^{\infty}_d(\Lambda T^{\ast}M)$ (resp. $H^{\infty}_{d^{\ast}}(\Lambda T^{\ast}M)$) as the dual space of $H^1_d(\Lambda T^{\ast}M)$ (resp. $H^1_{d^{\ast}}(\Lambda T^{\ast}M)$). Lemma \ref{dualh1} provides another description of $H^{\infty}_d(\Lambda T^{\ast}M)$ and $H^{\infty}_{d^{\ast}}(\Lambda T^{\ast}M)$. Namely, fix  $\widetilde{\psi}\in \Psi_{1,\beta+1}(\Sigma^0_{\theta})$. For all $G\in T^{\infty,2}(\Lambda T^{\ast}M)$ and all $f\in H^1(\Lambda T^{\ast}M)\cap {\mathcal R}(d)$, define
\[
T_{d^{\ast},G}(f)= \iint \langle (\mathcal 
Q_{d^{\ast},\psi}f)_{t}(x),G(x,t)\rangle dx\frac{dt}t.
\]
Then, $T_{d^{\ast},G}$ is a bounded linear functional on $H^1_d(\Lambda T^{\ast}M)$, and conversely, any bounded linear functional $U$ on $H^1_d(\Lambda T^{\ast}M)$ is equal to $T_{d^{\ast},G}$ for some $G\in T^{\infty,2}(\Lambda T^{\ast}M)$. Furthermore, $\left\Vert T_{d^{\ast},G}\right\Vert \sim \left\Vert G\right\Vert_{T^{\infty,2}(\Lambda T^{\ast}M)}$. The description of $(H^1_{d^{\ast}}(\Lambda T^{\ast}M))^{\prime}$ is similar.

\bigskip

As a consequence we have a more precise statement for the  Hodge-Riesz transforms.

\begin{theo} \label{rieszboundter}
Let $1\leq p\leq +\infty$. Then  $d\Delta^{-1/2}$ extends to a continuous isomorphism from 
$H^p_{d^{\ast}}(\Lambda T^{\ast}M)$ onto $H^p_{d}(\Lambda T^{\ast}M)$ and $d^{\ast}\Delta^{-1/2}$ to a continuous isomorphism from 
$H^p_{d}(\Lambda T^{\ast}M)$ onto $H^p_{d^{\ast}}(\Lambda T^{\ast}M)$. These operators are inverse to one another. 
\end{theo}
\noindent{\bf Proof: }Assume first that $1\leq p<+\infty$. Theorem \ref{rieszboundbis} shows that 
$d\Delta^{-1/2}$ extends to a $H^p(\Lambda T^{\ast}M)$-bounded linear map, and the very definition of 
$H^p_{d}(\Lambda T^{\ast}M)$ therefore ensures that it is $H^p_{d^{\ast}}(\Lambda T^{\ast}M)-H^p_{d}(\Lambda T^{\ast}M)$ bounded. 
Similarly, $d^{\ast}\Delta^{-1/2}$ extends to a $H^p_{d}(\Lambda T^{\ast}M)-H^p_{d^{\ast}}(\Lambda T^{\ast}M)$ 
bounded map. Next, for $f\in H^p_{d}(\Lambda T^{\ast}M)\cap{\mathcal R}(d)$, 
$(d\Delta^{-1/2})(d^{\ast}\Delta^{-1/2})f=dd^{\ast}\Delta^{-1}f=f$ since 
$d^{\ast}df=0$, which shows that $d\Delta^{-1/2}$ is onto $H^p_d(\Lambda T^{\ast}M)$ and $d^{\ast}\Delta^{-1/2}$ is one-to-one from $H^p_{d}(\Lambda T^{\ast}M)$. Symmetrically, $d^{\ast}\Delta^{-1/2}$ is onto $H^p_{d^{\ast}}(\Lambda T^{\ast}M)$ and $d\Delta^{-1/2}$ is one-to-one from $H^p_{d^{\ast}}(\Lambda T^{\ast}M)$. 

Finally, the conclusion for $p=+\infty$ follows from the case $p=1$ by 
duality. The proof is straightforward and relies on the fact that 
$d^{\ast}\Delta^{-1/2}=\Delta^{-1/2}d^{\ast}$ on a dense subspace of 
$H^{1}_{d}(\Lambda T^{\ast}M)$. \hfill\fin

\bigskip

To finich, let us specialize  the above to $k$-forms. For $1\leq p\leq +\infty$, one can also naturally define $H^p_{d}(\Lambda^kT^{\ast}M)$ for all 
$1\leq k\leq n$ and $H^p_{d^{\ast}}(\Lambda^kT^{\ast}M)$ for all 
$0\leq k\leq n-1$ (first for $1\leq p<+\infty$, then using duality for $p=+\infty$), and Theorem \ref{rieszboundter} shows the following result:
\begin{theo} \label{rieszboundfour}
Let $1\leq p\leq +\infty$.
\begin{itemize}
\item[$(a)$]    
For all $0\leq k\leq n-1$, 
$d\Delta^{-1/2}$ is a continuous isomorphism from 
$H^p_{d^{\ast}}(\Lambda^kT^{\ast}M)$ onto $H^p_{d}(\Lambda^{k+1}T^{\ast}M)$. 
\item[$(b)$]
For all $1\leq k\leq n$, $d^{\ast}\Delta^{-1/2}$ is a continuous isomorphism from
$H^p_{d}(\Lambda^kT^{\ast}M)$ onto $H^p_{d^{\ast}}(\Lambda^{k-1}T^{\ast}M)$.
\end{itemize}
\end{theo}

\SE{The decomposition into molecules} \label{mol}

As recalled in the introduction, an essential feature of the 
classical $H^1(\R^n)$ space is that every function in $H^1(\R^n)$ 
admits an atomic decomposition. Recall that an 
atom in $H^1(\R^n)$ is a measurable function $a\in L^2(\R^n)$, 
supported in a ball $B$, with zero integral and satisfying 
$\left\Vert a\right\Vert_2\leq \left\vert B\right\vert^{-1/2}$.  The Coifman-Latter theorem says that  an integrable 
function $f$ belongs to $H^1(\R^n)$ if and only if it can be written 
as
\[
f=\sum_{k\geq 1} \lambda_ka_k
\]
where $\sum_k\left\vert \lambda_k\right\vert<+\infty$ and the $a_k$'s 
are atoms. Moreover, $\left\Vert f\right\Vert_{H^1(\R^n)}$ is 
comparable with the infimum of 
$\sum\left\vert \lambda_k\right\vert$ over all such decompositions.

In \cite{loumcintosh}, Lou and the second author establish an 
atomic decomposition for $H^1_d(\R^n,\Lambda^k)$ for all $1\leq k\leq 
n$. In this context, an atom is a 
form $a\in L^2(\R^n,\Lambda^k)$ such that there exists $b\in 
L^2(\R^n,\Lambda^{k-1})$ supported in a ball $B\subset \R^n$ with radius 
$r$, $a=db$ and 
$\left\Vert a\right\Vert_2 + r^{-1} \left\Vert b\right\Vert_2\leq 
\left\vert B\right\vert^{-1/2}$. Note that the cancellation condition 
(in the case of functions) is replaced by the fact 
that an atom is the image of some other form under $d$ (that is, $a$ is exact), which 
implies in particular that $da=0$ whenever $a$ in an atom. The proof 
relies on a classical 
result due to Necas (\cite{necas}, Lemma 7.1, Chapter 3) and on 
(\cite{sc}, Theorem 3.3.3, Chapter 3).

In the present section, we prove a ``molecular'' decomposition for 
$H^1(\Lambda T^{\ast}M)$ inspired by the Coifman-Weiss terminology (see the introduction).  In our context, we do not know how to get atoms with compact support. Roughly speaking, a ``molecule'' is a form $f$ 
in $L^2(\Lambda T^{\ast}M)$ which is the image under $D^N$ of some $g\in 
L^2(\Lambda T^{\ast}M)$, with $L^2$ decay for $f$ and $g$, and for some integer $N$ large enough. 

To be more precise, we adopt the following terminology. Fix $C>0$. If 
$B\subset M$ is a ball with radius $r$ and if $(\chi_k)_{k\geq 0}$ is 
a sequence of nonnegative $C^{\infty}$ functions on $M$ with bounded 
support, say that $(\chi_k)_{k\geq 0}$ is adapted to $B$ if $\chi_0$ 
is supported in $4B$, $\chi_k$ is supported in $2^{k+2}B\setminus 
2^{k-1}B$ for all $k\geq 1$,
\begin{equation} \label{cutoff}
\sum_{k\geq 0} \chi_k=1\mbox{ on }M\mbox{ and }\left\Vert \left\vert 
\nabla \chi_k\right\vert\right\Vert_{\infty} \leq \frac C{2^kr},
\end{equation}
where $C>0$ only depends on $M$. Note that, when $C>0$ is large 
enough, there exist sequences adapted to any fixed ball.

Let $N$  be a positive integer. If $a\in L^2(\Lambda T^{\ast}M)$, $a$ is called an $N$-molecule if and 
only if there exists a ball 
$B\subset M$ with radius $r$, $b\in L^2(\Lambda T^{\ast}M)$ such that 
$a=D^Nb$, and a sequence $(\chi_k)_{k\geq 0}$ adapted to $B$ such 
that, for all $k\geq 0$,
\begin{equation} \label{molecestim}
\left\Vert \chi_k a\right\Vert_{L^2(\Lambda T^{\ast}M)}\leq 
2^{-k}V^{-1/2}(2^kB)\mbox{ and }\left\Vert \chi_k 
b\right\Vert_{L^2(\Lambda T^{\ast}M)}\leq 2^{-k}r^{N}V^{-1/2}(2^kB).
\end{equation}
Note that the first set of estimates in (\ref{molecestim}) imply
\begin{equation} \label{molecL2norm}  
\left\Vert a\right\Vert_{L^2(\Lambda T^{\ast}M)}\leq 2V^{-1/2}(B) \mbox{ and }
\left\Vert b\right\Vert_{L^2(\Lambda T^{\ast}M)}\leq 2 r^N V^{-1/2}(B)\ . 
\end{equation}
Thus $a\in \mathcal R(D)\subset H^2(\Lambda T^{\ast}M)$.   Furthermore, any $N$-molecule $a$ belongs to $L^1(\Lambda T^{\ast}M)$ and one has
\begin{equation} \label{molecL1norm}
\left\Vert a\right\Vert_{L^1(\Lambda T^{\ast}M)}\leq 2C
\end{equation}
where $C$ is the constant in (\ref{D}).

\begin{defi} Say that a section $f$ belongs to $H^1_{mol,N}(\Lambda T^{\ast}M)$ if there exists a 
sequence $(\lambda_j)_{j\geq 1}\in l^1$ and a sequence of $N$-molecules
$(a_j)_{j\geq 1}$ such that
\[
f=\sum_{j\geq 1} \lambda_ja_j,
\]
and define $\left\Vert f\right\Vert_{H^1_{mol,N}(\Lambda T^{\ast}M)}$ as the infimum of 
$\sum\left\vert \lambda_j\right\vert$ over all such decompositions. 
\end{defi}
It is plain to see that $H^1_{mol,N}(\Lambda T^{\ast}M)$ is a Banach space. We 
prove in this section the following
\begin{theo} \label{decompo}
Assume  (\ref{D}) and let  $\D$ is given by (\ref{Dbis}). Then, for 
integers $N> \frac {\D}2+1$,  
$H^1_{mol,N}(\Lambda T^{\ast}M)=H^1(\Lambda T^{\ast}M)$. As a consequence, $H^1_{mol,N}(\Lambda T^{\ast}M)$ 
is 
independent of $N$ provided that $N> \frac {\D}2+1$.
\end{theo}
\begin{cor} \label{Lp}
\begin{itemize}
\item[$(a)$]
For $1\leq p\leq 2$, $H^p(\Lambda T^{\ast}M)\subset \overline{{\mathcal R}(D)\cap L^p(\Lambda T^{\ast}M)}^{L^p(\Lambda T^{\ast}M)}$.
\item[$(b)$]
For $2\leq p<+\infty$, $\overline{{\mathcal R}(D)\cap L^p(\Lambda T^{\ast}M)}^{L^p(\Lambda T^{\ast}M)} \subset H^p(\Lambda T^{\ast}M)$.
\end{itemize}
\end{cor}
\noindent{\bf Proof of Corollary \ref{Lp}: } For assertion $(a)$, the inclusion $H^{1}(\Lambda T^{\ast}M)\subset L^{1}(\Lambda T^{\ast}M)$ is an 
immediate consequence of Theorem \ref{decompo} and of 
(\ref{molecL1norm}). Since 
$H^{2}(\Lambda T^{\ast}M)\subset L^{2}(\Lambda T^{\ast}M)$, we obtain by 
interpolation (Theorem \ref{interpohardy}) that $H^p(\Lambda T^{\ast}M)\subset L^p(\Lambda T^{\ast}M)$. 
Therefore 
\begin{align*} H^p(\Lambda T^{\ast}M) &= \overline{\mathcal R(D)\cap H^p(\Lambda T^{\ast}M)}^{H^p(\Lambda T^{\ast}M)}\\ &\subset
\overline{\mathcal R(D)\cap L^p(\Lambda T^{\ast}M)}^{L^p(\Lambda T^{\ast}M)}\ .
\end{align*}
For assertion $(b)$, observe first that, for all $1\leq p^{\prime}\leq 2$, there exists $C>0$ such that, for all $G\in T^{p^{\prime},2}(\Lambda T^{\ast}M)\cap T^{2,2}(\Lambda T^{\ast}M)$,
\begin{equation} \label{Scontrol}
\left\Vert {\mathcal S}_{\psi}G\right\Vert_{L^{p^{\prime}}(\Lambda T^{\ast}M)}\leq C\left\Vert G\right\Vert_{T^{p^{\prime},2}(\Lambda T^{\ast}M)}
\end{equation}
where   $\psi\in \Psi_{\beta,2}(\Sigma^0_{\theta})$ with $\beta=\left[\frac{\D}2\right]+1$.
Indeed, ${\mathcal S}_{\psi}G\in H^{p^{\prime}}(\Lambda T^{\ast}M)$, and therefore belongs to $L^{p^{\prime}}(\Lambda T^{\ast}M)$ by assertion $(a)$.

Let $p\geq 2$ and $f\in {\mathcal R}(D)\cap L^p(\Lambda T^{\ast}M)$. For all $G\in T^{p^{\prime},2}(\Lambda T^{\ast}M) \cap T^{2,2}(\Lambda T^{\ast}M)$ where $1/p+1/p^{\prime}=1$, since ${\mathcal Q}_\psi^{\ast}={\mathcal S}_{\overline{\psi}}$ (see Section \ref{h2}), one obtains, using (\ref{Scontrol}),
\[
\begin{array}{lll}
\displaystyle \left\vert\iint \langle ({\mathcal Q}_{\psi}f)_t(x),G(x,t)\rangle dx\frac{dt}t\right\vert&= & \displaystyle \left\vert \int_M \langle f(x),{\mathcal S}_{\overline{\psi}}G(x)\rangle dx\right\vert\\
& \leq & \left\Vert f\right\Vert_{L^p(\Lambda T^{\ast}M)} \left\Vert {\mathcal S}_{\overline{\psi}}G\right\Vert_{L^{p^{\prime}}(\Lambda T^{\ast}M)}\\
& \leq &  C\left\Vert f\right\Vert_{L^p(\Lambda T^{\ast}M)} \left\Vert G\right\Vert_{T^{p^{\prime},2}(\Lambda T^{\ast}M)},
\end{array}
\]
which shows that $\left\Vert {\mathcal Q}_{\psi}f\right\Vert_{T^{p,2}(\Lambda T^{\ast}M)}\leq C\left\Vert f\right\Vert_{L^p(\Lambda T^{\ast}M)}$ (remember that, by Proposition \ref{density}, $T^{p^{\prime},2}(\Lambda T^{\ast}M)\cap T^{2,2}(\Lambda T^{\ast}M)$ is dense in $T^{p^{\prime},2}(\Lambda T^{\ast}M)$), therefore $\left\Vert f\right\Vert_{H^p(\Lambda T^{\ast}M)}\leq C\left\Vert f\right\Vert_{L^p(\Lambda T^{\ast}M)}$. Next, if $f\in \overline{{\mathcal R}(D)\cap L^p(\Lambda T^{\ast}M)}^{L^p(\Lambda T^{\ast}M)}$, there exists a sequence $(f_j)_{j\geq 1}\in {\mathcal R}(D)\cap L^p(\Lambda T^{\ast}M)$ which converges to $f$ in the $L^p(\Lambda T^{\ast}M)$ norm, therefore in the $H^p(\Lambda T^{\ast}M)$ norm, which shows that $f\in H^p(\Lambda T^{\ast}M)$. \hfill\fin

\begin{rem} Note that the inclusion 
$H^{1}(\Lambda T^{\ast}M)\subset 
L^{1}(\Lambda T^{\ast}M)$ did not seem to be an immediate consequence of the 
definition 
of $H^{1}(\Lambda T^{\ast}M)$.
\end{rem}

\begin{rem}
What assertion $(b)$ in Corollary \ref{Lp} tells us is that, for all $\psi\in \Psi_{\beta,2}(\Sigma^{\theta}_0)$, for all $2\leq p< +\infty$, there exists $C>0$ such that, for all $f\in {\mathcal R}(D)\cap L^p(\Lambda T^{\ast}M)$, ${\mathcal Q}_{\psi}f\in T^{p,2}(\Lambda T^{\ast}M)$ and
\begin{equation} \label{Qcontrol}
\left\Vert {\mathcal Q}_{\psi}f\right\Vert_{T^{p,2}(\Lambda T^{\ast}M)}\leq C\left\Vert f\right\Vert_{L^p(\Lambda T^{\ast}M)}.
\end{equation}
In the Euclidean case and with the Laplacian on functions, this inequality is nothing but the well-known $L^p$-boundedness of the so-called Lusin area integral (for $p\geq 2$, it follows directly from the $L^p$-boundedness of the vertical quadratic $g$ function and the $L^p$ boundedness of the Hardy-Littlewood maximal function, see for instance \cite{stein70}, p.~91). In the context of spaces of homogeneous type, the $L^p$ boundedness (for all $1<p<+\infty$) of the area integral associated to an operator $L$ was proved in \cite{unpublished} under the following assumptions: $L$ is the generator of a holomorphic semigroup acting on $L^2$, the kernel of which satisfies Gaussian upper bounds, and $L$ has a bounded holomorphic calculus on $L^2$. Note that, in the framework of the present paper, we do not require any Gaussian upper estimate for the heat kernel of the Hodge Laplacian to obtain (\ref{Qcontrol}) for $p\geq 2$.
\end{rem}
As a consequence of Corollary \ref{Lp} and Theorem \ref{rieszbound}, 
we 
obtain the last part of Corollary \ref{riesz}: 
\begin{cor} \label{H1L1}
Assume  (\ref{D}). Then 
$D\Delta^{-1/2}$ is $H^{1}(\Lambda T^{\ast}M)-L^{1}(\Lambda T^{\ast}M)$ bounded.
\end{cor}
In Section \ref{appl} below, this theorem will be compared 
with previously known results for the Riesz transform on manifolds.
 
The proof of Theorem \ref{decompo} will be divided into two 
subsections, each corresponding to one inclusion.

\subsection{$H^1(\Lambda T^{\ast}M)\subset H^1_{mol,N}(\Lambda T^{\ast}M)$ for all $N\geq1$. } \label{firstinclusion}
Fix $N\geq 1$ and set $\beta=\left[\frac{\D}2\right]+1$ as usual. Choose $M\geq
\max\{\beta,N\}$ and define $\psi(z)=z^M(1+iz)^{-M-2}\in \Psi_{M,2}(\Sigma^0_+)$ and let $\phi(z)=z^{M-N}(1+iz)^{-M-2}$ so that $\psi(z)=z^N\phi(z)$. It is enough to prove that $E^1_{D,\psi}(\Lambda T^{\ast}M)\subset H^1_{mol,N}(\Lambda T^{\ast}M)$, 
which means that, if 
\[
f=\int_0^{+\infty}\psi_t(D) F_{t}\frac 
{dt}t
\]
with $F\in T^{1,2}(\Lambda T^{\ast}M)\cap T^{2,2}(\Lambda T^{\ast}M)$, then $f\in 
H^1_{mol,N}(\Lambda T^{\ast}M)$. According to the atomic decomposition for tent 
spaces (Theorem \ref{atomic}), one may assume that $F=A$ is a 
$T^{1,2}(\Lambda T^{\ast}M)$-atom, supported in $T(B)$ where $B$ is a ball in $M$ 
with radius $r$. Let
\[
g=\int_0^{+\infty} t^{N}\phi_t(D)A_{t}\frac {dt}t
\]
so that $f=D^Ng$, and let $(\chi_k)_{k\geq 0}$ be a sequence 
adapted to $B$. We claim that, up to a multiplicative constant, $f$ 
is a $N$-molecule, which gives the desired conclusion. To show this, we 
just have to establish the following $L^2$ estimates for $f$ and $g$: 
\begin{lem} \label{fgL2}
There exists $C>0$ only depending on $M$ such that, for all $k\geq 0$,
\[
\left\Vert \chi_k f\right\Vert_{L^2(\Lambda T^{\ast}M)}\leq 
C2^{-k}V^{-1/2}(2^kB),\ \left\Vert \chi_k 
g\right\Vert_{L^2(\Lambda T^{\ast}M)}\leq Cr^{N}2^{-k}V^{-1/2}(2^kB).
\]
\end{lem}
\noindent{\bf Proof: }We first deal with the estimates for $f$. First, since $D$ is self-adjoint, one has
\[
\left\Vert f\right\Vert_{L^2(\Lambda T^{\ast}M)}^2\leq C\int_M\int_0^{+\infty} 
\left\vert A(x,t)\right\vert^2 dx\frac{dt}t\leq CV(B)^{-1}.
\]
This shows that $\left\Vert \chi_0 f\right\Vert_{L^2(\Lambda T^{\ast}M)}\leq 
\left\Vert f\right\Vert_2\leq CV(B)^{-1/2}$.

Fix now $k\geq 1$ and $\displaystyle m\geq \frac {\D}2+1$. Lemma 
\ref{offdiag1} and the fact that $A$ is supported in $T(B)$ yield
\[
\begin{array}{lll}
\displaystyle \left\Vert \chi_kf\right\Vert_{L^2(\Lambda T^{\ast}M)} & \leq & \displaystyle \int_0^r \left\Vert 
\chi_{2^{k+1}B\setminus 2^{k-1}B}\psi_t(D)
A_{t}\right\Vert_{L^2(\Lambda T^{\ast}M)}\frac{dt}t\\
& \leq & \displaystyle C\int_0^r \left(\frac 
t{2^kr}\right)^m\left\Vert A_{t}\right\Vert_{L^2(\Lambda T^{\ast}M)} \frac{dt}t\\
& \leq & \displaystyle C\left(\int_0^r \left(\frac 
t{2^kr}\right)^{2m}\frac{dt}t\right)^{1/2} \left(\int_0^r \left\Vert 
A_{t}\right\Vert_{L^2(\Lambda T^{\ast}M)}^2
\frac{dt}t\right)^{1/2}\\
& \leq & \displaystyle C(2^kr)^{-m}r^m V(B)^{-1/2}\\
& \leq & \displaystyle C2^{-k\left(m-\frac {\D}2\right)} 
V^{-1/2}(2^kB).
\end{array}
\]
We now turn to the estimates on $g$. First,
\[
\begin{array}{lll}
\displaystyle \left\Vert g\right\Vert_{L^2(\Lambda T^{\ast}M)} & \leq & \displaystyle \int_0^r t^{N} \left\Vert A_{t}\right\Vert_{L^2(\Lambda T^{\ast}M)} 
\frac{dt}t\\
& \leq & \displaystyle \left(\int_0^r  t^{2N}\frac{dt}t\right)^{1/2} 
\left(\int_0^r \left\Vert A_{t}\right\Vert_{L^2(\Lambda T^{\ast}M)}^2 
\frac{dt}t\right)^{1/2}\\
& \leq  &\displaystyle C r^{N} V(B)^{-1/2},
\end{array}
\]
which shows that $\left\Vert \chi_0 g\right\Vert_{L^2(\Lambda T^{\ast}M)}\leq 
\left\Vert g\right\Vert_2\leq C r^{N} V(B)^{-1/2}$. 

Fix now $k\geq 1$ and $\displaystyle m\geq \frac {\D}2+1$. Lemma 
\ref{offdiag1} yields
\[
\begin{array}{lll}
\displaystyle \left\Vert \chi_k g\right\Vert_{L^2(\Lambda T^{\ast}M)} & \leq & \displaystyle \int_0^r  t^{N} \left\Vert 
\chi_{2^{k+1}B\setminus 
2^{k-1}B} \phi_t(D)A_{t}\right\Vert_{L^2(\Lambda T^{\ast}M)} \frac{dt}t\\
& \leq & \displaystyle C\int_0^r  t^{N}\left(\frac t{2^kr}\right)^m 
\left\Vert A_{t}\right\Vert_{L^2(\Lambda T^{\ast}M)} \frac{dt}t\\
& \leq & \displaystyle C(2^kr)^{-m} \left(\int_0^r 
 t^{2N+2m} \frac{dt}t\right)^{1/2} \left(\int_0^r \left\Vert 
A_{t}\right\Vert_{L^2(\Lambda T^{\ast}M)}^2\frac{dt}t\right)^{1/2}\\
& \leq & \displaystyle C2^{-km}  r^N V^{-1/2}(B)\\
& \leq & \displaystyle C r^{N} V^{-1/2}(2^kB)2^{-k\left(m-\frac 
{\D}2\right)}.
\end{array}
\]
This completes the proof of Lemma \ref{fgL2}, and provides the 
desired inclusion. \hfill\fin

\subsection{$H^1_{mol,N}(\Lambda T^{\ast}M)\subset H^1(\Lambda T^{\ast}M)$ for all $N>\frac\kappa2+1$.} \label{secondinclusion}
For the converse inclusion, it is enough to prove that there exists 
$C>0$ such that, for every $N$- molecule $f$ in $H^1_{mol,N}(\Lambda T^{\ast}M)$, $f\in 
H^1(\Lambda T^{\ast}M)$ with $\left\Vert 
f\right\Vert_{H^1(\Lambda T^{\ast}M)}\leq C$.

Let $f$ be such a $N$-molecule. Since $f\in {\mathcal R}(D)$, 
according to Lemma \ref{equivalent}, it suffices to show that, if 
 $F(x,t)=({\mathcal Q}_{\psi}f)_t(x)$ with $\psi(z)=z(1+iz)^{-\gamma -2}\in \Psi_{1,\beta+1}(\Sigma^0_\theta)$ for some $\gamma>N-1$, then
\begin{equation} \label{Ftent}
\left\Vert F\right\Vert_{T^{1,2}(\Lambda T^{\ast}M)}\leq C,
\end{equation}
There exists a ball $B$, a section $g\in L^{2}(\Lambda T^{\ast}M)$ and a sequence 
$(\chi_k)_{k\geq 0}$ adapted to $B$ such that $f=D^Ng$ and 
(\ref{molecestim}) holds. Define
\[
\eta_0=\chi_{2B\times \left(0,2r\right)}
\]
and, for all $k\geq 1$,
\[
\eta_k=\chi_{(2^{k+1}B\setminus 2^kB)\times \left(0,r\right)},\ 
\eta^{\prime}_k=\chi_{(2^{k+1}B\setminus 2^kB)\times 
\left(r,2^{k+1}r\right)},\ \eta^{\prime\prime}_k=\chi_{2^kB\times 
\left(2^kr,2^{k+1}r\right)},
\]
where these functions $\chi_S$ are the (un-smoothed) characteristic functions
of $S\subset M\times(0,\infty)$. Write
\[
F=\eta_0 F+\sum\limits_{k\geq 1} \eta_k F+\sum\limits_{k\geq 1} 
\eta^{\prime}_k F + \sum\limits_{k\geq 1} \eta^{\prime\prime}_k F.
\]
The estimate (\ref{Ftent}) will be an immediate consequence of the 
following
\begin{lem} \label{Ftentestimates}
\begin{itemize}
\item[$(a)$]    
For each $k\geq 0$, $\left\Vert 
\eta_{k}F\right\Vert_{T^{1,2}(\Lambda T^{\ast}M)}\leq 
C2^{-k}$.
\item[$(b)$]
For each $k\geq 1$, $\left\Vert 
\eta^{\prime}_{k}F\right\Vert_{T^{1,2}(\Lambda T^{\ast}M)}\leq 
C2^{-k}$.
\item[$(c)$]
For each $k\geq 1$, $\left\Vert 
\eta^{\prime\prime}_{k}F\right\Vert_{T^{1,2}(\Lambda T^{\ast}M)}\leq 
C2^{-k}$.
\end{itemize}
\end{lem}
\noindent{\bf Proof: }\par
{\bf Assertion $(a)$: } Since $\eta_{k}F$ is supported in 
the box ${\mathcal B}(2^{k+1}B)$ (see Remark \ref{box}), we just have 
to prove that its $T^{2,2}(\Lambda T^{\ast}M)$ norm is 
controlled by $C2^{-k}V^{-1/2}(2^kB)$ (recall that the $T^{2,2}(\Lambda T^{\ast}M)$ 
norm is equivalent to the norm in ${\mathcal H}$, see Section 
\ref{mainestim}), which will prove that $\frac 
1C 2^k 
\eta_{k}F$ is an atom in $T^{1,2}(\Lambda T^{\ast}M)$. First, by the spectral theorem, one has
\[
\begin{array}{lll}
\displaystyle \left\Vert \eta_{0}F\right\Vert_{T^{2,2}(\Lambda T^{\ast}M)}^2 & \leq & \displaystyle 
\left\Vert F\right\Vert_{T^{2,2}(\Lambda T^{\ast}M)}^2\\
& \leq & \displaystyle C \int_0^{+\infty} \left\Vert 
 \psi_t(D) f\right\Vert_{L^2(\Lambda T^{\ast}M)}^2\frac{dt}t\\
& \leq & \displaystyle C\left\Vert f\right\Vert_{L^2(\Lambda T^{\ast}M)}^2\\
& \leq & \displaystyle CV(B)^{-1}.
\end{array}
\]
Fix now $k\geq 1$. One has
\[
\begin{array}{lll}
\displaystyle \left\Vert \eta_k F\right\Vert_{T^{2,2}(\Lambda T^{\ast}M)} & \leq  & \displaystyle \sum_{l\geq 0} 
\left\Vert 
\chi_{(2^{k+1}B\setminus 
2^kB)\times 
(0,r)}  \psi_t(D) (\chi_lf)\right\Vert_{T^{2,2}(\Lambda T^{\ast}M)}\\
& := & \displaystyle \sum_{l\geq 0} I_{l}.
\end{array}
\] 
Assume that $0\leq l\leq k-2$. Then, using (\ref{Dbis}), Lemma \ref{offdiag1}, (\ref{molecestim}) for $f$ and 
the fact that $\rho\left(\mbox{supp }\chi_{l}f,2^{k+1}B\setminus 
2^kB\right)\geq c(2^k-2^l)r$ and choosing $m\geq \frac {\D}2+1$,
\[
\begin{array}{lll}
I_{l}^2 & =& \displaystyle \int_0^r \left(\int_{2^{k+1}B\setminus 
2^kB} 
\left\vert 
 \psi_t(D) (\chi_lf)(x)\right\vert^2dx\right)\frac{dt}t\\
& \leq & \displaystyle C\int_0^r \left(\frac 
t{(2^k-2^l)r}\right)^{2m}\left\Vert \chi_l 
f\right\Vert_{L^2(\Lambda T^{\ast}M)}^2\frac{dt}t\\
& \leq  & \displaystyle C2^{k({\D}-2m)}2^{-l({\D}+2)}V^{-1}(2^kB).
\end{array}
\]
It follows that
\[
\sum\limits_{l=0}^{k-2} \left\Vert 
\chi_{(2^{k+1}B\setminus 
2^kB)\times (0,r)} \psi_t(D) (\chi_lf)\right\Vert_{T^{2,2}(\Lambda T^{\ast}M)}\leq C2^{-k}V^{-1/2}(2^kB).
\]
Assume now that $k-2\leq l\leq k+2$. Then, by the spectral theorem,
\[
\begin{array}{lll}
I_{l}^2 & \leq & \displaystyle \int_0^{+\infty} \left\Vert 
 \psi_t(D) (\chi_lf)\right\Vert_{L^2(\Lambda T^{\ast}M)}^2\frac{dt}t\\
& \leq & \displaystyle C\left\Vert \chi_l f\right\Vert_{L^2(\Lambda T^{\ast}M)}^2\\
& \leq & \displaystyle C2^{-2k}V^{-1}(2^kB).
\end{array}
\]
Assume finally that $l\geq k+3$. Then, using Lemma \ref{offdiag1} and 
$\rho\left(\mbox{supp }\chi_{l}f,2^{k+1}B\setminus 
2^kB\right)\geq c2^lr$, and choosing $m\geq \frac {\D}2+1$,
\[
\begin{array}{lll}
I_{l}^2 & \leq & \displaystyle \int_0^r \left(\frac 
t{2^lr}\right)^{2m}\left\Vert \chi_lf\right\Vert_{L^2(\Lambda T^{\ast}M)}^2 
\frac{dt}t\\
& \leq & \displaystyle C2^{-2lm}2^{-2l}V^{-1}(2^kB).
\end{array}
\]
As a consequence,
\[
\sum\limits_{l=k+3}^{\infty} \left\Vert 
\chi_{(2^{k+1}B\setminus 
2^kB)\times (0,r)} \psi_t(D) (\chi_lf)\right\Vert_{T^{2,2}(\Lambda T^{\ast}M)}\leq C2^{-k}V^{-1/2}(2^kB).
\]
This ends the proof of assertion $(a)$ in Lemma \ref{Ftentestimates}.

{\bf Assertion $(b)$: }Similarly, we now estimate
\[
\begin{array}{lll}
\displaystyle \left\Vert \eta^{\prime}_kF\right\Vert_{T^{2,2}(\Lambda T^{\ast}M)}  & \leq  & \displaystyle \sum\limits_{l\geq 0} \left\Vert 
\chi_{(2^{k+1}B\setminus 2^kB)\times 
(r,2^{k+1}r)} \psi_t(D)D^N
(\chi_l g)\right\Vert_{T^{2,2}(\Lambda T^{\ast}M)} \\
& := & \displaystyle \sum\limits_{l\geq 0} J_{l}.
\end{array}
\] 
Define now $\widetilde{\psi}(z)=z^N\psi(z)\in \Psi_{N+1,\gamma +1-N}(\Sigma^0_\theta)$, so that
\[
J_l^2= \int_{r}^{2^{k+1}r} \left\Vert 
\chi_{2^{k+1}B\setminus 
2^kB}  \widetilde{\psi}_t(D) (\chi_{l}g)\right\Vert_{L^{2}(\Lambda T^{\ast}M)}^{2}
\frac{dt}{t^{ 2N +1}}.
\]
Assume first that $0\leq l\leq k-2$. Then, (\ref{molecestim}) applied to $g$, Lemma \ref{offdiag1} and 
the support conditions on $\chi_l$ yield
\[
\begin{array}{lll}
\displaystyle 
J_{l}^{2} 
& \leq & \displaystyle C\int_{r}^{2^{k+1}r} \left(\frac 
t{(2^k-2^l)r}\right)^{2m} \left\Vert 
\chi_{l}g\right\Vert_{L^2(\Lambda T^{\ast}M)}^{2} 
\frac{dt}{t^{ 2N +1}}\\
& \leq & \displaystyle C2^{-l(\D+2)}2^{-k(2m-\D)}V^{-1}(2^kB),
\end{array}
\]
if $m$ is chosen so that $\D+2\leq 2m< 2N$, which is possible 
since $N> \frac {\D}2+1$. This yields
\[
\sum\limits_{l=0}^{k-2} \left\Vert 
\chi_{(2^{k+1}B\setminus 
2^kB)\times (r,2^{k+1}r)} \psi_t(D)D^{N} (\chi_{l}g)\right\Vert_{L^{2}(\Lambda T^{\ast}M)} \leq
C2^{-k}V^{-1/2}(2^kB).
\]
Assume now that $k-1\leq l\leq k+1$. Then one has
\[
\begin{array}{lll}
\displaystyle J_{l}^{2}  
& \leq & \displaystyle C2^{-2k}V^{-1}(2^kB).
\end{array}
\]
Assume finally that $l\geq k+2$. Then, using Lemma \ref{offdiag1} and 
the support conditions again,
\[
\begin{array}{lll}
\displaystyle 
J_{l}^{2} 
& \leq & \displaystyle \int_{r}^{2^{k+1}r} \left(\frac 
t{2^lr}\right)^{2m} \left\Vert \chi_{l}g\right\Vert_{L^2(\Lambda T^{\ast}M)}^{2} 
\frac{dt}{ t^{2N+1}} \\
& \leq & \displaystyle C2^{-2(m+1)l}V^{-1}(2^kB),
\end{array}
\]
provided that  $m<N$. Thus, 
\[
\sum\limits_{l=k+2}^{+\infty} \left\Vert 
\chi_{(2^{k+1}B\setminus 
2^kB)\times (r,2^{k+1}r)} \psi_t(D)D^N(\chi_{l}g)\right\Vert_{T^{2,2}(\Lambda T^{\ast}M)} \leq C2^{-k}V^{-1/2}(2^kB)
\]
and assertion $(b)$ is proved.

{\bf Assertion $(c)$: }Finally, as in assertions $(a)$ and $(b)$, we have to estimate
\[
\begin{array}{lll}
\displaystyle \left\Vert 
\eta^{\prime\prime}_kF\right\Vert_{T^{2,2}(\Lambda T^{\ast}M)} & \leq & \displaystyle 
\sum\limits_{l\geq 0} \left\Vert 
\chi_{2^kB\times (2^kr,2^{k+1}r)} \psi_t(D)D^N
(\chi_l g)\right\Vert_{T^{2,2}(\Lambda T^{\ast}M)}\\
& := & \displaystyle \sum\limits_{l\geq 0} K_l.
\end{array}
\] 
Similarly, one has
\[
K_l^2=\int_{2^kr}^{2^{k+1}r} \left\Vert 
\chi_{2^kB} \widetilde{\psi}_t(D)
(\chi_l g)\right\Vert_{L^2(\Lambda T^{\ast}M)}^2\frac{dt}{t^{ 2N+1}}.
\]
Assume first that $0\leq l\leq k$. Then, 
one obtains
\[
\begin{array}{lll}
K_l^2 
& \leq & \displaystyle  C\int_{2^kr}^{2^{k+1}r} \frac{dt}{t^{ 2N+1}} 
\left\Vert \chi_lg\right\Vert_{L^2(\Lambda T^{\ast}M)}^2 \\
& \leq & \displaystyle C2^{-(\D+2)l} V^{-1}(2^kB)2^{ k(\D-2N)}.
\end{array}
\]
Since $2N\geq \D+1$, it follows that
\[
\sum\limits_{l=0}^k  \left\Vert 
\chi_{2^kB\times (2^kr,2^{k+1}r)} \psi_t(D)D^N
(\chi_l g)\right\Vert_{T^{2,2}(\Lambda T^{\ast}M)}\leq 
C2^{-k}V^{-1/2}(2^kB).
\]
Assume now that $l\geq k+1$. Then, using Lemma \ref{offdiag1} once 
more,
\[
\begin{array}{lll}
\displaystyle K_l^2 
& \leq & \displaystyle C\int_{2^kr}^{2^{k+1}r} \left(\frac 
t{2^lr}\right)^{2m} \frac{dt}{t^{ 2N+1}} \left\Vert 
\chi_lg\right\Vert_{L^2(\Lambda T^{\ast}M)}^2\\
& \leq & \displaystyle C2^{-(2m+2)l}2^{ k(2m-2N)}V^{-1}(2^kB),
\end{array}
\]
provided that  $m<N$. It follows that
\[
\sum\limits_{l\geq k+1} \left\Vert 
\chi_{2^kB\times (2^kr,2^{k+1}r)} \psi_t(D)D^N
(\chi_l g)\right\Vert_{T^{2,2}(\Lambda T^{\ast}M)} \leq 
C2^{-k}V^{-1/2}(2^kB).
\]
Assertion $(c)$ is therefore proved. \hfill\fin\par
\begin{rem} \label{deltamolec}
Let $N>\frac12(\frac\kappa2+1)$  be an integer.  Say that $a\in L^2(\Lambda T^{\ast}M)$ is an $N$-molecule for $\Delta$  if there exists $b\in L^2(\Lambda T^{\ast}M)$ such that $a=\Delta^Nb$, a ball $B\subset M$ with radius $r$ and a sequence $(\chi_k)_{k\geq 0}$ adapted to $B$ such that, for all $k\geq 0$,
\begin{equation} \label{molecestimdelta}
\left\Vert \chi_ka\right\Vert_2\leq 2^{-k}V^{-1/2}(2^kB)\mbox{ and }\left\Vert \chi_kb\right\Vert_2\leq 2^{-k}r^{2N}V^{-1/2}(2^kB).
\end{equation}
Moreover, define $H^1_{\Delta,mol,N}$ as the space of all sections $f$ such that $f=\sum_{j\geq 1} \lambda_ja_j$, where $\sum_j\left\vert \lambda_j\right\vert<+\infty$ and the $a_j$'s are $N$-molecules for $\Delta$,  and equip $H^1_{\Delta,mol,N}$ with the usual norm. Then, using Theorem \ref{decompo}, it is plain to see that, for all $N>\frac12(\frac\kappa2+1)$, $H^1_{\Delta,mol,N}=H^1(\Lambda T^{\ast}M)$.
\end{rem}
\par

Finally, we also have a decomposition into molecules for $H^{1}_{d}(\Lambda T^{\ast}M)$ 
and $H^{1}_{d^{\ast}}(\Lambda T^{\ast}M)$. Let $N> \frac{\D}2+1$. An $N$-molecule for $d$ is a 
section $a\in L^2(\Lambda T^{\ast}M)$, such that there exists a ball 
$B\subset M$ with radius $r$, $b\in L^2(\Lambda T^{\ast}M)$ with 
  $a=dD^{N-1}b$  and a sequence $(\chi_k)_{k\geq 0}$ adapted to $B$ such 
that, for all $k\geq 0$,
\begin{equation} \label{molecestimbis}
\left\Vert \chi_k a\right\Vert_{L^2(\Lambda T^{\ast}M)}\leq 
2^{-k}V^{-1/2}(2^kB)\mbox{ and }\left\Vert \chi_k 
b\right\Vert_{L^2(\Lambda T^{\ast}M)}\leq 2^{-k} r^{N} V^{-1/2}(2^kB).
\end{equation}
Then, $f\in H^{1}_{d}(\Lambda T^{\ast}M)$ if and only if $f=\sum_{j}\lambda_ja_{j}$ 
where the $a_{j}$'s are atoms in $H^{1}_{d}(\Lambda T^{\ast}M)$ and $\sum\left\vert 
\lambda_{j}\right\vert<+\infty$.  The proof is analogous and uses the characterization of $H^1_d(\Lambda T^{\ast}M)$ by means of ${\mathcal S}_{d,\psi}$ and ${\mathcal Q}_{d^{\ast},\psi}$ given in Section \ref{hodgedecompo}.   One obtains a similar decomposition for 
$H^{1}_{d^{\ast}}(\Lambda T^{\ast}M)$, defining an $N$-molecule for $d^{\ast}$ similarly to an $N$-molecule for $d$.
\begin{rem} It turns out that, under some Gaussian upper estimates for the heat kernel of the Hodge-de Rham Laplacian, we can take $N=1$ in Theorem \ref{decompo} and other similar results. We will come back to this in Section \ref{molecgauss}.
\end{rem}

\SE{The maximal characterization} \label{max}
In this section, we provide a characterization of $H^1(\Lambda T^{\ast}M)$ in 
terms of maximal functions. Recall that, for classical Hardy spaces 
of functions in the Euclidean case, such maximal functions are 
defined, for instance, in the following way: if $\int_{\R^n} \left\vert f(y)\right\vert (1+\left\vert y\right\vert^2)^{-(n+1)/2}<+\infty$  and $x\in \R^n$, define
\[
f^{\ast}(x)=\sup_{\left\vert y-x\right\vert<t} \left\vert 
e^{-t\sqrt{\Delta}}f(y)\right\vert.
\]
Then, a possible characterization of $H^{1}(\R^n)$ is the following 
one: $f\in H^1(\R^n)$ if and only if $f^{\ast}\in L^1(\R^n)$. \par
In the present context, such a definition has to be adapted. Let us 
explain the main lines before coming to the details. First, the lack 
of pointwise estimates forces us to replace the value at $(y,t)$ by 
an $L^{2}$ average on a ball centered at $(y,t)$. Secondly, the 
Poisson semigroup (on forms) $e^{-t\sqrt\Delta}$ only satisfies 
$OD_{t}(1)$ estimates in general, which is unsufficient to carry out 
the argument in \cite{feffermanstein} or its adaptation in \cite{ar}. 
Hence, we abandon in the maximal function the Poisson semigroup in favor of 
the heat semigroup. Thirdly, the good-$\lambda$ argument of 
\cite{feffermanstein} or \cite{ar} with the heat semigroup produces 
then uncontrolled error terms involving the time derivatives due to 
the parabolic nature of the equation associated with. The trick is to 
modify the maximal function to incorporate the errors in the very 
definition of the maximal function (see the function 
$\tilde{f}^{\ast}_{\alpha,c}$ below) so that they are under control in 
the argument.

In the sequel, if $x\in M$ and $0<r<t$, $B((x,t),r)=B(x,r)\times 
\left(t-r,t+r\right)$. For all $x\in M$ and all $\alpha>0$, recall 
that
\[
\Gamma_{\alpha}(x)=\left\{(y,t)\in M\times (0,+\infty);\ y\in 
B(x,\alpha t)\right\}.
\]
Let $0<\alpha$. Fix $c>0$ such that, for all $x\in M$, whenever 
$(y,t)\in \Gamma_{\alpha}(x)$, $B((y,t),ct)\subset 
\Gamma_{2\alpha}(x)$. Elementary geometry shows that $c\le 
\alpha/(1+2\alpha)$ works. For $f\in L^2(\Lambda T^{\ast}M)$ and all $x\in M$, define
\[
f^{\ast}_{\alpha,c}(x)=\sup_{(y,t)\in \Gamma_{\alpha}(x)} \left( 
\frac 1{tV(y,t)} \iint_{B((y,t),ct)} \left\vert 
e^{-s^2\Delta}f(z)\right\vert^2 dzds\right)^{1/2}
\]
Define  $H^1_{max}(\Lambda T^{\ast}M)$  as the completion of 
$\{ f\in {\mathcal R}(D);  \left\Vert 
f^{\ast}_{\alpha,c}\right\Vert_{L^1(M)}<\infty\}$ for that norm and 
set 
\[
\left\Vert f\right\Vert_{H^1_{max}(\Lambda T^{\ast}M)}=\left\Vert 
f^{\ast}_{\alpha,c}\right\Vert_{L^1(M)}.
\]
The norm depends a priori on $\alpha,c$. However, the doubling 
condition 
(\ref{D}) allows us to compare them. For fixed $\alpha$,  the 
pointwise bound $f^{\ast}_{\alpha,c} \le 
C(1+ c/c')^{\D/2} f^{\ast}_{2\alpha,c'}$  holds if $0<c\le 
\alpha/(1+2\alpha)$ and $0<c^{\prime}\leq 2\alpha/(1+4\alpha)$. Next, if $0<\alpha\le \beta$ and $c \le  
\alpha/(1+2\alpha)$  then $f^{\ast}_{\alpha,c} \le 
f^{\ast}_{\beta,c}$ while if $0<\beta<\alpha$ and 
$c\le \beta/(1+2\beta)$ (hence $c\le \alpha/(1+2\alpha)$), we have 
$\|f^{\ast}_{\alpha,c}\|_{1}\le C 
(\beta/\alpha)^{\D}\|f^{\ast}_{\beta,c}\|_{1}$ by a variant of the 
Fefferman-Stein argument in \cite{feffermanstein} which is skipped. 
Hence, 
the space $H^1_{max}(\Lambda T^{\ast}M)$ is independent from the choice of 
$\alpha,c$. 
Notice also that, in the definition of $f^{\ast}_{\alpha,c}$, because 
of the doubling property again,
replacing $V(y,t)$ by $V(y,ct)$ yields an equivalent norm.

 The following characterization holds as part of Theorem \ref{equality}
\begin{theo} \label{maxchara}
Assume (\ref{D}). Then 
$H^1(\Lambda T^{\ast}M)=H^1_{max}(\Lambda T^{\ast}M).$\end{theo}

\begin{rem}
The average in $s$ in the definition of $f^{\ast}_{\alpha,c}$ is 
useful 
only in the proof of the inclusion $H^1_{max}(\Lambda T^{\ast}M)\subset H^1(\Lambda T^{\ast}M)$. Equivalent norms occur without the average is $s$ in the definition. 
\end{rem}

As a consequence of this result and Corollary \ref{H1L1}, we have
\begin{cor} \label{H1L1bis}
Assume  (\ref{D}). Then 
$D\Delta^{-1/2}$ is $H^{1}_{max}(\Lambda T^{\ast}M)-L^{1}(\Lambda T^{\ast}M)$ bounded.
\end{cor}

For the proof of the theorem, we introduce an auxiliary space  $\widetilde 
H^1_{max}(\Lambda T^{\ast}M)$ and show the following chain of inclusions: 
$H^1(\Lambda T^{\ast}M)\subset H^1_{\max}(\Lambda T^{\ast}M)= \widetilde H^1_{max}(\Lambda T^{\ast}M)\subset 
H^1(\Lambda T^{\ast}M)$.  This space is built as $ H^1_{max}(\Lambda T^{\ast}M)$ with $ 
f^{\ast}_{\alpha,c}$ changed to 
\begin{equation}
\label{ftildestar}
\tilde f^{\ast}_{\alpha,c}(x)=\sup_{(y,t)\in \Gamma_{\alpha}(x)} 
\left( \frac 1{tV(y,t)} \iint_{B((y,t),ct)} \left\vert 
e^{-s^2\Delta}f(z)\right\vert^2 + \left\vert s \frac 
\partial{\partial s}e^{-s^2\Delta}f(z)\right\vert^2 dzds\right)^{1/2}.
\end{equation}

\subsection{Proof of $H^1(\Lambda T^{\ast}M)\subset H^1_{max}(\Lambda T^{\ast}M)$}
Fix $\alpha=1/2$ and $0<c\le 1/4$ and set $f^\ast=f^{\ast}_{1/2,c}$. 
In 
view of Theorem \ref{decompo}, it is enough to show that any $N$-molecule for $\Delta$ (for suitable $N$) in 
$H^1(\Lambda T^{\ast}M)$ belongs to $H^1_{max}(\Lambda T^{\ast}M)$.
We denote by ${\mathcal M}$ the usual Hardy-Littlewood maximal 
function:
\[
{\mathcal M}f(x)=\sup_{B\ni x} \frac 1{V(B)}\int_B \left\vert 
f(y)\right\vert dy,
\]
where the supremum is taken over all the balls $B\subset M$ 
containing $x$. 
Here is our first technical lemma:
\begin{lem} \label{techmax}
Assume that $(T_t)_{t>0}$ satisfies $OD_t(N)$ estimates with 
$N>{\D}/2$. Then, there exists $C>0$ such that, for all $f\in 
L^2_{loc}(\Lambda T^{\ast}M)$, all $x\in M$ and all $(y,t)\in \Gamma_{\alpha}(x)$,
\[
\frac 1{tV(y,ct)}\iint_{B((y,t),ct)} \left\vert 
T_sf(z)\right\vert^2dzds \leq C{\mathcal M}(\left\vert 
f\right\vert^2)(x).
\]
\end{lem}
\noindent{\bf Proof: }Decompose $f=\sum\limits_k f_k$, with 
$f_0=\chi_{B(y,2ct)}f$ and $f_k=\chi_{B(y,2^{k+1}ct)\setminus 
B(y,2^kct)}f$ for all $k\geq 1$ (where $\chi_A$ stands for the 
characteristic function of $A$). For $k=0$, the $L^2$-boundedness of 
$T_s$ and the fact that $s\sim t$ and $V(x, (2c+1)t) \sim V(y,ct)$ 
yield
\[
\displaystyle \frac 1{tV(y,ct)}\iint_{B((y,t),ct)} \left\vert 
T_sf_0(z)\right\vert^2dzds  \leq  \displaystyle 
C{\mathcal M}(\left\vert f\right\vert^2)(x).
\]
For $k\geq 1$, the estimate $OD_t(N)$ and the fact that $s\sim t$  
and 
${V(x,(2^{k+1}c+1)t)} \le C2^{k\D} {V(y,ct)}$ give us
\[
\begin{array}{lll}
\displaystyle \frac 1{tV(y,ct)}\iint_{B((y,t),ct)} \left\vert 
T_sf_k(z)\right\vert^2dzds  &\leq & \displaystyle \frac C{V(y,ct)} 
\left(\frac 1{2^k}\right)^{2N} \int \left\vert f_k(z)\right\vert^2 
dz\\
& \leq & \displaystyle \frac {C2^{k\D}}{2^{2kN}}{\mathcal 
M}(\left\vert f\right\vert^2)(x).
\end{array}
\]
Since $2N>\D$, one can sum up these estimates by the Minkowski 
inequality. \hfill\fin

We prove now that $H^1(\Lambda T^{\ast}M)\subset H^1_{max}(\Lambda T^{\ast}M)$. Let $a=\Delta^{N_0}b \in {\mathcal R}(D)$   
be a $N_0$-molecule for $\Delta$ in $H^1(\Lambda T^{\ast}M)$, for some $N_0\geq \frac {\D}2+1$, $r>0$ 
the radius of the ball $B$ associated with $a$, $(\chi_j)_{j\geq 0}$ 
a sequence adapted to $B$  such that (\ref{molecestimdelta}) holds (see Remark \ref{deltamolec}) . For 
each $j\geq 0$, set $a_j=\chi_j a$.
First, the Kolmogorov inequality (\cite{meyer}, p. 250), Lemma \ref{techmax}, the maximal 
theorem and  the doubling property (\ref{D}) show that
\[
\int_{2B} a^{\ast}(x)dx\leq CV(2B)^{1/2}\left\Vert {\mathcal 
M}\left(\vert a\vert^2\right)\right\Vert_{1,\infty}^{1/2}\leq 
CV(B)^{1/2} \left\Vert a\right\Vert_2\leq C.
\]
We next show that, for some $\delta>0$ only depending on doubling 
constants and all $k\geq 1$,
\begin{equation} \label{astar}
\int_{2^{k+1}B\setminus 2^kB} a^{\ast}(x)dx\leq C2^{-k\delta}.
\end{equation}
Fix $k\geq 1$ and write $a^{\ast}\leq 
a^{\ast}_{low}+a^{\ast}_{medium}+a^{\ast}_{high}$, where 
$a^{\ast}_{low}$ (resp. $a^{\ast}_{medium}$,$a^{\ast}_{high}$) 
correspond to the constraint $t<r$ (resp. $r\leq t<2^{k-1}r$,\ $t\geq 
2^{k-1}r$) in the supremum defining $a^{\ast}$.
We first deal with $a^{\ast}_{low}$. According to the definition of 
$a_{j}$ we have $a^{\ast}_{low}\leq \sum_{j\geq 0} a^{\ast}_{j,low}$. 
Fix 
$(y,t)\in \Gamma_{\alpha}(x)$, $t<r$, $x\in 2^{k+1}B\setminus 2^kB$. 
If 
$j\leq k-2$, so that $\rho(B(y,ct),\mbox{supp }a_j)\sim 2^kr$, then, 
using 
the fact that $s\sim t$ and off-diagonal estimates for $e^{-s^2\Delta}$ (Lemma \ref{gaffneyheat})  and arguing as in Lemma \ref{techmax} (using 
the 
fact that $t<r$), one obtains, provided that $2N>{\D}$,
\[
\begin{array}{lll}
\displaystyle \frac 1{t}\frac 1{V(y,ct)} \iint_{B(y,t),ct)} 
\left\vert 
e^{-s^2\Delta}a_j(z)\right\vert^2dzds & \leq & \displaystyle 
C\left(\frac 
t{2^kr}\right)^{2N-{\D}} {\mathcal M}\big(\left\vert 
a_j\right\vert^2\big)(x)\\
& \leq & \displaystyle \frac C{2^{k(2N-{\D})}} {\mathcal 
M}\big(\left\vert a_j\right\vert^2\big)(x).
\end{array}
\]
If $j\geq k+2$, so that $\rho(B(y,ct),\mbox{supp }a_j)\sim 2^jr$, 
then, one has similarly
\[
\frac 1{t}\frac 1{V(y,ct)} \iint_{B(y,t),ct)} \left\vert 
e^{-s^2\Delta}a_j(z)\right\vert^2dzds \leq \frac C{2^{j(2N-{\D})}} 
{\mathcal M}\big(\left\vert a_j\right\vert^2\big)(x).
\]
Setting
\[
c_{j,k}=\left\{
\begin{array}{ll}
2^{-k(N-{\D}/2)} &\mbox{if }j\leq k-2,\\
1 & \mbox{if }k-1\leq j\leq k+1,\\
2^{-j(N-{\D}/2)}& \mbox{if }j\geq k+2,
\end{array}
\right.
\]
one therefore has, using the Kolmogorov inequality and 
(\ref{molecestim}) again,
\[
\begin{array}{lll}
\displaystyle \int_{2^{k+1}B\setminus 2^kB} a^{\ast}_{low}& \leq & \displaystyle \sum_{j\geq 0} \int_{2^{k+1}B\setminus 2^kB} 
a^{\ast}_{j,low}\\
& \leq & \displaystyle C\sum_{j\geq 0} c_{j,k} 
V^{1/2}(2^{k+1}B)2^{-j}V^{-1/2}(2^{j}B)\\
& \leq & \displaystyle C\sum_{j\geq 0} c_{j,k} 
\sup\left(1,2^{(k-j){\D}/2}\right)2^{-j}\\
& \leq & \displaystyle C2^{-k\delta}
\end{array}
\]
if $N>{\D}+\delta$.

To estimate $a^{\ast}_{medium}$ on $2^{k+1}B\setminus 2^kB$, write 
$a=\Delta^{N_0} b$ and $b=\sum_{j\geq 0} b_j$ where $b_j=\chi_j b$. 
Let 
$(y,t)\in \Gamma_{\alpha}(x)$, $x\in 2^{k+1}B\setminus 2^kB$ and 
$r\leq 
t<2^{k-1}r$. Then, $\rho(B(y,ct),\mbox{supp }b_j)\sim 2^kr$ if $j\leq 
k-2$ and $\rho(B(y,ct),\mbox{supp }b_j)\sim 2^jr$ if $j\geq k+2$. 
Hence, 
arguing as before, using off-diagonal estimates for $(s^2\Delta)^{N} 
e^{-s^2\Delta}$ (see Lemma \ref{gaffneyheat} again) 
and $s\sim t$, one has, if $2N>{\D}$,
\[
\frac 1{t}\frac 1{V(y,ct)} \iint_{B((y,t),ct)} \left\vert 
\Delta^{N_0}e^{-s^2\Delta}b_j(z)\right\vert^2dzds\leq \left\{
\begin{array}{ll}
\displaystyle \frac 1{t^{4N_0}}\left(\frac t{2^kr}\right)^{2N-{\D}} 
{\mathcal M}\big(\left\vert b_j\right\vert^2\big)(x) &\mbox{if }j\leq 
k-2,\\
\displaystyle \frac 1{t^{4N_0}}{\mathcal M}\big(\left\vert 
b_j\right\vert^2\big)(x)&\mbox{if }|j-k|\leq1,\\
\displaystyle \frac 1{t^{4N_0}}\left(\frac 
t{2^jr}\right)^{2N-{\D}}{\mathcal M}\big(\left\vert 
b_j\right\vert^2\big)(x) &\mbox{if }j\geq k+2,
\end{array}
\right.
\]
and, if we choose $N$ such that $N< \frac {\D}2+2N_0$, one obtains
\[
\frac 1{t}\frac 1{V(y,ct)}\iint_{B((y,t),ct)} \left\vert 
\Delta^{N_0}e^{-s^2\Delta}b_j(z)\right\vert^2dzds \leq \frac 
1{r^{4N_0}}c_{j,k}^2{\mathcal M}\big(\left\vert 
b_j\right\vert^2\big)(x)
\]
where $c_{j,k}$ was defined above. Thus, provided that $N>{\D}$,
\[
\int_{2^{k+1}B\setminus 2^kB} a^{\ast}_{medium} \leq C2^{-k\delta}
\]
for some $\delta>0$. Note that this choice of $N$ is possible since 
$N_{0}\geq \frac {\D}2+1$.

It remains to look at $a^{\ast}_{high}$. Let $x\in 2^{k+1}B\setminus 
2^kB$, $(y,t)\in \Gamma_{\alpha}(x)$ and $t>2^{k-1}r$. Using 
$a=\Delta^{N_0} b$ and the $L^2$-boundedness of 
$(s^2\Delta)^{N_0}e^{-s^2\Delta}$, one obtains
\[
\begin{array}{lll}
\displaystyle \frac 1{t}\frac 1{V(y,ct)}\iint_{B((y,t),ct)}\left\vert 
\Delta^{N_0}e^{-s^2\Delta}b(z)\right\vert^2dzds & \leq & \displaystyle \frac C{t^{4N_0}}\int \left\vert b\right\vert^2 \frac 
1{V(y,ct)}\\
& \leq &\displaystyle C2^{-k(4N_0-\D)}\frac 1{V^2(2^{k+1}B)}.
\end{array}
\]
It follows that
\[
\int_{2^{k+1}B\setminus 2^kB} a^{\ast}_{high}\leq C2^{-k(2N_0-\D/2)}.
\]
Finally, (\ref{astar}) is proved since $N_0\geq \D/2+1$. Thus, 
$\int_M 
a^{\ast}\leq C$, which ends the proof of the inclusion 
$H^1(\Lambda T^{\ast}M)\subset 
H^1_{max}(\Lambda T^{\ast}M)$. \hfill\fin
\subsection{$H^1_{max}(\Lambda T^{\ast}M)=  \widetilde H^1_{max}(\Lambda T^{\ast}M)$}
As $f^{\ast}_{\alpha,c} \le \tilde f^{\ast}_{\alpha,c}$, it follows 
that  
$H^1_{max}(\Lambda T^{\ast}M)\supset  \widetilde H^1_{max}(\Lambda T^{\ast}M)$. We turn to the 
opposite 
inclusion. For that argument we fix $\alpha=1$ and $c\le 1/12$ and 
write $f^\ast$, $\tilde f^\ast$. The term in $\tilde f^\ast$ coming 
from $e^{-s^2\Delta}f(z)$ is immediately controlled by $f^\ast$. 
Next, fix $x\in M$ and $(y,t) \in \Gamma_{1}(x)$. Let $(z,s) \in 
B((y,t),ct)$. Write
$$
s\partial_{s}e^{-s^2\Delta}f(z) = -2s^2\Delta e^{-s^2\Delta}f(z)= 
(-2s^2\Delta e^{-(\sqrt 3 s/2)^2\Delta})e^{-(s/2)^2\Delta}f(z)
$$
and observe that $-2s^2\Delta e^{-(\sqrt 3 s/2)^2\Delta}$ satisfies 
$OD_{s}(N)$ for any $N$  (see Lemma \ref{gaffneyheat}). 
Since $s\sim t$,  
$$
\left(\frac 1 {tV(y,t)}\iint_{B((y,t), ct)} \left\vert 
s\partial_{s}e^{-s^2\Delta}f(z)\right\vert^2 dzds\right)^{1/2}
$$
is controlled by 
\[
\begin{array}{ll}
      & \displaystyle  C\left(\frac 1 
{tV(y,t)}\int_{t-ct}^{t+ct}\int_{B(y,2ct)} \left\vert 
e^{-(s/2)^2\Delta}f(z)\right\vert^2 dzds\right)^{1/2}     \\
      &\displaystyle  \quad  +C \sum_{k=1}^\infty  2^{-kN} 
\left(\frac 1 
{tV(y,t)}\int_{t-ct}^{t+ct}\int_{B(y,2^{k+1}ct)\setminus B(y,2^kct)} 
\left\vert e^{-(s/2)^2\Delta}f(z)\right\vert^2 dzds\right)^{1/2}. 
\end{array}
\]
 Change $s/2$ to $s$. This first term is controlled by $f^\ast_{2, 
4c}(x)$.  For the $k$th term in the series, one covers 
$B(y,2^{k+1}ct)$ by balls $B(y_{j}, 2ct)$ 
 with $y_{j}\in B(y,2^{k+1}ct)$ and  $B(y_{j},ct)$ pairwise disjoint. 
By doubling, the balls $B(y_{j}, 2ct)$ have bounded overlap. Observe 
also that each point 
 $(y_{j}, t/2)$ belongs to $\Gamma_{2+2^{k+2}c}(x)$. Hence, the $k$th 
term is bounded by
\begin{align*}
\label{}
   C   2^{-kN} \left(\sum_{j} \frac{V(y_{j}, 2t)}{V(y,t)} 
\right)^{1/2} 
   f^\ast_{2+2^{k+2}c, 4c}(x) & \le C \left(\frac {V(y, 
   2^{k+1}ct+2t)}{V(y,t)}\right)^{1/2} f^\ast_{2+2^{k+2}c, 4c}(x)
   \\
   & \le C 2^{k\D/2} f^\ast_{2+2^{k+2}c, 4c}(x).
\end{align*}
 Hence, using the comparisons between the $\|f^\ast_{\alpha,c}\|_{1}$ 
norms, we obtain
  \begin{align*}
\label{}
   \| \tilde f^\ast\|_{1}&\le \|f^\ast\|_{1}+  C\|f^\ast_{2, 
4c}\|_{1}+ 
   C\sum_{k=1}^\infty 2^{-k(N-\D/2)}\|f^\ast_{2+2^{k+2}c, 4c}\|_{1}.
   \\
    &      \le \|f^\ast\|_{1}(1+ C+C\sum_{k=1}^\infty 
    2^{-k(N-\D/2)}2^{k\D})
\end{align*} 
so that if $N>3\D/2$, we obtain $ \| \tilde f^\ast\|_{1}\le C 
\|f^\ast\|_{1}$. \hfill\fin
 
\subsection{$\widetilde H^1_{max}(\Lambda T^{\ast}M)\subset H^1(\Lambda T^{\ast}M)$}
We fix some notation. If $f\in L^2(\Lambda T^{\ast}M)$,  the section 
$u(y,t)=e^{-t^2\Delta}f(y)$  for all $y\in M$ and all $t>0$
 satisfies the equation
\begin{equation} \label{harmonic}
D^2u= \Delta u=-\frac{1}{2t} \frac{\partial u}{\partial t}.
\end{equation}
For all $\alpha>0$, $0\leq \varepsilon<R\leq +\infty$,  $x\in M$ 
set the truncated cone
\[
\Gamma_{\alpha}^{\varepsilon,R}(x)=\left\{(y,t)\in M\times 
\left(\varepsilon,R\right);\ y\in B(x,\alpha t)\right\} = 
\left\{(y,t)\in \Gamma_{\alpha}(x);\varepsilon<t<R\right\}.
\]
and for all $f\in L^2(\Lambda T^{\ast}M)$,
\[
S_{\alpha}^{\varepsilon,R}f(x)=\left(\iint_{\Gamma_{\alpha}^{\varepsilon,R}(x)} 
\frac{\left\vert 
t{D}u(y,t)\right\vert^2}{V(y,t)}dy\frac{dt}t\right)^{1/2}
\]
where $\left\vert {D}u(y,t)\right\vert^2=\langle 
Du(y,t),Du(y,t)\rangle$ (remember that $\langle.,.\rangle$ denotes the 
inner product in $\Lambda T^{\ast}M$). Note that, contrary to the definition of $\S F$ given in Section \ref{tent}, we divide here by $V(y,t)$ instead of $V(x,t)$. By (\ref{D}), this amounts to the same, since, when $(y,t)\in \Gamma_{\alpha}(x)$, $d(y,x)\leq \alpha t$, and it turns out that $V(y,t)$ is more handy here. For our 
purpose, it is enough to show that, for all $f\in  L^2(\Lambda T^{\ast}M)\cap 
\widetilde H^{1}_{max}(\Lambda T^{\ast}M)$ and all 
$0<\varepsilon<R$,
\begin{equation} \label{Lusin}
\left\Vert S_{\alpha}^{\varepsilon,R}f\right\Vert_1\leq C\left\Vert 
f\right\Vert_{\widetilde H^1_{max}(\Lambda T^{\ast}M)}.
\end{equation}
Indeed,  if furthermore  $f \in {\mathcal R}(D)\cap \widetilde H^1_{max}(\Lambda T^{\ast}M)$, letting 
$\varepsilon$ go to $0$ and $R$ to $+\infty$, this means that 
$(y,t)\mapsto tDe^{-t^2\Delta}f(y)\in T^{1,2}(\Lambda T^{\ast}M)$ and, since
\[
f=a\int_0^{+\infty} 
(tD)^{2N_1+1}(I+t^{2}D^2)^{-\alpha_1}tDe^{-t^2D^{2}}f\frac{dt}t
\]
 for suitable integers $N_1,\alpha_1$ and constant $a$ (we use 
 $\Delta=D^2$), this yields $f\in H^1(\Lambda T^{\ast}M)$ by definition of the 
Hardy space. 
    
 The proof of (\ref{Lusin}) is inspired by the one of Proposition 7 
in \cite{ar} where the Poisson semigroup is changed to the heat 
semigroup. We first need the following inequality (see Lemma 8 in 
\cite{ar}):
\begin{lem} \label{L1}
There exists $C>0$ such that, for all $f\in L^{2}(\Lambda T^{\ast}M)$, all $0<\varepsilon<R<+\infty$ 
and all $x\in M$,
\[
S^{\varepsilon,R}_{1/20}(x)\leq 
C\left(1+\ln(R/\varepsilon)\right)\widetilde{f}^{\ast}_{1,c}(x).
\]
\end{lem}

\noindent{\bf Proof: }Fix any $0<c<1/3$. Let $x\in M$. One 
can cover $\Gamma^{\varepsilon,R}_{1/20}(x)$ by balls in $M\times (0,+\infty)$ in 
the following way: for each 
$l\in \Z$, let $(B(x_{j,l},\frac c2\tau^l))_{j\in \Z}$ be a covering of $M$ 
by 
balls (where $\tau=\frac{1+c/2}{1-c/2}$) so that the balls 
$B\left(x_{j,l},\frac c4\tau^l\right)$ are pairwise disjoint. Let 
$K_{j,l}=B(x_{j,l},\frac c2 \tau^l)\times \left[\tau^l-\frac c2 
\tau^l,\tau^l+\frac c2\tau^l\right]$ and 
$\widetilde{K_{j,l}}=B(x_{j,l},c\tau^l)\times 
\left[\tau^l-c\tau^l,\tau^l+c\tau^l\right]=B((x_{j,l},\tau^l),c\tau^l)$.  
Since, for $(y,t)\in K_{j,l}$, if $K_{j,l}\cap 
\Gamma^{\varepsilon,R}_{1/20}(x)\neq \emptyset$, one has $t\sim\tau^l$ 
and $V(y,t)\sim V(x_{j,l},\tau^l)$, we obtain
\[
\begin{array}{lll}
\displaystyle S^{\varepsilon,R}_{1/20}f(x)^{2} &\leq & \displaystyle \sum_{l,j;\ K_{j,l}\cap 
\Gamma^{\varepsilon,R}_{1/20}(x)\neq \emptyset} \iint_{K_{j,l}} 
\frac{\left\vert tDu(y,t)\right\vert^{2}}{V(y,t)}dy\frac{dt}t\\
& \leq & \displaystyle C\sum_{l,j;\ K_{j,l}\cap 
\Gamma^{\varepsilon,R}_{1/20}(x)\neq\emptyset} 
\frac{\tau^l}{V(x_{j,l},\tau^l)} \iint_{K_{j,l}} \left\vert 
Du(y,t)\right\vert^{2}dydt.
\end{array}
\]
At this stage, we need the following parabolic Caccioppoli inequality 
recalling that $u(y,t)=e^{-t^2\Delta}f(y)$: for some constant $C>0$ 
only depending on $M$, but not on $j$, $l$, 
\begin{equation} \label{caccio}
\iint_{K_{j,l}} \left\vert {D}u(y,t)\right\vert^2dydt\leq 
C\tau^{-2l}\iint_{\widetilde K_{j,l}} \left\vert u(y,t)\right\vert^2dydt 
+ C\tau^{-2l}\iint_{\widetilde K_{j,l}} \left\vert  t \partial 
_{t}u(y,t)\right\vert^2dydt.
\end{equation}
The proof of this inequality is classical and will therefore be 
skipped (see for instance \cite{asterisque}, Chapter 1). Since, by the 
choice of $c$, $x_{j,l}\in 
\Gamma^{\varepsilon,R}_{1}(x)$ whenever $K_{j,l}\cap 
\Gamma^{\varepsilon,R}_{1/20}(x)\neq \emptyset$, it follows 
that
\[
\begin{array}{lll}
\displaystyle  S^{\varepsilon,R}_{1/20}f(x)^{2} &\leq &  \displaystyle 
C\sum_{l,j;\ K_{j,l}\cap 
\Gamma^{\varepsilon,R}_{1/20}(x)\neq\emptyset}  
\frac{1}{\tau^lV(x_{j,l},\tau^l)} \iint_{\widetilde K_{j,l}}  \left(\left\vert 
u(y,t)\right\vert^2 + \left\vert  t \partial 
_{t}u(y,t)\right\vert^2\right) dydt\\
& \leq & \displaystyle C \sharp\left\{(l,j);\ K_{j,l}\cap 
\Gamma^{\varepsilon,R}_{1/20}(x)\neq\emptyset\right\}\tilde 
f^{\ast}_{1,c}(x)^{2}.
\end{array}
\]
For fixed $l\in \Z$, the bounded overlap property of the balls 
$B(x_{j,l},c\tau^l)$ implies that the number of $j$'s such that $K_{j,l}\cap 
\Gamma^{\varepsilon,R}_{1/20}(x)\neq\emptyset$ is uniformly bounded 
with respect to $l$. Now, if $l\in \Z$ is such that $K_{j,l}\cap 
\Gamma^{\varepsilon,R}_{1/20}(x)\neq\emptyset$, then one has 
$(1+c)^{-1}\varepsilon<\tau^l\leq R(1-c)^{-1}$, which yields the desired 
conclusion. \hfill\fin

\medskip

\noindent{\bf Proof of (\ref{Lusin}): }it suffices to 
establish the following ``good $\lambda$'' inequality:
\begin{lem} \label{goodlambda}
There exists $C>0$ such that, for all $0\leq \varepsilon<R\leq 
+\infty$, 
all $f\in L^2(\Lambda T^{\ast}M)$, all $0<\gamma<1$ and all $\lambda>0$,
\begin{equation} \label{goodineq}
\mu\left(\left\{x\in M;\ S_{1/20}^{\varepsilon,R}f(x)>2\lambda,\ 
\tilde f^{\ast}(x)\leq \gamma\lambda\right\}\right) \leq 
C\gamma^2\mu\left(\left\{x\in M;\ 
S_{1/2}^{\varepsilon,R}f(x)>\lambda\right\}\right)
\end{equation}
where $\tilde f^*=\tilde f^{\ast}_{1,c}$ is as in (\ref{ftildestar})  
with $0<c<1/3$  to be chosen in the proof. 
\end{lem}
Indeed, assume that Lemma \ref{goodlambda} is proved. Then, if $f\in 
L^{2}(\Lambda T^{\ast}M)\cap\tilde{H}^{1}_{max}(\Lambda T^{\ast}M)$, integrating (\ref{goodineq}) with 
respect to $\lambda$ and using Remark
\ref{aperture} yield
\[
\left\Vert S_{1/20}^{\varepsilon,R}f\right\Vert_{1}\leq 
C\gamma^{-1}\left\Vert \tilde f^{\ast}\right\Vert_{1}+ C\gamma^{2} 
\left\Vert S_{1/2}^{\varepsilon,R}f\right\Vert_{1}\leq 
C\gamma^{-1}\left\Vert \tilde f^{\ast}\right\Vert_{1}+ C^{\prime}\gamma^{2} 
\left\Vert S_{1/20}^{\varepsilon,R}f\right\Vert_{1}.
\]
Since $f\in \widetilde{H}^{1}_{max}(\Lambda T^{\ast}M)$, Lemma 
\ref{L1} ensures that $\left\Vert 
S_{1/20}^{\varepsilon,R}f\right\Vert_{1}<+\infty$, and (\ref{Lusin}) 
follows at once if $\gamma$ is chosen small enough.\par

\bigskip

\noindent{\bf Proof of Lemma \ref{goodlambda}: }Assume first that $M$ 
is unbounded, which, by (\ref{D}), implies that $\mu(M)=+\infty$ (see 
\cite{martellthesis}). Let $O=\left\{x\in M; 
S_{1/2}^{\varepsilon,R}f(x)>\lambda\right\}$. Observe that, since 
$f\in L^{2}(\Lambda T^{\ast}M)$, $S^{\varepsilon,R}_{1/2}f\in L^{2}(\Lambda T^{\ast}M)$, whence 
$\mu(O)<+\infty$, therefore $O\neq M$. Morever, $O$ is open since the 
map $x\mapsto S_{1/2}^{\varepsilon,R}f(x)$ is continuous. Let $(B_k)_{k\geq 1}$ be a Whitney decomposition of $O$, so that 
$2B_{k}\subset O$ and $4B_k\cap (M\setminus O)\neq \emptyset$\footnote{To be correct 
$4B_{k}$ should be $c_{1}B_{k}$ where $c_{1}$ depends on the doubling 
property. To avoid too many constants, we set $c_{1}=4$ to fix 
ideas.} for all $k\geq 1$. For all $k\geq 1$, define
\[
E_k=\left\{x\in B_k;\ S_{1/20}^{\varepsilon,R}f(x)>2\lambda,\ \tilde 
f^{\ast}(x)\leq \gamma\lambda\right\}.
\]
Because of the bounded overlap property of the $B_k$'s, and since 
$\left\{S^{\varepsilon,R}_{1/20}f>2\lambda\right\}\subset 
\left\{S^{\varepsilon,R}_{1/2}f>\lambda\right\}$, it is enough to 
prove that
\begin{equation} \label{enoughgoodlambda}
\mu(E_k)\leq C\gamma^2\mu(B_k).
\end{equation}
Observe first that, if $\varepsilon\geq 20 r(B_k)$ (where $r(B_k)$ is 
the radius of $B_k$), then $E_k=\emptyset$. Indeed, there exists 
$x_k\in 4B_k$ such that $S^{\varepsilon,R}_{1/2}f(x_k)\leq \lambda$. 
Let now $x\in B_k$ and $(y,t)\in \Gamma_{1/20}^{\varepsilon,R}(x)$. 
Then,
\[
d(x_k,y)\leq d(x_k,x)+d(x,y)\leq 5r(B_k)+\frac t{20}\leq \frac 
{5\varepsilon}{20}+\frac t{20}\leq \frac{6t}{20}<\frac t2,
\]
so that $(y,t)\in \Gamma_{1/2}^{\varepsilon,R}(x_k)$. As a 
consequence, $S^{\varepsilon,R}_{1/20}f(x)\leq 
S^{\varepsilon,R}_{1/2}f(x_k)\leq \lambda$.
We may therefore assume that $\varepsilon<20r(B_k)$. Since one has 
$S^{20r(B_k),R}_{1/20}f(x)\leq\lambda$ by similar arguments, we 
deduce that
\[
E_k\subset\widetilde{E_k}=\left\{x\in B_k\cap F;\ 
S_{1/20}^{\varepsilon,20r(B_k)}f(x)>\lambda\right\}
\]
where
\[
F=\left\{x\in M\, ;\, \tilde f^{\ast}(x)\leq \gamma\lambda\right\}
\]
(note that $F$ is closed). By Tchebycheff inequality,
\[
\begin{array}{lll}
\displaystyle \mu\left(\widetilde{E_k}\right) & \leq & \displaystyle 
\frac 1{\lambda^2} \int_{B_k\cap F} \left\vert 
S^{\varepsilon,20r(B_k)}_{1/20}f(x)\right\vert^2dx\\
& = & \displaystyle \frac 1{\lambda^2}\int_{x\in B_k\cap F} 
\left(\iint_{\varepsilon<t<20r(B_k),\ y\in B(x,t/20)} 
\frac{\left\vert 
t{D}u(y,t)\right\vert^2}{V(y,t)}dy\frac{dt}t\right)dx\\
& \leq &\displaystyle \frac C{\lambda^2} 
\iint_{\Omega^{\varepsilon}_k} t^{2}\left\vert 
{D}u(y,t)\right\vert^2dy \frac{dt}t,
\end{array}
\]
where $\Omega^{\varepsilon}_k$ is the region in $M\times 
\left(0,+\infty\right)$ defined by the following conditions:
\[
\varepsilon<t<20r(B_k),\ \psi(y)<t/20
\]
with
\[
\psi(y)=\rho(y,B_k\cap F).
\]
Note that $\Omega^{\varepsilon}_k\subset \widetilde{\Omega_k}$, where 
the region $\widetilde{\Omega_k}$ is defined by $0<\psi(y)<t$. By 
definition of $F$, one has 
\begin{equation} \label{omegaktilde}
\frac 1{tV(y,t)}\iint_{B((y,t),ct)}\left\vert u(s,z)\right\vert^2+  
\left\vert s {\partial_{s} u(z,s)}\right\vert^2 dzds\leq 
\gamma^2\lambda^2
\end{equation}
for all $(y,t)\in \widetilde{\Omega_k}$.
To avoid the use of surface measure on 
$\partial\widetilde{\Omega_k}$, let us introduce
\[
\zeta(y,t)=\eta^2\left(\frac{\psi(y)}t\right)\chi_1^2\left(\frac 
t{\varepsilon}\right)\chi_2^2\left(\frac t{20r(B_k)}\right),
\]
where $\eta,\chi_1$ and $\chi_2$ are nonnegative $C^{\infty}$ 
functions on $\R$, $\eta$ is supported in $\left[0,\frac 
1{10}\right]$ and is equal to $1$ on $\left[0,\frac 1{20}\right]$,  
$\chi_1$ is supported in $\left[\frac 9{10},+\infty\right[$ and is 
equal to $1$ on $\left[1,+\infty\right[$, and $\chi_2$ is supported 
in $\left[0,\frac{11}{10}\right]$ and is equal to $1$ on 
$\left[0,1\right]$. One therefore has
\[
\mu\left(\widetilde{E_k}\right)\leq \frac 1{\lambda^2} \iint 
\zeta(y,t)\left\vert {D}u(y,t)\right\vert^2 dydt:=\frac 1{\lambda^2}I.
\]
The integral is over $M\times (0,\infty)$. An integration by parts in space 
using $D^*=D$   yields 
\[
\begin{array}{lll}
I & = & \Re \displaystyle \iint \langle {D} 
\left(\zeta(y,t){D}u(y,t)\right),u(y,t)\rangle  t dydt\\
& = & \Re \displaystyle \iint \langle {D}\zeta(y,t) Du(y,t),  
u(y,t)\rangle tdydt\\
& + & \Re \displaystyle \iint \langle \zeta(y,t) {D}^2u(y,t),
u(y,t)\rangle tdydt.
\end{array}
\]
In the first integral in the right hand side $D\zeta Du$ is the 
Clifford product (its exact expression is not relevant as we merely 
use $|D\zeta Du| \le |D\zeta||Du|$).
In the second, we use (\ref{harmonic}) and since $\zeta$ is 
real-valued, $\Re \langle\zeta D^2u,u\rangle=  -\frac 1 {4t} 
\zeta \partial_{t}\vert u \vert ^2$. Then integration by parts in $t$ 
gives us
\[
\begin{array}{lll}
I & = & \Re  \displaystyle \iint t \langle {D}\zeta(y,t){D}u(y,t), 
u(y,t)\rangle dydt\\
& + & \frac 1 4  \displaystyle \iint \partial _{t 
}\zeta(y,t)\left\vert u(y,t)\right\vert^2dydt\\
& := & I_1+I_2.
\end{array}
\]
\
\noindent{\bf Estimates and support considerations:}
Observe that, because of the support conditions on $\eta,\chi_1$ and 
$\chi_2$ and  since $\psi$ is a Lipschitz function with Lipschitz 
constant 1,
\begin{equation} \label{gradient}
\left\vert {D}\zeta(y,t)\right\vert + \left\vert 
\partial_{t}\zeta(y,t)\right\vert \leq \frac Ct
\end{equation}
independently of $\varepsilon$ and $k$. 
Now let us look more closely at the supports of $D\zeta$ and 
$\partial_{t}\zeta$. Examination shows that they are both supported
in the region  $G_k^{\varepsilon}$ of $M\times (0,\infty)$ defined by 
$0<\psi(y)<\frac {t}{10}$, $\frac {9\varepsilon}{10}<t<22r(B_k)$ and 
\[
\frac t{20}<\psi(y)<\frac t{10}\ \mbox{ or }\ \frac 
{9\varepsilon}{10}<t<\varepsilon \ \mbox{ or }\ 20r(B_k)<t<22r(B_k).
\]
Observe that $G_k^{\varepsilon}\subset \widetilde{\Omega_k}$. 
Consider again the balls $K_{j,l}$ introduced above. 
It is possible to choose $c$ small enough such that,   for all $k, j, 
l, \varepsilon$,  if
$K_{j,l}\cap G_k^{\varepsilon}\neq\emptyset$, then 
$\widetilde{K_{j,l}}\subset \widetilde{\Omega_k}$. Thus, 
(\ref{omegaktilde}) yields
$$
\iint_{\widetilde{K_{j,l}}}\left\vert u(z,s)\right\vert^2+  
\left\vert s  {\partial_{s} u(z,s)}\right\vert^2 dzds\leq 
\gamma^2\lambda^2 \tau^l V(x_{j,l},\tau^l) \leq C\gamma^2\lambda^2 \tau^l 
V(x_{j,l},c\tau^l)
$$
where we used the doubling property in the last inequality.

\medskip

\noindent {\bf Estimate of $I_2$:}  Using the considerations above 
and $t\sim \tau^l$ on $K_{j,l}$
\[
\begin{array}{lll}
\left\vert I_2\right\vert  & \leq & \displaystyle 
C\iint_{G_k^{\varepsilon}} \left\vert u(y,t)\right\vert^2 
dy\frac{dt}t\\
& = & \displaystyle C\sum_{l,j;\ K_{j,l}\cap 
G_k^{\varepsilon}\neq\emptyset} \iint_{K_{j,l}} \left\vert 
u(y,t)\right\vert^2dy\frac{dt}t\\
& \leq & \displaystyle C \sum_{l,j;\ K_{j,l}\cap 
G_k^{\varepsilon}\neq\emptyset} \tau^{-l} \iint_{K_{j,l}} \left\vert 
u(y,t)\right\vert^2dydt\\
& \leq & \displaystyle C\gamma^2\lambda^2 \sum_{l,j;\ K_{j,l}\cap 
G_k^{\varepsilon}\neq\emptyset} V(x_{j,l},c\tau^l)\\
& \leq & \displaystyle C\gamma^2\lambda^2 
\iint_{\widetilde{G_k^{\varepsilon}}} dy\frac{dt}t,
\end{array}
\]
where $\widetilde{G_k^{\varepsilon}}$ is the region defined by 
$0<\psi(y)<5t,\ \frac{9\varepsilon}{100}<t<100r(B_k)$ and
\[
\frac t{40}<\psi(y)<5t\ \mbox{ or }\ \frac 
{9\varepsilon}{20}<t<10\varepsilon \ \mbox{ or } \ 2r(B_k)<t<100r(B_k).
\]
The last inequality is due to the bounded overlap property of the 
balls $B(x_{j,l},c\tau^l)$ for each $l\in \Z$ and $t\sim \tau^l$ on each of 
them. Thus,
\[
\vert I_2\vert\leq C\gamma^2\lambda^2 \mu(H_k^{\varepsilon}),
\]
where $H_k^{\varepsilon}=\left\{y\in M;\ \exists{t>0 },\ (y,t)\in 
\widetilde{G_k^{\varepsilon}}\right\}$. It remains to observe that 
$H_k^{\varepsilon}\subset 221B_k$. Indeed, if $y\in H_k^{\varepsilon}$ 
and $t>0$ is such that $(y,t)\in \widetilde{G_k^{\varepsilon}}$, one 
has $\psi(y)<5t$, so that there exists $z\in B_k\cap F$ such that 
$\rho(y,z)<5t<220r(B_k)$. Thus, $y\in 221B_k$. Using the doubling 
property, we have therefore obtained
\[
\vert I_2\vert \leq C\gamma^2\lambda^2 V(B_k).\ \hfill\fin
\]

\medskip

{\bf Estimate of $I_1$: } Using the same notation, one has
\[
\begin{array}{lll}
\left\vert I_1\right\vert & \leq & \displaystyle 
C\iint_{G_k^{\varepsilon}} \left\vert Du(y,t)\right\vert 
\left\vert u(y,t)\right\vert dydt\\
& \leq & \displaystyle C\sum_{l,j;\ K_{j,l}\cap G_k^{\varepsilon}\neq 
\emptyset} \iint_{K_{j,l}} \left\vert Du(y,t)\right\vert 
\left\vert u(y,t)\right\vert dydt\\
& \leq & \displaystyle C\sum_{l,j;\ K_{j,l}\cap G_k^{\varepsilon}\neq 
\emptyset} \left(\iint_{K_{j,l}} \left\vert 
Du(y,t)\right\vert^2dydt\right)^{1/2} 
\left(\iint_{K_{j,l}} \left\vert u(y,t)\right\vert^2dydt\right)^{1/2}.
\end{array}
\]
The Caccioppoli inequality (\ref{caccio}) yields
\[
\begin{array}{lll}
\left\vert I_1\right\vert & \leq & \displaystyle  C\sum_{l,j;\ 
K_{j,l}\cap G_k^{\varepsilon}\neq \emptyset} \tau^{-l} 
\iint_{\widetilde{K_{j,l}}} \left\vert u(y,t)\right\vert^2 dydt + 
\tau^{-l} \iint_{\widetilde{K_{j,l}}} \left\vert 
t\partial_{t}u(y,t)\right\vert^2 dydt\\
& \leq & \displaystyle C\gamma^2\lambda^2 \sum_{l,j;\ K_{j,l}\cap 
G_k^{\varepsilon}\neq \emptyset} V(x_{j,l},c\tau^l),
\end{array}
\]
and the same computations as before yield
\[
\left\vert I_1\right\vert \leq C\gamma^2\lambda^2 V(B_k).\ \hfill\fin
\]
Finally, (\ref{enoughgoodlambda}) holds and Lemma \ref{goodlambda} is 
proved when $M$ is unbounded.\par

\medskip

When $M$ is bounded, call $\delta$ the diameter of $M$. We claim that there exists a constant $C>0$ such that, for all $R\geq 20\delta$ and all $x\in M$,
\begin{equation} \label{claimconst}
S_{1/20}^{20\delta,R}f(x)\leq C\tilde f^{\ast}(x).
\end{equation}
Assume that (\ref{claimconst}) is proved. It is enough to prove Lemma \ref{goodlambda} for $\gamma$ small, say $\gamma\leq 1/C$, where $C$ is the constant in (\ref{claimconst}). In this case, (\ref{claimconst}) ensures that, 
if $S_{1/20}^{20\delta,R}f(x)>\lambda$, then $\tilde f^{\ast}(x)>\gamma\lambda$, so that it remains to establish that
\begin{equation} \label{goodineqbounded}
\mu\left(\left\{x\in M;\ S_{1/20}^{\varepsilon,20\delta}f(x)>\lambda,\ 
\tilde f^{\ast}(x)\leq \gamma\lambda\right\}\right) \leq 
C\gamma^2\mu\left(\left\{x\in M;\ 
S_{1/2}^{\varepsilon,R}f(x)>\lambda\right\}\right).
\end{equation}
If $O=\left\{x\in M;\ 
S_{1/2}^{\varepsilon,R}f(x)>\lambda\right\}$  is a proper subset of $M$, argue as before, using the Whitney decomposition. If $O=M$, then $O$ is a ball itself, 
and we argue directly, without the Whitney decomposition, replacing 
$\mu\left(\widetilde{E_k}\right)$ by the left-hand side of (\ref{goodineqbounded}). \par
It remains to prove (\ref{claimconst}). First, if $t\geq 20\delta$, $B(x,t/20)=M$ and $V(y,t)=\mu(M)$ for any $y\in M$, so that
\[
S_{1/20}^{20\delta,R}f(x)^2=\mu(M)^{-1}\int_M\int_{20\delta}^R \left\vert t Du(y,t)\right\vert^2 dy\frac{dt}t.
\]
Next, computations similar to the estimates of $I$ above yield, for all $t>0$,
\[
\int_M \left\vert t Du(y,t)\right\vert^2 dy=-\frac t 4 \int_M \partial_t \left\vert u(y,t)\right\vert^2dy,
\]
so that, integrating by parts with respect to $t$, we obtain
\[
S_{1/20}^{20\delta,R}f(x)^2\leq \frac 1{4\mu(M)} \int_M \left\vert u(y,20\delta)\right\vert^2dy.
\]
But, for any $s\leq 20\delta$, the semigroup contraction property shows that
\[
 \int_M \left\vert u(y,20\delta)\right\vert^2dy \leq \int_M \left\vert u(y,s)\right\vert^2dy.
\]
It follows that
\[
S_{1/20}^{20\delta,R}f(x)^2\leq \frac 1{8\delta\mu(M)} \int_M\int_{\delta (1/c-1)}^{\delta (1/c+1)}  \left\vert u(y,s)\right\vert^2dyds.
\]
Noticing that $M=B(x,\delta)=B(x,c(\delta/c))$, one concludes, by definition of the maximal function, that
\[
S_{1/20}^{20\delta,R}f(x)^2\leq \frac 1{8c} \tilde f^{\ast}(x)^2,
\]
which is (\ref{claimconst}). The proof of Lemma \ref{goodlambda} is now complete. \hfill\fin
 \begin{rem} \label{sumofsquares}
The same proof shows that (\ref{Lusin}) holds if $\left\vert 
Du(y,t)\right\vert^{2}$ is replaced by the sum $\left\vert 
du(y,t)\right\vert^{2}+ \left\vert d^{\ast}u(y,t)\right\vert^{2}$ in 
the definition of $S_{\alpha}^{\varepsilon,R}$ (the sum is important). This is a stronger 
fact, since $\left\vert Du(y,t)\right\vert^{2}\leq 2\left\vert 
du(y,t)\right\vert^{2}+ 2\left\vert d^{\ast}u(y,t)\right\vert^{2}$ 
(but observe that, if one restricts to $k$-forms for fixed $0\leq 
k\leq \mbox{dim }M$, the two versions are equal). We could then 
conclude using the $H^1_{d}(\Lambda T^{\ast}M)$ and $H^1_{d^{\ast}}(\Lambda T^{\ast}M)$ spaces.
\end{rem}

\SE{Further examples and applications} \label{appl}
\subsection{The Coifman--Weiss Hardy space}
In this section, we focus on the case of functions, {\it i.e.} $0$-forms. 
Assuming that $M$ 
satisfies (\ref{D}), we may compare
$H^{1}_{d^{\ast}}(\Lambda^{0}T^{\ast}M)=H^{1}_D(\Lambda^{0}T^{\ast}M)=
H^{1}_\Delta(\Lambda^{0}T^{\ast}M)$ with the Coifman-Weiss Hardy 
space, 
{\it i.e.} the $H^{1}$ space defined in the general context of a space of 
homogeneous type in \cite{coifmanweiss}. 

We first recall what this space is. A (Coifman-Weiss) atom is a 
function 
$a\in L^{2}(M)$ supported in a ball $B\subset M$ and satisfying
\[
\int_{M}a(x)dx=0\mbox{ and } \left\Vert a\right\Vert_{2}\leq 
V(B)^{-1/2}.
\]
A complex-valued function $f$ on $M$ belongs to $H^{1}_{CW}(M)$ if 
and 
only if it can be written as
\[
f=\sum_{k\geq 1} \lambda_ka_k
\]
where $\sum_k\left\vert \lambda_k\right\vert<+\infty$ and the $a_k$'s 
are 
Coifman-Weiss atoms. Define
\[
\left\Vert f\right\Vert_{H^{1}_{CW}(M)}=\inf\sum_{k\geq 1}\left\vert 
\lambda_{k}\right\vert,
\]
where the infimum is taken over all such decompositions of $f$. 
Equipped with this norm, $H^{1}_{CW}(M)$ is a Banach space.
The link between this Coifman-Weiss space and the 
$H^{1}_{d^{\ast}}(\Lambda^0T^{\ast}M)$ space 
is as follows:
\begin{theo} \label{cw}
Assume  (\ref{D}). Then 
$H^{1}_{d^{\ast}}(\Lambda^0T^{\ast}M)\subset 
H^{1}_{CW}(M)$.
\end{theo}
\noindent{\bf Proof: }
Recall that a 
(Coifman-Weiss) molecule is a function $f\in L^{1}(M)\cap L^{2}(M)$ 
such that 
\[
\int_{M}f(x)dx=0
\]
and there exist $x_{0}\in M$ and $\varepsilon>0$ with
\begin{equation} \label{conditionmolec}
\left(\int_{M} \left\vert f(x)\right\vert^{2}dx\right)\left(\int_{M} 
\left\vert 
f(x)\right\vert^{2}m(x,x_{0})^{1+\varepsilon}dx\right)^{\frac 
1{\varepsilon}}\leq 1,
\end{equation}
where $m(x,x_{0})$ is the infimum of the measures of the balls both 
containing $x$ and $x_{0}$. It is shown in \cite{coifmanweiss} 
(Theorem C, p. 594) that such a 
molecule belongs to $H^{1}_{CW}(M)$ with a norm only depending on the 
constant in (\ref{D}) 
and $\varepsilon$. Note that condition (\ref{conditionmolec}) is 
satisfied if
\[
\left(\int_{M} \left\vert f(x)\right\vert^{2}dx\right)\left(\int_{M} 
\left\vert 
f(x)\right\vert^{2}V(x_{0},d(x,x_{0}))^{1+\varepsilon}dx\right)^{\frac 
1{\varepsilon}}\leq 1.
\]
Let $f\in {\mathcal R}(D)\cap H^{1}_{d^{\ast}}(\Lambda^0T^{\ast}M)$ and $F\in 
T^{1,2}(\Lambda^{0}T^{\ast}M)$ with 
$\left\Vert 
F\right\Vert_{T^{1,2}(\Lambda^{0}T^{\ast}M)} \sim \left\Vert 
f\right\Vert_{H^{1}_{d^{\ast}}(\Lambda^0T^{\ast}M)}$ and
\[
f=\int_{0}^{+\infty} 
(tD)^{N}(I+itD)^{-\alpha}F_{t}\frac{dt}t
\]
with $N> \D/2+1$ and $\alpha=N+2$. Since $F$ has an atomic 
decomposition in $T^{1,2}(\Lambda^{0}T^{\ast}M)$, it is enough to show that, whenever 
$A$ 
is a (scalar-valued) atom in $T^{1,2}(\Lambda^{0}T^{\ast}M)$ supported in $T(B)$ for 
some ball $B\subset M$,
\[
a=\int_{0}^{+\infty} 
(tD)^{N}(I+itD)^{-\alpha}A_{t}\frac{dt}t
\]
belongs to $H^{1}_{CW}(M)$ and satisfies
\begin{equation} \label{molec}
\left\Vert a\right\Vert_{H^{1}_{CW}(M)}\leq C.
\end{equation}
To that purpose, it suffices to check that, up to a multiplicative 
constant, $a$ is a Coifman-Weiss molecule. First, since $a\in 
{\mathcal R}(D)$ and $a$ is a function, it is clear 
that $a$ has zero integral. Furthermore, the spectral theorem  shows 
that
\[
\left\Vert a\right\Vert_2\leq CV^{-1/2}(B).
\]
Moreover, if $B=B(x_0,r)$ and if $\varepsilon>0$,
\[
\begin{array}{lll}
\displaystyle \int_M \left\vert a(x)\right\vert^2 
V^{1+\varepsilon}(x_0,d(x,x_0))dx & = & \displaystyle \int_{2B} 
\left\vert a(x)\right\vert^2 V^{1+\varepsilon}(x_0,d(x,x_0))dx\\
& + & \displaystyle \sum_{k\geq 1} \int_{2^{k+1}B\setminus 2^kB} 
\left\vert a(x)\right\vert^2 
V^{1+\varepsilon}(x_0,d(x,x_0))dx \\
&= & \displaystyle A_{0}+\sum_{k\geq 1} A_k.
\end{array}
\]
On the one hand, by the doubling property,
\[
A_{0} \leq CV^{1+\varepsilon}(B)\left\Vert a\right\Vert_2^2\leq 
CV(B)^{\varepsilon}.
\]
On the other hand, if $k\geq 1$, using Lemma \ref{offdiag1} and 
choosing $N^{\prime}$ such that $N^{\prime}>\D (1+\varepsilon)/2$, 
one has
\[
\begin{array}{lll}
A_{k}^{1/2}& \leq & \displaystyle V^{(1+\varepsilon)/2}(2^{k+1}B) \int_0^{r} \left\Vert 
(tD)^N 
(I+itD)^{-\alpha}A_{t}\right\Vert_{L^2(2^{k+1}B\setminus 
2^kB)}\frac{dt}t\\
& \leq & \displaystyle V^{(1+\varepsilon)/2}(2^{k+1}B) \int_0^r 
\left(\frac t{2^kr}\right)^{N^{\prime}} \left\Vert 
A_{t}\right\Vert_2 \frac{dt}t\\
& \leq & \displaystyle 2^{k(\D(1+\varepsilon)/2-N^{\prime})} 
V(B)^{\varepsilon/2}.
\end{array}
\]
Finally,
\[
\int_{M} \left\vert a(x)\right\vert^2 
V^{1+\varepsilon}(x_0,d(x,x_0))dx \leq CV^{\varepsilon}(B),
\]
which ends the proof of (\ref{molec}), therefore of Theorem 
\ref{cw}.\hfill\fin

We will focus on the converse inclusion in Theorem \ref{cw} in the 
following section. 

\subsection{Hardy spaces and Gaussian estimates}
In this section, we give further results about $H^p(\Lambda T^{\ast}M)$ spaces 
assuming some ``Gaussian'' upper bounds for the heat kernel of the 
Hodge-de Rham Laplacian on $M$. Denote by $n$ the dimension of $M$. For 
each $0\leq k\leq n$, let $p_{t}^k$ be the kernel of 
$e^{-t\Delta_{k}}$, where $\Delta_{k}$ is the Hodge-de Rham Laplacian 
restricted to $k$-forms. Say that $(G_{k})$ holds if there exist 
$C,c>0$ such that, for all $t>0$ and all $x,y\in M$,
\begin{equation} \label{Gk}
\left\vert p_{t}^k(x,y)\right\vert \leq \frac 
C{V(x,\sqrt{t})}e^{-cd^{2}(x,y)/t}.
\end{equation}
Say that $(G)$ holds if $(G_{k})$ holds for all $0\leq k\leq n$. See 
the introduction for comments on the validity of $(G_{(k)})$ when 
$k\geq 1$.

\subsubsection{The Coifman-Weiss Hardy space and Gaussian estimates}
Under the assumptions of Theorem \ref{cw}, and even if one assumes 
furthermore that $(G_0)$ holds, the inclusion $H^{1}_{d^{\ast}}(\Lambda^{0}T^{\ast}M)\subset H^{1}_{CW}(M)$, proved in Theorem \ref{cw}, is strict 
in general. This can be seen by considering the example where $M$ is 
the union of two copies of 
$\R^n$ ($n\geq 2$) glued smoothly together by a cylinder. First, on this 
manifold, (\ref{D}) and (\ref{Gk}) clearly hold (see \cite{couduo}). 
Moreover, Theorem \ref{rieszboundfour} asserts that the Riesz transform $d\Delta^{-1/2}$ is 
$H^1_{d^{\ast}}(\Lambda^0T^{\ast}M)-L^1(\Lambda^{1}T^{\ast}M)$ bounded on 
this 
manifold. But, as was kindly explained to us by A. Hassell, it is 
possible to prove, using arguments analogous to those contained in 
\cite{carcouhas}, that the 
Riesz transform is not $H^1_{CW}(M)-L^1(M)$ bounded (while it is 
shown 
in \cite{carcouhas} that the Riesz transform is $L^p(M)$-bounded for 
all $1<p<n$). 

However, under a stronger assumption on $M$, the spaces $H^1_{CW}(M)$ 
and 
$H^1_{d^{\ast}}(\Lambda^0T^{\ast}M)$ do coincide. Say that $M$ satisfies 
an $L^2$ Poincar\'e inequality on balls if there exists $C>0$ such 
that, for any ball $B\subset M$ and any function $f\in 
C^{\infty}(2B)$,
\begin{equation} \label{Poincare}
\int_B \left\vert f(x)-f_B\right\vert^2 dx \leq Cr^2 \int_{2B} 
\left\vert \nabla f(x)\right\vert^2 dx,
\end{equation}
where $f_B$ denotes the mean-value of $f$ on $B$ and $r$ the radius 
of $B$. Then we have:
\begin{theo} \label{identification}
Assume (\ref{D}) and  (\ref{Poincare}). Then 
$H^1_{d^{\ast}}(\Lambda^0T^{\ast}M)=H^1_{CW}(M)$.
\end{theo}
Indeed, as recalled in the introduction, these assumptions on $M$ 
imply 
that $p_t$ satisfies the estimates (\ref{G}), and these estimates, in 
turn, easily imply that any atom in $H^1_{CW}(M)$ belongs to 
$H^1_{max}(\Lambda^{0}T^{\ast}M)$ with a controlled norm. See, for 
instance, 
\cite{ar}. \hfill\fin

As a consequence of Theorem \ref{identification} and of Theorem 
\ref{rieszboundfour}, we recover the following result, already obtained in 
\cite{russ}:
\begin{cor} \label{rieszhardyfunctions}
Assume (\ref{D}) and (\ref{Poincare}). Then,  the Riesz transform on functions
 $d\Delta^{-1/2}$ is $H^{1}_{CW}(M)-L^{1}(M)$ bounded.
\end{cor}
 
Moreover, under the assumptions of Corollary 
\ref{rieszhardyfunctions}, some kind of $H^1$-boundedness result for 
the Riesz transform had been proved by M. Marias and the third author in 
\cite{mariasruss}. Namely,
if $u$ is a harmonic function on $M$ (in the sense that $\Delta u=0$ in $M$) with a growth at most linear 
(which means that $\left\vert u(x)\right\vert \leq C(1+d(x_0,x))$ for 
some $x_0\in M$), the 
operator $R_uf= du \cdot d\Delta^{-1/2}f$ is 
$H^1_{CW}(M)$-bounded (here and after in this section, $\cdot$ stands 
for the real scalar product on $1$-forms). Actually, we can also recover this result 
using the Hardy spaces defined in the present paper. Indeed, since, by 
Theorem \ref{rieszboundfour}, 
$d\Delta^{-1/2}$ is 
$H^{1}_{d^{\ast}}(\Lambda^{0}T^{\ast}M)-H^{1}_{d}(\Lambda^{1}T^{\ast}M)$ 
bounded, it suffices to prove that the map $g\mapsto du\cdot g$ is 
$H^{1}_{d}(\Lambda^{1}T^{\ast}M)-H^{1}_{d^{\ast}}(\Lambda^{0}T^{\ast}M)$ 
bounded. 

To that purpose, because of the decomposition into molecules for 
$H^{1}_{d}(\Lambda^{1}T^{\ast}M)$, one may assume that $g=a$ is a $1$-molecule for $d^{\ast}$ in $H^1_d(\Lambda^1T^{\ast}M)$, see Section \ref{firstinclusion}. Namely, one has $a=db$ where $b\in 
L^{2}(\Lambda^{0}T^{\ast}M)$ and there exists a ball $B$ and a 
sequence $(\chi_{k})_{k\geq 0}$ adapted to $B$ such that, for each 
$k\geq 0$,
\[
\left\Vert \chi_{k}a\right\Vert_{2}\leq 2^{-k}V^{-1/2}(2^kB)\mbox{ 
and }\left\Vert \chi_{k}b\right\Vert_{2}\leq r2^{-k}V^{-1/2}(2^kB).
\]
But it is plain to see that, up to a constant, $du\cdot a$ is a $1$-molecule for $d^{\ast}$ in 
$H^1_{d^{\ast}}(\Lambda^{0}T^{\ast}M)$. Indeed, since $du$ is bounded 
on $M$, one has, for each $k\geq 0$,
\[
\left\Vert \chi_{k}du\cdot a\right\Vert_{2}\leq C\left\Vert \chi_{k}a\right\Vert_{2}\leq C2^{-k}V^{-1/2}(2^kB).
\]
Moreover, since $\Delta u=0$ on $M$, one has
\[
du\cdot a=du\cdot db=-d^{\ast}(bdu)+b\Delta u=-d^{\ast}(bdu),
\]
and, for each $k\geq 0$,
\[
\left\Vert \chi_{k}bdu\right\Vert_{2}\leq C\left\Vert 
\chi_{k}b\right\Vert_{2}\leq Cr2^{-k}V^{-1/2}(2^kB).
\]
This ends the proof. \hfill\fin
\subsubsection{The decomposition into molecules and Gaussian estimates} \label{molecgauss}
We state here an improved version of Theorem \ref{decompo}, assuming furthermore some Gaussian upper estimates:
\begin{theo} \label{decompogauss}
Assume  (\ref{D}).
\begin{itemize}
\item[$(a)$]
If $(G)$ holds, then $H^1(\Lambda T^{\ast}M)=H^1_{mol,1}(\Lambda T^{\ast}M)$.
\item[$(b)$]
If $1\leq k\leq n$ and $(G_{k-1})$ holds, then $H^1_d(\Lambda^kT^{\ast}M)=H^1_{d,mol,1}(\Lambda^kT^{\ast}M)$.
\item[$(c)$]
If $0\leq k\leq n-1$ and $(G_{k+1})$ holds, then $H^1_{d^{\ast}}(\Lambda^kT^{\ast}M)=H^1_{d^{\ast},mol,1}(\Lambda^kT^{\ast}M)$.
\end{itemize}
\end{theo}
This theorem roughly says that, assuming Gaussian estimates, any section of $H^1(\Lambda T^{\ast}M)$ can be decomposed by means of $1$-molecules instead of $N$-molecules 
for $N>\frac{\D}2+1$. Observe that, in assertion $(c)$, if $M=\R^n$, the conclusion for $k=0$ is nothing 
but the usual atomic decomposition for functions in $H^{1}(\R^n)$. \par

\medskip

\noindent{\bf Proof: }We just give a sketch, which follows the same lines as the one of Theorem \ref{decompo}, focusing on assertion $(a)$. The inclusion $H^1(\Lambda T^{\ast}M)\subset H^1_{mol,1}(\Lambda T^{\ast}M)$ was proved in Section \ref{firstinclusion} and does not require Gaussian estimates. As for the converse inclusion, consider a molecule $f=Dg$ where $f$ and $g$ satisfy (\ref{molecestim}), and define $F(x,t)=tDe^{-t^2\Delta}f(x)$. We argue exactly as in Section \ref{secondinclusion} for $\eta_kF$. For $\eta^{\prime}_kF$, we use the fact that, if $0\leq l\leq k-2$, $y\in 2^{k+1}B\setminus 2^kB$ and $r<t<2^{k+1}r$, assumption $(G)$ yields
\[
\begin{array}{lll}
\displaystyle \left\vert tDe^{-t^2\Delta}tD(\chi_lg)(y)\right\vert & = & \displaystyle \left\vert t^2\Delta e^{-t^2\Delta}(\chi_lg)(y)\right\vert\\
& \leq & \displaystyle \frac C{V(y,t)}e^{-c\frac{2^{2k}r^2}{t^2}}\int_{2^{l+1}B} \left\vert \chi_lg(z)\right\vert dz\\
& \leq & \displaystyle \frac C{V(y,t)}r2^{-l}e^{-c\frac{2^{2k}r^2}{t^2}}.
\end{array}
\]
Using this estimate, one concludes for $\eta^{\prime}_kF$ in the same way as in Section \ref{secondinclusion}. The other terms in the proof of Theorem \ref{decompo} are dealt with in a similar way. This kind of argument can easily be transposed for assertions $(b)$ and $(c)$. \hfill\fin

An observation related to Theorem \ref{decompogauss} is that the elements of $H^1_{mol,1}(\Lambda T^{\ast}M)$ actually have an {\bf atomic} decomposition. More precisely, a section $a\in L^2(\Lambda T^{\ast}M)$ is called an atom if there exist a ball $B\subset M$ with radius $r$ and a section $b\in L^2(\Lambda T^{\ast}M)$ such that $b$ is supported in $B$, $a=Db$ and
\begin{equation} \label{atomestim}
\left\Vert a\right\Vert_{L^2(\Lambda T^{\ast}M)}\leq V^{-1/2}(B)\mbox{ and }\left\Vert b\right\Vert_{L^2(\Lambda T^{\ast}M)}\leq rV^{-1/2}(B).
\end{equation}
Say that a section $f$ of $\Lambda T^{\ast}M$ belongs to $H^1_{at}(\Lambda T^{\ast}M)$ if and only if there exist a sequence $(\lambda_j)_{j\geq 1}\in l^1$ and a sequence $(a_j)_{j\geq 1}$ of atoms such that $f=\sum_j \lambda_ja_j$, and equip $H^1_{at}(\Lambda T^{\ast}M)$ with the usual norm. We claim that $H^1_{mol,1}(\Lambda T^{\ast}M)=H^1_{at}(\Lambda T^{\ast}M)$. Indeed, an atom is clearly a $1$-molecule up to a multiplicative constant. Conversely, let $a=Db$ be a $1$-molecule, $B$ a ball and $(\chi_k)_{k\geq 0}$ a sequence of $C^{\infty}(M)$ functions adapted to $B$, such that (\ref{molecestim}) holds with $N=1$. Notice that, for some universal constant $C^{\prime}>0$, one has
\[
\left\Vert \chi_{2^{k+2}B\setminus 2^{k-1}B}b\right\Vert_{L^2(\Lambda T^{\ast}M)}\leq C^{\prime}r2^{-k}V^{-1/2}(2^{k+2}B)
\]
(this fact is a consequence of the support properties of the $\chi_j$'s). Define now $C^{\prime\prime}=\max(CC^{\prime},1)$ where $C>0$ is the constant in (\ref{cutoff}). For all $k\geq 0$, set 
\[
b_k=\frac{2^{k-1}}{C^{\prime\prime}}\chi_kb\mbox{ and }a_k=Db_k.
\]
It is obvious that $b_k$ is supported in $2^{k+2}B$ and that $\left\Vert 
b_k\right\Vert_{L^2(\Lambda T^{\ast}M)}\leq 2^{k+2}rV^{-1/2}(2^{k+2}B)$. Moreover, 
by (\ref{development}),
\begin{equation} \label{ak}
a_k=\frac{2^{k-1}}{C^{\prime\prime}}\left(\chi_kDb+d\chi_{k}\wedge 
b-d\chi_{k}\vee b\right)=\frac{2^{k-1}}{C^{\prime\prime}}\left(\chi_ka+d\chi_{k}\wedge 
b-d\chi_{k}\vee b\right),
\end{equation}
which implies $\left\Vert a_k\right\Vert_{L^2(\Lambda T^{\ast}M)}\leq \frac{2^{k-1}}{C^{\prime\prime}}\left(2^{-k}V^{-1/2}(2^{k+2}B)+
\frac C{2^kr}C^{\prime}r2^{-k}V^{-1/2}(2^{k+2}B)\right)\leq 
V^{-1/2}(2^{k+2}B)$. Thus, for each $k\geq 0$, $a_k$ is an atom. Moreover, 
since $\displaystyle \sum_{k}\chi_{k}=1$, one has
\[
\sum_{k} d\chi_{k}\wedge b=\sum_{k} d\chi_{k}\vee b=0,
\]
and (\ref{ak}) therefore yields 
\[
a=\sum_{k\geq 0} \frac{C^{\prime\prime}}{2^{k-1}}a_k,
\]
which shows that $a\in H^1_{at}(\Lambda T^{\ast}M)$. \par

One can similarly define $H^1_{d,at}(\Lambda T^{\ast}M)$, and  $H^1_{d,at}(\Lambda T^{\ast}M)=H^1_{d,mol,1}(\Lambda T^{\ast}M)$ holds by an analoguous argument. As a corollary of this fact and Theorem \ref{decompogauss}, we get that, when $M=\R^n$,  $H^1_d(\Lambda^k T^{\ast}\R^n)$ coincides with the ${\mathcal H}^1_d(\R^n,\Lambda^k)$ space introduced in \cite{loumcintosh}, as was claimed in the introduction.
\subsubsection{$H^p$ spaces and $L^p$ spaces} 
It turns out that, assuming that $(G)$ holds (which is the case of $M=\R^n$), one can compare precisely $H^p(\Lambda T^{\ast}M)$ and $L^p(\Lambda T^{\ast}M)$ for $1<p\leq 2$:
\begin{theo} \label{hplp}
Assume (\ref{D}). Let $1<p<2$.
\begin{itemize}
\item[$(a)$]
Assume $(G)$. Then,  $H^p(\Lambda T^{\ast}M)=\overline{{\mathcal R}(D)\cap L^p(\Lambda T^{\ast}M)}^{L^p(\Lambda T^{\ast}M)}$,  $H^p_d(\Lambda T^{\ast}M)=\overline{{\mathcal R}(d)\cap L^p(\Lambda T^{\ast}M)}^{L^p(\Lambda T^{\ast}M)}$ and $H^p_{d^{\ast}}(\Lambda T^{\ast}M)=\overline{{\mathcal R}(d^{\ast})\cap L^p(\Lambda T^{\ast}M)}^{L^p(\Lambda T^{\ast}M)}$
\item[$(b)$]
 Assume $(G_k)$ for some $0\leq k\leq n$. Then $H^p(\Lambda^kT^{\ast}M)=\overline{{\mathcal R}(D)\cap L^p(\Lambda^kT^{\ast}M)}^{L^p(\Lambda^kT^{\ast}M)}$, and the corresponding equalities for $H^p_d(\Lambda^kT^{\ast}M)$ and $H^p_{d^{\ast}}(\Lambda^kT^{\ast}M)$ also hold.
\end{itemize}
\end{theo}
\noindent{\bf Proof: }For assertion $(a)$, the inclusion $H^p(\Lambda T^{\ast}M)\subset
\overline{{\mathcal R}(D)\cap L^p(\Lambda T^{\ast}M)}^{L^p(\Lambda T^{\ast}M)}$ has already been proved (Corollary \ref{Lp}) and does not require assumption $(G)$. Conversely, it is enough to deal with $f\in {\mathcal R}(D)\cap L^p(\Lambda T^{\ast}M)$. Theorem 6 in \cite{unpublished} ensures that ${\mathcal Q}_{\psi}f\in T^{p,2}(\Lambda T^{\ast}M)$, where $\psi\in \Psi_{1,\beta+1}(\Sigma^{\theta}_0)$ with $\beta=\left[\frac{\D}2\right]+1$ (this is where we use Gaussian estimates). Now, for suitable $\widetilde{\psi}\in \Psi_{\beta,2}(\Sigma^{\theta}_0)$, one has $f={\mathcal S}_{\widetilde{\psi}}{\mathcal Q}_{\psi}f$ since $f\in {\mathcal R}(D)$, which shows that $f\in H^p(\Lambda T^{\ast}M)$. The other equalities, as well as assertion $(b)$, have a similar proof.\hfill\fin

\begin{rem} 
It is not known whether  equality $H^p(\Lambda^0T^{\ast}M)=\overline{{\mathcal R}(D)\cap L^p(\Lambda^0T^{\ast}M)}^{L^p(\Lambda^0T^{\ast}M)}$ holds for $1< p<2$ in general ({\it i.e.} without assuming Gaussian estimates for the heat kernel).
\end{rem}
\begin{rem}
What happens in Theorem \ref{hplp} for $p\geq 2$ ? We proved in Corollary \ref{Lp} that $\overline{{\mathcal R}(D)\cap L^p(\Lambda T^{\ast}M)}^{L^p(\Lambda T^{\ast}M)}\subset H^p(\Lambda T^{\ast}M)$ for $2\leq p<+\infty$. The converse inclusion cannot be true in general. Indeed, assume that $\overline{{\mathcal R}(D)\cap L^p(\Lambda T^{\ast}M)}^{L^p(\Lambda T^{\ast}M)}= H^p(\Lambda T^{\ast}M)$. Then, if $f\in L^p(M)\cap {\mathcal R}(D)$ is a function, one has  $f\in H^p_{d^{\ast}}(\Lambda^0T^{\ast}M)$ and Theorem \ref{rieszboundter} shows that $d\Delta^{-1/2}f\in H^p_d(\Lambda^0T^{\ast}M)$. Our assumption therefore implies $d\Delta^{-1/2}f\in L^p(M)$. In other words, the Riesz transform on functions $d\Delta^{-1/2}$ is $L^p$-bounded, which is false in general for $p>2$, even if Gaussian upper estimates for the heat kernel hold (\cite{couduo}), and even if the $L^2$ Poincar\'e inequality for balls is true (see \cite{couli}). 
\end{rem} 
As a corollary of Theorem \ref{hplp} and of Theorem \ref{rieszboundfour}, we obtain:
\begin{cor} \label{couduorecover}
Assume (\ref{D}) and $(G_0)$. Then, for all $1<p\leq 2$, the Riesz transform on functions  $d\Delta^{-1/2}$ is $L^p(\Lambda^0T^{\ast}M)-L^p(\Lambda^1T^{\ast}M)$ bounded.
\end{cor}
\noindent{\bf Proof: }It suffices to consider $f\in {\mathcal R}(d^{\ast})\cap L^p(\Lambda T^{\ast}M)$. Then, Theorem \ref{hplp} ensures that $f\in H^p_{d^\ast}(\Lambda^0T^{\ast}M)$, and Theorem \ref{rieszboundfour} yields that $d\Delta^{-1/2}f\in H^p_{d}(\Lambda^1T^{\ast}M)\subset L^p(\Lambda^1T^{\ast}M)$. \hfill\fin

Note that Corollary \ref{couduorecover} is not new and was originally proved in \cite{couduo}. Actually, in \cite{couduo}, the weak $(1,1)$ boundedness of $d\Delta^{-1/2}$ is established, and the $L^p$ boundedness for $1<p\leq 2$ follows at once by interpolation with the $L^2$ boundedness. Our approach by Hardy spaces does not allow us to recover the weak $(1,1)$ boundedness for $d\Delta^{-1/2}$.


\begin{thebibliography}{AAA}
\bibitem{auscher} Auscher, P., On necessary and sufficient conditions 
for $L^p$ Estimates of Riesz Transforms Associated to Elliptic 
Operators on $\R^n$ and Related Estimates, to appear in {\it Mem. A. M. S.}, vol. 
{\bf 186}, 871 (March 2007).
\bibitem{acdh} Auscher, P., Coulhon, T., Duong, X. T. and Hofmann, 
S., Riesz transforms 
on manifolds and heat kernel regularity, {\it Ann. Sci. Ecole Norm. 
Sup.} {\bf 37}, 6, 911-957 (2004).
\bibitem{unpublished} Auscher, P., Duong, X. T. and McIntosh, A., Boundedness of Banach space valued singular integral operators and applications to Hardy spaces, unpublished manuscript. 
\bibitem{ahlmt} Auscher, P., Hofmann, S., Lacey, M., McIntosh, A. and Tchamitchian, 
P., The solution of the Kato square root problem for
second order elliptic operators on $\R^n$, {\it Ann. Math.} {\bf 156}, 
633-654 (2002).
\bibitem{amr} Auscher, P., McIntosh, A. and Russ, E., Hardy spaces of differential 
forms and Riesz transforms on Riemannian manifolds, preprint. 
\bibitem{ar} Auscher, P. and Russ, E., Hardy spaces and divergence
operators on strongly Lipschitz domains
of $\R^n$, {\it J. Funct. Anal.} {\bf 201} (1), 148-184 (2003).
\bibitem{grenoble} Auscher, P. and Tchamitchian, P., Calcul fonctionnel 
pr\'ecis\'e pour des op\'erateurs elliptiques complexes en dimension 
un (et
applications \`a certaines \'equations elliptiques complexes en 
dimension
deux), {\it Ann. Inst. Fourier} {\bf 45}, 3, 721-778 (1995).
\bibitem{asterisque} Auscher, P. and Tchamitchian, P., {\it Square root
problem for divergence operators and related topics}, Ast\'erisque,
{\bf 249}, Soc. Math. France (1998).
\bibitem{akm} Axelsson, A., Keith, S. and McIntosh, A., Quadratic estimates 
and functional calculi of perturbed Dirac operators, {\it Inventiones 
Mathematicae} {\bf 163}, 455--497 (2006).
\bibitem{bakry} Bakry, D., Etude des transformations de Riesz dans les 
vari\'et\'es riemanniennes \`a courbure de Ricci minor\'ee, in {\it 
S\'eminaire de Probabilit\'es, XXI},  137--172, Lecture Notes in 
Math., {\bf 1247}, Springer, Berlin (1987). 
\bibitem{BC}  Bishop, R. and Crittenden, R., \textit{Geometry of 
manifolds},
Academic Press, N. York (1964).
\bibitem{carronsurvey} Carron, G., Formes harmoniques $L^{2}$ sur les 
vari\'et\'es non-compactes, {\it Rend. Mat. Appl. } {\bf 7} (21) 14, 
87-119 (2001). 
\bibitem{carcouhas} Carron, G., Coulhon, T. and Hassell, A., Riesz 
transform and $L^p$ cohomology 
for manifolds with Euclidean ends, {\it Duke Math. J.} {\bf 133},  no. 
1, 59--93 (2006).
\bibitem{coifman} Coifman, R., A real-variable characterization of 
$H^p$, {\it
Studia Math.} {\bf 51}, 269-274 (1974).
\bibitem{coifmanmeyerstein} Coifman, R., Meyer, Y. and Stein, E.M., Some 
new function spaces
and their
applications to harmonic analysis, {\it J. Funct. Anal.} {\bf 62}, 
304-335 (1985).
\bibitem{coifmanweissspringer} Coifman, R. and Weiss, G., {\it Analyse 
harmonique non commutative sur certains espaces homog\`enes}, Lect. 
Notes Math. {\bf 242}, Springer Verlag (1971).
\bibitem{coifmanweiss} Coifman, R. and Weiss, G., Extensions of Hardy
spaces and their use in analysis, {\it Bull. Amer. Math. Soc.} {\bf 
83}, 
569-645 (1977).
\bibitem{couduo} Coulhon, T. and Duong, X.T., Riesz transforms for $1\leq 
p\leq 2$, {\it Trans. Amer. Math. Soc.}  {\bf 351}, no. 3, 1151-1169 (1999).
\bibitem{couli} Coulhon, T. and Li, H. Q., Estimations inf\'erieures du noyau de la chaleur sur les vari\'et\'es coniques et transform\'ee de Riesz, 
{\it Archiv der Mathematik} {\bf 83}, 229-242 (2004).
\bibitem{CSC}  Coulhon, T. and Saloff-Coste, L., Vari\'{e}t\'{e}s 
Riemanniennes
isom\'{e}triques \`{a} l'infini, \textit{Rev. Mat. Iberoamericana} {\bf 
11}, 687-726 (1995).
\bibitem{coulhonzhang} Coulhon, T. and Zhang, Qi S., Large time behaviour 
of heat kernels on forms, to appear in {\it J. Diff. Geom.}, available at http://fr.arxiv.org/abs/math.DG/0506252.
\bibitem{davidsemmes} David, G., Journ\'e, J. L. and Semmes, S.,  
Op\'erateurs de Calder\'on-Zygmund, fonctions para-accr\'etives et 
interpolation, {\it Rev. Mat. Iberoamer.} {\bf 1}, 1-56 (1985).
\bibitem{davies1} Davies, E. B., Heat kernel bounds, conservation of 
probability and the Feller property, {\it J. Analyse Math.} {\bf 58}, 
99-119 (1992).
\bibitem{davies2} Davies, E. B., Uniformly elliptic operators with 
measurable 
coefficients, {\it J. Funct. Anal.}  {\bf 132},  no. 1, 141--169 
(1995).
\bibitem{derham} De Rham, G., {\it Vari\'et\'es diff\'erentiables, 
formes, 
courants, formes harmoniques}, 3\`eme \'edition, Hermann, Paris (1973).
\bibitem{feffermanstein} Fefferman, C. and Stein, E. M., $H^p$ spaces of
several variables, {\it Acta Math.} {\bf 129}, 137-195 (1972).
\bibitem{gaffney} Gaffney, M. P., The conservation property of the 
heat 
equation on Riemannian manifolds, {\it Comm. Pure Appl. Math.} {\bf 
12}, 
1-11 (1959).
\bibitem{ghl} Gilbert, J. E., Hogan, J. A. and Lakey, J. D., Atomic decomposition of divergence-free Hardy spaces, {\it Mathematica Moraviza}, Special Volume, Proc. 
IWAA, 33-52 (1997). 
\bibitem{grigoryan} Grigor'yan, A., Heat equation on a non-compact 
Riemannian manifold, {\it Math. USSR Sb.} {\bf 72} (1), 47-77 (1992).
\bibitem{HM} Hofmann, S. and Mayboroda, S., Hardy and BMO spaces associated to divergence form elliptic operators, in preparation.
\bibitem{latter} Latter, R. H., A characterization of $H^p(\R^n)$
in terms of atoms, {\it Studia Math.} {\bf 62}, 1, 93-101 (1978).
\bibitem{loumcintosh} Lou, Z. and McIntosh, A., Hardy spaces of exact 
forms 
on $\R^n$, {\it Trans. Amer. Math. Soc.} {\bf 357}, 4, 1469-1496 (2005).
\bibitem{mariasruss} Marias, M. and Russ, E., $H^1$-boundedness of Riesz 
transforms and imaginary powers of the Laplacian on Riemannian 
manifolds, {\it Ark. Mat.} {\bf 41}, 
115-132 (2003).
\bibitem{martellthesis} Martell, J.-M., Desigualdades con pesos en el 
An\'alisis 
de Fourier: de los espacios de tipo homog\'eneo a las medidas no 
doblantes, 
Ph. D., Universidad Aut\'onoma de Madrid (2001).
\bibitem{macfunct} McIntosh, A., Operators which have an $H^{\infty}$
functional
calculus, in {\it Miniconference on operator
theory and partial differential equations} (Canberra), Centre for
Math. and
Appl., vol. {\bf 14}, Australian National Univ., 210-231 (1986).
\bibitem{meyer} Meyer, Y., {\it Ondelettes et op\'erateurs}, t. II, 
Hermann (1990).
\bibitem{necas} Necas, J., {\it Les m\'ethodes directes en th\'eorie 
des \'equations elliptiques}, Masson et Cie, Eds.;  Paris, Academia, 
Editeurs, Prague (1967).
\bibitem{russ} Russ, E., $H^1-L^1$ boundedness of Riesz transforms on 
Riemannian manifolds and on graphs, {\it Pot. Anal.} {\bf 14}, 
301-330 (2001).
\bibitem{russtent} Russ, E., The atomic decomposition for tent spaces on spaces of homogeneous type, preprint.
\bibitem{salofflie} Saloff-Coste, L., Analyse sur les groupes de Lie 
\`a 
croissance polynomiale, {\it Ark. Mat.}  {\bf 28}, no. 2, 315--331 (1990).
\bibitem{saloff-coste}  Saloff-Coste, L., Parabolic Harnack inequality
for divergence form second order differential operators, {\it Pot. 
Anal.} {\bf 4}, 4, 429-467 (1995).
\bibitem{sc} Schwarz, G., {\it Hodge decomposition, a method for 
solving boundary value problems}, Lecture Notes in Mathematics {\bf 
1607}, 
Springer-Verlag, Berlin Heidelberg (1985).
 
\bibitem{semmes} Semmes, S., A primer on Hardy spaces and some remarks 
on a theorem of Evans and M\"uller, {\it Comm. Partial Differential 
Equations} {\bf 19}, 277-319 (1994).
\bibitem{stein70} Stein, E. M., {\it Singular integrals and
differentiability of functions}, Princeton University
Press (1970).
\bibitem{stein} Stein, E. M., {\it Harmonic analysis: Real-variable
methods, Orthogonality, and Oscillatory Integrals}, Princeton 
University
Press (1993).
\bibitem{steinweiss} Stein, E. M. and Weiss, G., On the theory of harmonic 
functions of several variables, I: the theory of $H^p$ spaces, {\it 
Acta Math.} {\bf 103}, 25-62 (1960). 
\bibitem{strihardy} Strichartz, R. S., The Hardy space $H^1$ on manifolds 
and submanifolds, {\it Canad. J. Math.} {\bf 24}, 915-925 (1972).
\bibitem{stri} Strichartz, R. S., Analysis of the Laplacian on 
the complete Riemannian manifold, {\it J. Funct. Anal.} {\bf 52}, 1, 
48-79 (1983).
\bibitem{wilson} Wilson, J., On the atomic decomposition for Hardy 
spaces, {\it Pacific J. Math.} {\bf 116}, 201-207 (1985).
\bibitem{zygmund} Zygmund, A., {\it Trigonometric series}, Cambridge 
Univ. Press (1959).
\end{thebibliography}
\end{document}